\documentclass[10pt]{amsart}
\usepackage[dvipdfmx]{graphicx}
\usepackage{amsmath,amscd,amsthm,amsxtra,yhmath,stackrel}
\usepackage{epsfig,graphics,color,colortbl}
\usepackage{amssymb,latexsym}
\usepackage{mathrsfs}
\usepackage{bm}
\usepackage[poly,all]{xy}
\usepackage{hyperref}
\usepackage{tikz-cd}
\usepackage[draft]{todonotes}
\usepackage{upgreek}
\usepackage{ytableau}
\usepackage{stmaryrd}
\usepackage{verbatim}
\usepackage{mathtools}
\usepackage[title]{appendix}
\usepackage{enumerate}
\usepackage[normalem]{ulem}
\usepackage{datetime}
\usepackage[pagewise]{lineno} 

\newtheorem{thm}{\bf Theorem}[section]
\newtheorem{df}[thm]{\bf Definition}
\newtheorem{prop}[thm]{\bf Proposition}
\newtheorem{cor}[thm]{\bf Corollary}
\newtheorem{lem}[thm]{\bf Lemma}
\newtheorem{rem}[thm]{\bf Remark}
\newtheorem{ex}[thm]{\bf Example}

\numberwithin{equation}{section}

\newcommand{\mc}{\mathcal}
\newcommand{\mf}{\mathfrak}
\newcommand{\ms}{\mathscr}

\newcommand{\pf}{\noindent{\bfseries Proof. }}
\newcommand{\ov}{\overline}

\newcommand{\U}{{\mc U}}

\newcommand{\cP}{\mathscr{P}}
\newcommand{\cO}{\mc{O}}
\newcommand{\I}{\mathbb{I}}

\newcommand{\ttq}{\texttt{q}}

\newcommand{\Z}{\mathbb{Z}}
\newcommand{\Q}{\mathbb{Q}}

\newcommand{\e}{\epsilon}
\newcommand{\de}{\delta}

\newcommand{\te}{\tilde{e}}
\newcommand{\tf}{\tilde{f}}
\newcommand{\td}{\widetilde}

\newcommand{\gl}{\mf{gl}}

\newcommand{\g}{\mf{g}}

\newcommand{\La}{\Lambda}

\newcommand{\la}{\lambda}

\newcommand{\red}[1]{{\color{red}#1}}

\newcommand{\ot}{\otimes}

\newcommand{\sa}{(-1)^{\e_a}}

\newcommand{\bq}{{\bf q}}
\newcommand{\bff}{{\bf f}}

\newcommand{\bfF}{{\bf F}}

\newcommand{\db}[1]{\llbracket #1 \rrbracket}

\begin{document}
\title
[Parabolic Verma modules of the quantum orthosymplectic superalgebra]{Crystal bases of parabolic Verma modules over the quantum orthosymplectic superalgebras}

\author{IL-SEUNG JANG}
\address{Department of Mathematics, Incheon National University, Incheon 22012, Korea}
\email{ilseungjang@inu.ac.kr}

\author{JAE-HOON KWON}
\address{Department of Mathematical Sciences and RIM, Seoul National University, Seoul 08826, Korea}
\email{jaehoonkw@snu.ac.kr}

\author{AKITO URUNO}
\address{Department of Mathematical Sciences, Seoul National University, Seoul 08826, Korea}
\email{aki926@snu.ac.kr}

\thanks{This work is supported by the National Research Foundation of Korea(NRF) grant funded by the Korea government(MSIT) (No.\,2019R1A2C1084833 and 2020R1A5A1016126).}

\begin{abstract}
We show that there exists a unique crystal base of a parabolic Verma module over a quantum orthosymplectic superalgebra, which is induced from a $q$-analogue of a polynomial representation of a general linear Lie superalgebra.
 
\end{abstract}

\maketitle
\setcounter{tocdepth}{1}

\noindent

\section{Introduction}

Let $\mf{g}$ be a classical Lie superalgebra. It is a natural question whether there is a crystal base of a highest weight representation of a quantized enveloping algebra $U_q(\mf{g})$ in the sense of Kashiwara \cite{Kas91}. It is shown that there exists a family of irreducible representations with crystal bases when $\g$ is a general linear Lie superalgebra $\mf{gl}(m|n)$ in \cite{BKK}, a queer Lie superalgebra $\mf{q}(n)$ in \cite{GJKKK}, and an orthosymplectic superalgebra $\mf{osp}(m|2n)$ with $m\ge 2$ in \cite{K15,K16}. 

Unlike symmetrizable Kac-Moody algebras, those families of representations do not seem to yield natural inductive systems to give a crystal base of $U_q(\g)^-$, the negative half of $U_q(\g)$. Instead, a crystal base of  $U_q(\g)^-$ for $\gl(m|n)$ is given in \cite{JKU} by taking a limit of the crystal bases of $q$-deformed Kac modules \cite{K14}. For $\mf{q}(n)$, only a candidate of the abstract crystal of a lowest weight Verma module is given in \cite{SS}.

We continue to follow the approach in \cite{JKU,K14} to have a better understanding of crystal base theory for orthosymplectic Lie superalgebras. As a first step, we construct in this paper a crystal base of a parabolic Verma module of $U_q(\g)$ for an orthosymplectic Lie superalgebra $\g$, which is induced from a $q$-analogue of an irreducible polynomial representation of $\mf{gl}(m|n)$.

For positive integers $m$ and $n$ with $m\ge 2$, let $\g$ be an orthosymplectic Lie superalgebra with a fundamental system $\Pi=\{\alpha_0,\dots,\alpha_{m+n-1}\}$ or a Dynkin diagram as follows:
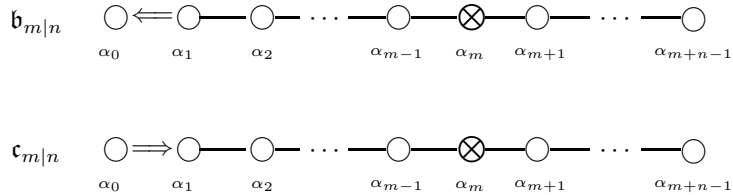
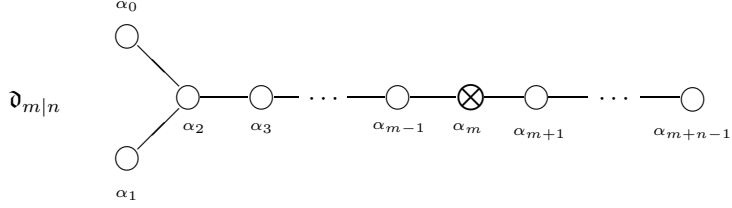
\begin{center}
\hskip -3cm \setlength{\unitlength}{0.16in}
\begin{picture}(24,4)
\put(5.6,2){\makebox(0,0)[c]{$\bigcirc$}}
\put(8,2){\makebox(0,0)[c]{$\bigcirc$}}
\put(10.4,2){\makebox(0,0)[c]{$\bigcirc$}}
\put(14.85,2){\makebox(0,0)[c]{$\bigcirc$}}
\put(17.25,2){\makebox(0,0)[c]{$\bigotimes$}}
\put(19.4,2){\makebox(0,0)[c]{$\bigcirc$}}
\put(24.5,1.95){\makebox(0,0)[c]{$\bigcirc$}}
\put(8.35,2){\line(1,0){1.5}}
\put(10.82,2){\line(1,0){0.8}}
\put(13.2,2){\line(1,0){1.2}}
\put(15.28,2){\line(1,0){1.45}}
\put(17.7,2){\line(1,0){1.25}}
\put(19.81,2){\line(1,0){1.28}}
\put(22.8,2){\line(1,0){1.28}}
%
\put(6.8,2){\makebox(0,0)[c]{$\Longleftarrow$}}
\put(12.5,1.95){\makebox(0,0)[c]{$\cdots$}}
\put(22,1.95){\makebox(0,0)[c]{$\cdots$}}
\put(3,1.8){\makebox(0,0)[c]{$\mf{b}_{m|n}$}}
\put(5.4,0.8){\makebox(0,0)[c]{\tiny $\alpha_{0}$}}
\put(7.8,0.8){\makebox(0,0)[c]{\tiny $\alpha_{1}$}}
\put(10.4,0.8){\makebox(0,0)[c]{\tiny $\alpha_{2}$}}
\put(14.8,0.8){\makebox(0,0)[c]{\tiny $\alpha_{m-1}$}}
\put(17.2,0.8){\makebox(0,0)[c]{\tiny $\alpha_{m}$}}
\put(19.5,0.8){\makebox(0,0)[c]{\tiny $\alpha_{m+1}$}}
\put(24.5,0.8){\makebox(0,0)[c]{\tiny $\alpha_{m+n-1}$}}
\end{picture}
\end{center}
\begin{center}
\hskip -3cm \setlength{\unitlength}{0.16in}
\begin{picture}(24,4)
\put(5.6,2){\makebox(0,0)[c]{$\bigcirc$}}
\put(8,2){\makebox(0,0)[c]{$\bigcirc$}}
\put(10.4,2){\makebox(0,0)[c]{$\bigcirc$}}
\put(14.85,2){\makebox(0,0)[c]{$\bigcirc$}}
\put(17.25,2){\makebox(0,0)[c]{$\bigotimes$}}
\put(19.4,2){\makebox(0,0)[c]{$\bigcirc$}}
\put(24.5,1.95){\makebox(0,0)[c]{$\bigcirc$}}
\put(8.35,2){\line(1,0){1.5}}
\put(10.82,2){\line(1,0){0.8}}
\put(13.2,2){\line(1,0){1.2}}
\put(15.28,2){\line(1,0){1.45}}
\put(17.7,2){\line(1,0){1.25}}
\put(19.81,2){\line(1,0){1.28}}
\put(22.8,2){\line(1,0){1.28}}
%
\put(6.8,2){\makebox(0,0)[c]{$\Longrightarrow$}}
\put(12.5,1.95){\makebox(0,0)[c]{$\cdots$}}
\put(22,1.95){\makebox(0,0)[c]{$\cdots$}}
\put(3,1.8){\makebox(0,0)[c]{$\mf{c}_{m|n}$}}
\put(5.4,0.8){\makebox(0,0)[c]{\tiny $\alpha_{0}$}}
\put(7.8,0.8){\makebox(0,0)[c]{\tiny $\alpha_{1}$}}
\put(10.4,0.8){\makebox(0,0)[c]{\tiny $\alpha_{2}$}}
\put(14.8,0.8){\makebox(0,0)[c]{\tiny $\alpha_{m-1}$}}
\put(17.2,0.8){\makebox(0,0)[c]{\tiny $\alpha_{m}$}}
\put(19.5,0.8){\makebox(0,0)[c]{\tiny $\alpha_{m+1}$}}
\put(24.5,0.8){\makebox(0,0)[c]{\tiny $\alpha_{m+n-1}$}}
\end{picture}
\end{center}
\begin{center}
\hskip -3cm \setlength{\unitlength}{0.16in}
\begin{picture}(24,5)
\put(6,0){\makebox(0,0)[c]{$\bigcirc$}}
\put(6,4){\makebox(0,0)[c]{$\bigcirc$}}
\put(8,2){\makebox(0,0)[c]{$\bigcirc$}}
\put(10.4,2){\makebox(0,0)[c]{$\bigcirc$}}
\put(14.85,2){\makebox(0,0)[c]{$\bigcirc$}}
\put(17.25,2){\makebox(0,0)[c]{$\bigotimes$}}
\put(19.4,2){\makebox(0,0)[c]{$\bigcirc$}}
\put(24.5,1.95){\makebox(0,0)[c]{$\bigcirc$}}
\put(6.35,0.3){\line(1,1){1.35}} 
\put(6.35,3.7){\line(1,-1){1.35}}
\put(8.4,2){\line(1,0){1.55}} 
\put(10.82,2){\line(1,0){0.8}}
\put(13.2,2){\line(1,0){1.2}} 
\put(15.28,2){\line(1,0){1.45}}
\put(17.7,2){\line(1,0){1.25}} 
\put(19.8,2){\line(1,0){1.25}}
\put(22.8,2){\line(1,0){1.28}}
\put(12.5,1.95){\makebox(0,0)[c]{$\cdots$}}
\put(22,1.95){\makebox(0,0)[c]{$\cdots$}}
\put(3,1.8){\makebox(0,0)[c]{$\mf{d}_{m|n}$}}
\put(6,5){\makebox(0,0)[c]{\tiny $\alpha_{0}$}}
\put(6,-1.2){\makebox(0,0)[c]{\tiny $\alpha_{1}$}}
\put(8.2,1){\makebox(0,0)[c]{\tiny $\alpha_{2}$}}
\put(10.4,1){\makebox(0,0)[c]{\tiny $\alpha_{3}$}}
\put(14.9,1){\makebox(0,0)[c]{\tiny $\alpha_{m-1}$}}
\put(17.15,1){\makebox(0,0)[c]{\tiny $\alpha_m$}}
\put(19.5,0.8){\makebox(0,0)[c]{\tiny $\alpha_{m+1}$}}
\put(24.5,0.8){\makebox(0,0)[c]{\tiny $\alpha_{m+n-1}$}}
\end{picture}\vskip 5mm
\end{center}
where $\bigotimes$ denotes an isotropic odd simple root.
%
We denote by $\mf{l}$ the subalgebra of $\g$ associated with $\Pi\setminus\{\alpha_0\}$, which is isomorphic to a general linear Lie superalgebra $\mf{gl}(m|n)$. 
Let $\U(\g)$ be a quantized enveloping algebra associated with $\mf{g}$ introduced in \cite{KOS}, which is an associative algebra over $\Bbbk=\Q(q)$ with $q$ indeterminate, and isomorphic to the one in  \cite{Ya94} under a mild extension, and let $\U(\mf{l})$ be the subalgebra corresponding to $\mf{l}$. 

For each $(m|n)$-hook partition $\la$, let $V_{\mf{l}}(\la)$ be an irreducible polynomial $\U(\mf{l})$-module which corresponds to a Specht module over a Hecke algebra of type $A$ corresponding to $\la$, in the super-analogue of Schur-Weyl-Jimbo duality \cite{J}. It is shown in \cite{BKK} that $V_{\mf{l}}(\la)$ has a unique crystal base, whose associated connected crystal is realized as the set of $(m|n)$-hook semistandard tableaux of shape $\la$ (cf.~\cite{BR}). 
We consider an induced representation 
\begin{equation*}
 P(\la)=\U(\g)\ot_{\U(\mf{p})}V_{\mf{l}}(\la),
\end{equation*}
where $\U(\mf{p})$ is the subalgebra of $\U(\g)$ associated to a parabolic subalgebra $\mf{p}$ generated by $\mf{l}$ and the Borel subalgebra of $\g$ with respect to $\Pi$. 
To construct a crystal base of $P(\la)$, we first regard $P(\la)$ as a polynomial $\U(\mf{l})$-module, which decomposes as a direct sum of $V_{\mf{l}}(\mu)$'s.
Then we have as a $\U(\mf{l})$-module
\begin{equation*}
 P(\la)\cong \mc{N}\ot_{\Bbbk} V_{\mf{l}}(\la),
\end{equation*}
where $\mc{N}$ is the subalgebra of $\U(\g)^-$ generated by the root  vectors corresponding to the roots of $\mf{g}$ but not of $\mf{p}$, and the action of $\U(\mf{l})$ on $\mc{N}$ is given by a quantum adjoint. 
Here we use a PBW type basis of the negative half $\U(\g)^-$ following \cite{CHW,Le}, since we do not have an analogue of braid symmetry \cite{Lu} for the odd isotropic simple root $\alpha_m$ on $\U(\g)$.
Note that
we have as a $\U(\mf{l})$-module 
\begin{equation*}\label{eq:decomp of crystal of N}
 \mc{N} \cong \bigoplus_{\mu\in \cP(\g)}V_{\mf l}(\mu),
\end{equation*}
where $\cP(\mf{b}_{m|n})$ is the set of $(m|n)$-hook partitions, $\cP(\mf{c}_{m|n})$ is the set of $(m|n)$-hook even partitions, and $\cP(\mf{d}_{m|n})$ is the set of $(m|n)$-hook partitions whose conjugates are even.

We show that the $A_0$-span of a PBW type basis of $\mc{N}$, after a suitable scaling on each monomial, together with its image at $q=0$ forms a crystal base of $\mc{N}$, where $A_0$ is the ring of regular functions $f(q)\in \Q(q)$ at $q=0$. We present an explicit form of Levendorskii–Soibelman type formula \cite{LS} for all the root vectors in $\mc{N}$, which are crucially used in the proof. 
Then a crystal base of $P(\la)$ as a $\U(\g)$-module is given by taking tensor product of the crystal bases of $\mc{N}$ and $V_{{\mf l}}(\la)$, where the crystal operators for $\alpha_0$ are naturally induced from the one on the negative half of the usual quantum group of type $X_m$ ($X=B,C,D$) in the Dynkin diagram of $\g$.
The crystal of $P(\la)$ also has a simple combinatorial description, which is connected and hence implies the uniqueness of the crystal base of $P(\la)$.

When the highest weight of $P(\la)$, say $\La_\la$, is replaced by $\La_\la+d\La_0$ for some non-negative integer $d$ (where $\La_0$ is the fundamental weight for $\alpha_0$), its maximal quotient is isomorphic to the irreducible highest weight $\U(\g)$-module studied in \cite{K15,K16}, which belongs to a semisimple category and has a crystal base.
We expect that the crystal base of $P(\la)$ is compatible with that of these irreducible $\U(\g)$-modules under the canonical projection, and also admits a natural limit to give a crystal base of $\U(\g)^-$. This will be discussed in the forthcoming paper.

The paper is organized as follows. In Section \ref{sec:QOSP}, we review necessary background on an orthosymplectic Lie superalgebra $\g$ and its quantized enveloping algebra $\U(\g)$. In Section \ref{sec:PBW type basis}, we recall a PBW type basis of $\U(\g)^-$ and its subalgebra $\mc{N}$. We also give explicit formulas for $q$-commutators of the root vectors in $\mc{N}$, and a $\U(\mf{l})$-module structure on $\mc{N}$ via quantum adjoint. In Section \ref{sec:main section}, we prove the existence of a crystal base of $P(\la)$, and the connectedness of its crystal.

\section{Quantum orthosymplectic superalgebra}\label{sec:quantum superalgebra}\label{sec:QOSP}

\subsection{Orthosymplectic Lie superalgebra}\label{subsec:notations}
Throughout the paper, we assume that $m$ and $n$ are positive integers with $m\ge 2$. Let us briefly recall necessary background on orthosymplectic Lie superalgebras.

Let $\I = \{\, 1, \dots ,m , m+1, \dots, m+n \,\}$ be given with a linear order $\{\, 1 <  \dots < m < m+ 1 < \dots < m+n \,\}$ and a $\mathbb{Z}_2$-grading $\I = \I_{\ov 0} \cup \I_{\ov 1}$ such that
$\I_{\ov 0} = \{\, 1, \dots, m \,\}$ and $\I_{\ov 1} = \{ \, m+1,\dots, m+n \,\}$. Let $\e=(\e_a)_{a\in \I}$ be a sequence with $\e_a=\varepsilon$ for $a\in \I_{\varepsilon}$ ($\varepsilon={\ov 0},{\ov 1}$). Let $P$ be a free abelian group with a basis $\{\,\delta_a\,|\,a\in \I\,\}$.

Let $I=\{\,0,1, \ldots,m,m+1,\dots, m+n-1\}$.
Let $\g=\mf{x}_{m|n}$ ($\mf{x}=\mf{b},\mf{c},\mf{d}$) denote an orthosymplectic Lie superalgebra associated to the following Dynkin diagram with the set of simple roots $\Pi=\{\,\alpha_i\,|\,i\in I\,\}\subset P$: \medskip

$\bullet$ $\mf{b}_{m|n}$
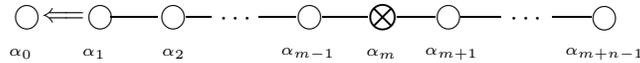
\begin{center}
\hskip -3cm \setlength{\unitlength}{0.16in}
\begin{picture}(24,4)
\put(5.6,2){\makebox(0,0)[c]{$\bigcirc$}}
\put(8,2){\makebox(0,0)[c]{$\bigcirc$}}
\put(10.4,2){\makebox(0,0)[c]{$\bigcirc$}}
\put(14.85,2){\makebox(0,0)[c]{$\bigcirc$}}
\put(17.25,2){\makebox(0,0)[c]{$\bigotimes$}}
\put(19.4,2){\makebox(0,0)[c]{$\bigcirc$}}
\put(24.5,1.95){\makebox(0,0)[c]{$\bigcirc$}}
\put(8.35,2){\line(1,0){1.5}}
\put(10.82,2){\line(1,0){0.8}}
\put(13.2,2){\line(1,0){1.2}}
\put(15.28,2){\line(1,0){1.45}}
\put(17.7,2){\line(1,0){1.25}}
\put(19.81,2){\line(1,0){1.28}}
\put(22.8,2){\line(1,0){1.28}}
%
\put(6.8,2){\makebox(0,0)[c]{$\Longleftarrow$}}
\put(12.5,1.95){\makebox(0,0)[c]{$\cdots$}}
\put(22,1.95){\makebox(0,0)[c]{$\cdots$}}
\put(5.4,0.8){\makebox(0,0)[c]{\tiny $\alpha_{0}$}}
\put(7.8,0.8){\makebox(0,0)[c]{\tiny $\alpha_{1}$}}
\put(10.4,0.8){\makebox(0,0)[c]{\tiny $\alpha_{2}$}}
\put(14.8,0.8){\makebox(0,0)[c]{\tiny $\alpha_{m-1}$}}
\put(17.2,0.8){\makebox(0,0)[c]{\tiny $\alpha_{m}$}}
\put(19.5,0.8){\makebox(0,0)[c]{\tiny $\alpha_{m+1}$}}
\put(24.5,0.8){\makebox(0,0)[c]{\tiny $\alpha_{m+n-1}$}}
\end{picture}
\end{center}
\begin{equation*}
\alpha_i=
\begin{cases}
-\de_{1} & \text{if $i=0$},\\
\de_{i}-\de_{i+1} & \text{if $1 \le i \le m+n-1$},
\end{cases}
\end{equation*}
\medskip

$\bullet$ $\mf{c}_{m|n}$
\begin{center}
\hskip -3cm \setlength{\unitlength}{0.16in}
\begin{picture}(24,4)
\put(5.6,2){\makebox(0,0)[c]{$\bigcirc$}}
\put(8,2){\makebox(0,0)[c]{$\bigcirc$}}
\put(10.4,2){\makebox(0,0)[c]{$\bigcirc$}}
\put(14.85,2){\makebox(0,0)[c]{$\bigcirc$}}
\put(17.25,2){\makebox(0,0)[c]{$\bigotimes$}}
\put(19.4,2){\makebox(0,0)[c]{$\bigcirc$}}
\put(24.5,1.95){\makebox(0,0)[c]{$\bigcirc$}}
\put(8.35,2){\line(1,0){1.5}}
\put(10.82,2){\line(1,0){0.8}}
\put(13.2,2){\line(1,0){1.2}}
\put(15.28,2){\line(1,0){1.45}}
\put(17.7,2){\line(1,0){1.25}}
\put(19.81,2){\line(1,0){1.28}}
\put(22.8,2){\line(1,0){1.28}}
%
\put(6.8,2){\makebox(0,0)[c]{$\Longrightarrow$}}
\put(12.5,1.95){\makebox(0,0)[c]{$\cdots$}}
\put(22,1.95){\makebox(0,0)[c]{$\cdots$}}
\put(5.4,0.8){\makebox(0,0)[c]{\tiny $\alpha_{0}$}}
\put(7.8,0.8){\makebox(0,0)[c]{\tiny $\alpha_{1}$}}
\put(10.4,0.8){\makebox(0,0)[c]{\tiny $\alpha_{2}$}}
\put(14.8,0.8){\makebox(0,0)[c]{\tiny $\alpha_{m-1}$}}
\put(17.2,0.8){\makebox(0,0)[c]{\tiny $\alpha_{m}$}}
\put(19.5,0.8){\makebox(0,0)[c]{\tiny $\alpha_{m+1}$}}
\put(24.5,0.8){\makebox(0,0)[c]{\tiny $\alpha_{m+n-1}$}}
\end{picture}
\end{center}
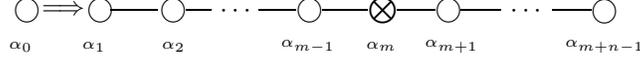
\begin{equation*}
\alpha_i=
\begin{cases}
-2\de_{1} & \text{if $i=0$},\\
\delta_i-\delta_{i+1} & \text{if $1\le i \le m+n-1$},
\end{cases}
\end{equation*}
\medskip

$\bullet$ $\mf{d}_{m|n}$
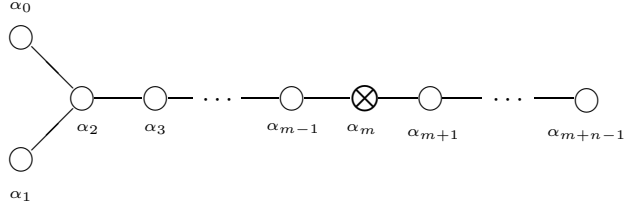
\begin{center}
\hskip -3cm \setlength{\unitlength}{0.16in} \medskip
\begin{picture}(24,5.8)
\put(6,0){\makebox(0,0)[c]{$\bigcirc$}}
\put(6,4){\makebox(0,0)[c]{$\bigcirc$}}
\put(8,2){\makebox(0,0)[c]{$\bigcirc$}}
\put(10.4,2){\makebox(0,0)[c]{$\bigcirc$}}
\put(14.85,2){\makebox(0,0)[c]{$\bigcirc$}}
\put(17.25,2){\makebox(0,0)[c]{$\bigotimes$}}
\put(19.4,2){\makebox(0,0)[c]{$\bigcirc$}}
\put(24.5,1.95){\makebox(0,0)[c]{$\bigcirc$}}
\put(6.35,0.3){\line(1,1){1.35}} 
\put(6.35,3.7){\line(1,-1){1.35}}
\put(8.4,2){\line(1,0){1.55}} 
\put(10.82,2){\line(1,0){0.8}}
\put(13.2,2){\line(1,0){1.2}} 
\put(15.28,2){\line(1,0){1.45}}
\put(17.7,2){\line(1,0){1.25}} 
\put(19.8,2){\line(1,0){1.25}}
\put(22.8,2){\line(1,0){1.28}}
\put(12.5,1.95){\makebox(0,0)[c]{$\cdots$}}
\put(22,1.95){\makebox(0,0)[c]{$\cdots$}}
\put(6,5){\makebox(0,0)[c]{\tiny $\alpha_{0}$}}
\put(6,-1.2){\makebox(0,0)[c]{\tiny $\alpha_{1}$}}
\put(8.2,1){\makebox(0,0)[c]{\tiny $\alpha_{2}$}}
\put(10.4,1){\makebox(0,0)[c]{\tiny $\alpha_{3}$}}
\put(14.9,1){\makebox(0,0)[c]{\tiny $\alpha_{m-1}$}}
\put(17.15,1){\makebox(0,0)[c]{\tiny $\alpha_m$}}
\put(19.5,0.8){\makebox(0,0)[c]{\tiny $\alpha_{m+1}$}}
\put(24.5,0.8){\makebox(0,0)[c]{\tiny $\alpha_{m+n-1}$}}
\end{picture}\vskip 8mm
\end{center}
\begin{equation*}
\alpha_i=
\begin{cases}
-\de_{1}-\de_{2} & \text{if $i=0$},\\
\delta_i-\delta_{i+1} & \text{if $1\le i \le m+n-1$}.
\end{cases}
\end{equation*}
\medskip

Here we assume that there is a symmetric bilinear form  $(\, \cdot\, |\, \cdot\, )$ on $P$ such that 
\begin{equation*}\label{eq:bilinear form on P}
 (\de_a|\de_b)=(-1)^{\e_a}r_{\g}\de_{ab} \quad (a,b\in \I),
\end{equation*} 
where 
\begin{equation*}
 r_{\mf{g}}=
\begin{cases}
 2 & \text{if $\mf{g}=\mf{b}_{m|n}$},\\
 1 & \text{if $\mf{g}=\mf{c}_{m|n},\mf{d}_{m|n}$}.
\end{cases}
\end{equation*}
We have 
\begin{equation*}\label{eq:norm}
 (\alpha_i|\alpha_i)\in \{0,\pm 2, \pm 4\}\quad  (i\in I).
\end{equation*}
In the above diagrams, $\bigotimes$ means that $\alpha_m$ is   isotropic, that is, $(\alpha_m|\alpha_m)=0$. We refer the reader to \cite{CW} for more detailed exposition of $\mf{g}$.

Let $\Phi^+$ be the set of reduced positive roots $\beta$ of $\g$, that is, $\beta$ is a positive root of $\g$ but $\beta/2$ is not, and let $\Phi^+_{\ov{0}}$ (resp.~$\Phi^+_{\ov{1}}$) be the set of even (resp.~odd) roots in $\Phi^+$. We have $\Phi^+_{\ov{1}}=\Phi^+_{\rm iso}\cup\Phi^+_{\text{n-iso}}$, where $\Phi^+_{\rm iso}$ is the set of  isotropic roots in $\Phi^+_{\ov 1}$ and $\Phi^+_{\text{n-iso}}=\Phi^+_{\ov 1}\setminus \Phi^+_{\rm iso}$.

Let us give an explicit description of $\Phi^+$.
Let $\mf{l}$ be the subalgebra of $\g$ corresponding to $\Pi\setminus\{\alpha_0\}$, which is isomorphic to the general linear Lie superalgebra $\gl(m|n)$.  
Let $\Phi^+(\mf{l})$ be the set of positive roots of $\mf{l}$. 
Then $\Phi^+(\mf{l})= \Phi^+(\mf{l})_{\ov{0}}\cup\Phi^+(\mf{l})_{\ov{1}}$, where
$\Phi^+(\mf{l})_{\varepsilon}=\Phi^+(\mf{l})\cap \Phi^+_{\varepsilon}$ ($\varepsilon=\ov{0}, \ov{1}$) given by
\begin{equation*}
\begin{split}
\Phi^+(\mf{l})_{\ov{0}}
&=\{\,\de_{i} - \de_{j}\,|\, 1\le i<j \le m \,\}\cup \{\,\de_{s} - \de_t\,|\, m+1\le s<t\le m+n \,\},\\
\Phi^+(\mf{l})_{\ov{1}}
&=\{\,\de_{i} - \de_s\,|\, 1\le i\le m,\ m+1\le s\le m+n  \,\}.
\end{split}
\end{equation*}

Let $\mf{u}$ be the subalgebra of $\g$ generated by the root vectors corresponding to $\Phi^+\setminus \Phi^+(\mf{l})$, and put $\Phi^+(\mf{u})= \Phi^+\setminus \Phi^+(\mf{l})$. 
Then $\Phi^+(\mf{u})=\Phi^+(\mf{u})_{\ov{0}}\cup  \Phi^+(\mf{u})_{\ov{1}}$, where 
$\Phi^+(\mf{u})_{\varepsilon}=\Phi^+(\mf{u})\cap \Phi^+_{\varepsilon}$ ($\varepsilon=\ov{0}, \ov{1}$) is given by
\begin{equation}\label{eq:roots in radical}
\begin{split}
 \Phi^+(\mf{u})_{\ov{0}} &=
 \begin{cases}
 \left\{\, -\de_{k},\ -\de_{i}-\de_{j},\ -\de_{s}-\de_{t}\,\right\} & \text{if $\mf{g}=\mf{b}_{m|n}$},\\
\left\{\, -2\de_{k},\ -\de_{i}-\de_{j},\ -\de_{s}-\de_{t}\,\right\} & \text{if $\mf{g}=\mf{c}_{m|n}$}, \\
\left\{\, -\de_{i}-\de_{j},\ -\de_{s}-\de_{t},\ -2\de_{u}\,\right\} & \text{if $\mf{g}=\mf{d}_{m|n}$},
\end{cases}\\
 \Phi^+(\mf{u})_{\ov{1}} &=
 \begin{cases}
\{\, -\de_{k}-\de_{u},\ -\de_{u} \,\} & \text{if $\mf{g}=\mf{b}_{m|n}$},\\
\{\, -\de_{k}-\de_{u} \,\}  & \text{if $\mf{g}=\mf{c}_{m|n}$, $\mf{d}_{m|n}$},
\end{cases}
 \end{split}
\end{equation}
where $1\le k\le m,\, 1\le i<j\le m$, $m+1\le s<t\le m+n$, and $m+1\le u \le m+n$. \smallskip




\subsection{Quantized enveloping algebra for $\mf{x}_{m|n}$}
Let us consider a quantized enveloping algebra for $\g=\mf{x}_{m|n}$.
We assume the following notations: 
\begin{itemize}
\item[$\bullet$] $\Bbbk=\Q(q)$ and $\Bbbk^{\times}=\Bbbk \setminus \left\{0\right\}$ where $q$ is an indeterminate,

\item[$\bullet$] ${\ttq}_a=\sa q^{\sa r_{\g}}$ $(a\in \I)$, that is,
\begin{equation*}
{\ttq}_a=
\begin{cases}
q^{r_{\mf{g}}} & \text{if $\e_a=0$}\\
-q^{-r_{\mf{g}}} & \text{if $\e_a=1$}\\
\end{cases} \quad (a\in \I),
\end{equation*}

\item[$\bullet$] ${\bq}(\,\cdot\,,\,\cdot\,)$: a symmetric biadditive function from $P\times P$ to $\Bbbk^{\times}$ given by
\begin{equation*}
\bq(\mu,\nu) = \prod_{a\in \I}{\tt q}_a^{\mu_a\nu_a},
\end{equation*}
for $\mu=\sum_{a\in\I}\mu_a\de_a$ and $\nu=\sum_{a\in\I}\nu_a\de_a$,

\item[$\bullet$] $q_i = q^{d_i}$ $(i\in I\setminus\{m\})$, where $d_i=|(\alpha_i|\alpha_i)|/2$, and $q_m = q^{r_{\mf g}}$.
\end{itemize}
\medskip

For $s\in\Z_{\ge 0}$ and $z\in q^{\Z_{>0}}$, we put
\begin{gather*}
[s]_z=\frac{z^s-z^{-s}}{z-z^{-1}},\quad
[s]_z!=[s]_z [s-1]_z \cdots [1]_z \quad (s\geq 1),\quad [0]_z!=1,
\end{gather*}
\noindent where we simply write $[s]$ when $z=q$, and $[s]_i=[s]_{q_i}$ when $z=q_i$ ($i\in I$). 
We also put
\begin{gather*}
\{ s \}_z = \frac{(-z)^s -z^{-s}}{-z-z^{-1}},\quad \{s\}_z!=\{s\}_z\{s-1\}_z\dots \{1\}_z\quad (s\geq 1),\quad  \{0\}_z!=1,
\end{gather*}
and write $\{s\}=\{s\}_q$ and $\{s\}_i=\{s\}_{q_i}$ ($i\in I$).

\begin{df}\label{def:U(e)}
{\rm
We define ${\U}(\g)$ to be an associative $\Bbbk$-algebra with $1$ 
generated by $k_\mu, e_i, f_i$ for $\mu\in P$ and $i\in I$ 
satisfying
{\allowdisplaybreaks
\begin{gather*}
k_{\mu}=1 \quad (\mu=0), \quad k_{\nu +\nu'}=k_{\nu}k_{\nu'} \quad (\nu, \nu' \in P),\label{eq:Weyl-rel-1} \\ 
k_\mu e_i k_{-\mu}=\bq(\mu,\alpha_i)e_i,\quad 
k_\mu f_i k_{-\mu}=\bq(\mu,\alpha_i)^{-1}f_i\quad (i\in I, \mu\in P), \label{eq:Weyl-rel-2} \\ 
e_if_j - f_je_i =\delta_{ij}\frac{k_{i} - k^{-1}_{i}}{q_i-q_i^{-1}}\quad (i,j\in I),\label{eq:Weyl-rel-3}\\
e_m^2= f_m^2 =0,\label{eq:Weyl-rel-4}
\end{gather*}
where $k_i=k_{\alpha_i}$ for $i\in I$, and 
\begin{gather*}
  x_i x_j -  x_j x_i =0
 \quad \text{($i,j \in I$ and $(\alpha_i|\alpha_j) = 0$)},\\
x_i^2 x_j- (-1)^{\e_i}[2]_i x_i x_j x_i + x_j x_i^2= 0
\quad \text{($i\in I\setminus\{0, m\}$, $j\in I\setminus\{0\}$, $i\neq j$ and $(\alpha_i|\alpha_j)\neq 0$)},    \\
\begin{array}{ll}
		  x_0^3 x_1- [3]_0 x_0^2 x_1 x_0 + [3]_0 x_0 x_1 x_0^2 - x_1 x_0^3 = 0\\ 
		  x_1^2 x_0- [2]_1 x_1 x_0 x_1 + x_0 x_1^2= 0
\end{array}
\quad (\g = \mf{b}_{m|n}),    \\
\begin{array}{ll}
		x_0^2 x_1- [2]_0 x_0 x_1 x_0 + x_1 x_0^2= 0 \\
		x_1^3 x_0- [3]_1 x_1^2 x_0 x_1 + [3]_1 x_1 x_0 x_1^2 - x_0 x_1^3 = 0
\end{array}
\quad (\g = \mf{c}_{m|n}), \\
\begin{array}{ll}
		x_0^2 x_2- [2]_0 x_0 x_2 x_0 + x_2 x_0^2= 0 \\
		x_2^2 x_0- [2]_2 x_2 x_0 x_2 + x_0 x_2^2= 0
\end{array}
	\quad (\g = \mf{d}_{m|n}\ \text{with $m\ge 3$}),\\
\begin{array}{ll}
  x_{m}x_{m'}x_{m}x_{m+1}  
- x_{m}x_{m+1}x_{m}x_{m'} 
+ x_{m+1}x_{m}x_{m'}x_{m} \\
\qquad \qquad \qquad \qquad \qquad - x_{m'}x_{m}x_{m+1}x_{m} 
+ [2]_{m} x_{m}x_{m'}x_{m+1}x_{m} =0 \\ 
\hskip 6cm (m'=0,1 \text{ if $\g={\mf d}_{2|n}$}, \text{ and } m'=m-1 \text{ otherwise} ),
\end{array}
\end{gather*}}
where $x=e,f$.}
\end{df}

There is a Hopf algebra structure on $\U(\g)$, where the comultiplication $\Delta$ and the antipode $S$ are given by 
\begin{equation}\label{eq:comult-1}
\begin{split}
& \Delta(k_\mu)=k_\mu\otimes k_\mu, \\ 
& \Delta(e_i)= 1\ot e_i + e_i\ot k_i^{-1}, \\
& \Delta(f_i)= f_i\ot 1 + k_i\ot f_i , \\  
S(k_i)= & k_i^{-1}, \quad S(e_i)=-e_i k_i, \quad S(f_i)=-k_i^{-1} f_i,
\end{split}
\end{equation}
for $\mu\in P$ and  $i\in I$.
%
%

Let $\U(\g)^+$ (resp. $\U(\g)^-$) be the subalgebra of $\U(\g)$ generated by $e_i$ (resp. $f_i$) for $i\in I$, and let $\U(\g)^0$ be the one generated by $k_\mu$ for $\mu\in P$.
We have $\U(\g)\cong \U(\g)^-\ot\,\U(\g)^0\ot\,\U(\g)^+$ as a $\Bbbk$-space, which follows from \cite[Theorem 10.5.1]{Ya94} and \eqref{eq:iso tau standard} below.
Let $\U(\g)^{\ge 0}$ be the subalgebra generated by $k_\mu$ and $e_i$ for $\mu\in P$ and $i\in I$.
We also have $\U(\g)^\pm =\bigoplus_{\alpha\in Q^\pm}\U(\g)^{\pm}_\alpha$, where $Q^\pm=\pm\sum_{i\in I}\Z_{\ge 0}\alpha_i$ and
$\U(\g)^{\pm}_\alpha=\{\,u\,|\,k_\mu u k_{-\mu}=\bq(\mu,\alpha)u\ (\mu\in P)\,\}$.
We put $|x|=\alpha$ for $x\in \U(\g)^\pm_\alpha$.  
For homogeneous $x,y\in \U(\g)^-$, we let
\begin{equation*} \label{eq: braket}
 [x,y]_{\bq} = xy - \bq(|x|,|y|)^{-1}yx.
\end{equation*}

For $i\in I$, let $e'_i$ and $e''_i$ denote the $\Bbbk$-linear maps on $\U(\g)^-$ defined by 
$e'_i(f_j)=e''_i(f_j)=\de_{ij}$ for $j\in I$ and 
\begin{equation} \label{eq: derivation}
\begin{split}
e'_i(uv)&=e'_i(u)v + \bq{(\alpha_i,|u|)}ue'_i(v),\\
e''_i(uv)&=e''_i(u)v + \bq{(\alpha_i,|u|)}^{-1}ue''_i(v),
\end{split}
\end{equation}
for $u,v\in \U(\g)^-$ with $u$ homogeneous.

For a $\U(\g)$-module $V$ and $\mu\in P$, let 
\begin{equation*}
V_\mu 
= \{\,u\in V\,|\,k_{\nu} u= \bq(\mu,\nu) u \ \ (\nu\in P) \,\}.
\end{equation*}
For $u\in V_\mu\setminus\{0\}$, we call $u$ a weight vector with weight $\mu$ and put ${\rm wt}(u)=\mu$. 
%

The algebra $\U(\g)$ was introduced in \cite{KOS} called a {\em generalized quantum group}. Let $U_q(\g)$ be the {\em quantum superalgebra} associated to $\g=\mf{x}_{m|n}$ \cite{Ya94} generated by $K_\mu, E_i, F_i$ for $\mu\in P$ and $i\in I$. 
As in the case of type $A$ \cite[Proposition 4.4]{KO}, these two algebras are isomorphic under mild extensions as follows:
Let $\Sigma$ be the bialgebra over $\Bbbk$ generated by $\sigma_a$ $(a\in \I)$  such that $\sigma_a\sigma_b=\sigma_b\sigma_a$, $\sigma_a^2=1$, and $\Delta(\sigma_a)=\sigma_a\otimes \sigma_a$ ($a,b\in \I$). 
Let $\U(\g)[\sigma]$ be the semidirect product of $\U(\g)$ and $\Sigma$, where $\Sigma$ acts on $\U(\g)$ by 
\begin{equation*}\label{eq:sigma action}
\begin{split}
&\sigma_a k_\mu =k_\mu,\quad
\sigma_ae_i=(-1)^{\e_a(\delta_a|\alpha_i)}e_i,\quad 
\sigma_af_i=(-1)^{\e_a(\delta_a|\alpha_i)}f_i,
\end{split}
\end{equation*}
for $a\in \I$, $\mu\in P$ and $i\in I$. 
%
Let $U_q(\g)[\sigma]$ be defined in the same way. Then there exists an isomorphism of $\Bbbk$-algebras 
\begin{equation}\label{eq:iso tau standard}
\xymatrixcolsep{2pc}\xymatrixrowsep{3pc}\xymatrix{
 U_q(\g)[\sigma] \ \ar@{->}^{\tau}[r] &\ \U(\g)[\sigma]},
\end{equation}
such that 
{\allowdisplaybreaks
\begin{equation*}\label{eq:iso phi}
\begin{split}
& \tau(\sigma_a) =\sigma_a \quad (1\le a\leq m+n),\qquad \tau(K_{\de_a})= 
\begin{cases}
k_{\de_a} & \text{$(1\leq a\leq m)$},\\ 
k_{\de_a}\sigma_a  & \text{$(m+1\le j\leq m+n)$},
\end{cases}\\
& \tau(E_i)= 
\begin{cases}
e_i  & \text{$(0\leq i\leq m)$},\\
-e_i(\sigma_i\sigma_{i+1})^{i+m}  & \text{$(m+1\le i\leq m+n-1)$},
\end{cases}\\
& \tau(F_i)= 
\begin{cases}
f_i  & \text{$(0\leq i\le m-1)$},\\
f_m\sigma_{m+1}  & \text{$(i=m)$},\\
f_i(\sigma_i\sigma_{i+1})^{i+m+1}  & \text{$(m+1\le i\leq m+n-1)$}.
\end{cases}
\end{split}
\end{equation*}}

\begin{rem}\label{rem:comparison of reps of two algs}
{\rm
Let $V$ be a $\U(\g)$-module with $V=\bigoplus_{\mu\in P}V_\mu$.
Then $V$ can be extended to a $\U(\g)[\sigma]$-module by letting
$\sigma_a u = (-1)^{\e_a\mu_a} u$
for $a\in \I$ and $u\in V_\mu$ with $\mu=\sum_{a\in \I}\mu_a\delta_a$. By \eqref{eq:iso tau standard}, we naturally have a $U_q(\g)$-module. For this reason, it does not make much difference when we consider $\U(\g)$ and its representations instead of $U_q(\g)$, while we remark the Hopf algebra structure on $\U(\g)[\sigma]$ is slightly different from that of $U_q(\g)[\sigma]$ (cf.~\cite[Remark 4.7]{KO}).
} 
\end{rem}

\section{PBW type basis} \label{sec:PBW type basis}
\subsection{Root vectors} \label{subsec:root vectors}

Let us construct a root vector $\bff_\beta$ ($\beta\in \Phi^+$) and a PBW type basis of $\U(\g)^-$ following \cite{CHW,Le}.

Let $\mc{A}=\{\,w_1<w_2<\cdots \,\}$ be a linearly ordered set of alphabets, and let $\mc{W}$ be the set of words with alphabets in $\mc{A}$. We put $w[i_1,\dots,i_k]=w_{i_1}\dots w_{i_k}$ for $i_1,\dots,i_k\in \Z_{>0}$.
We regard $\mc{W}$ as an ordered set with respect to the following lexicographic order:
\begin{equation*}
	w[i_1, \dots, i_a] < w[j_1, \dots, j_b]
\end{equation*}
if there exists an $r$ such that $w_{i_r} < w_{j_r}$ and $w_{i_s} = w_{j_s}$ for $s < r$, or if $a < b$ and $w_{i_s} = w_{j_s}$ for $1 \le s \le a$.
For a word $w=w[i_1,\dots,i_k] \in \mc{W}$, we say that $w$ is a {\it Lyndon word} if $w<w[i_{s},i_{s+1}\dots,i_k]$ for $2\le s\le k$. 
%

From now on, we assume that $\mc{A}=\{\,w_0<w_1<\dots<w_{m+n-1}\,\}$.
Consider the following set of Lyndon words $\mc{GL}(\mf{g})$ as follows:
{\allowdisplaybreaks
\begin{equation}\label{eq:good Lyndon words}
\begin{split}
 \mc{GL}(\mf{b}_{m|n})
 =&\{\,w[i,\dots,j]\,|\,0\le i\le j\le m+n-1 \,\}\\ 
  &\cup \{\,w[0,\dots,j,0,\dots,k]\,|\,0\le j<k \le m+n-1\,\},\\
 \mc{GL}(\mf{c}_{m|n})
 =&\{\,w[i,\dots, j] \,|\,0\le i\le j\le m+n-1 \,\}\\
 & \cup \{\,w[0,\dots,k,1,\dots,j]\,|\,1\le j<k\le m+n-1\,\}\\
 & \cup \{\,w[0,\dots,k,1,\dots,k]\,|\,1\le k\le m-1\,\},\\
 \mc{GL}(\mf{d}_{m|n}) 
 =&\{\,w[0]\,\}\cup \{\,w[0,2,\dots,i]\,|\,2\le i \le m+n-1\,\}\\
 & \cup \{\,w[i,\dots,j]\,|\,1\le i\le j\le m+n-1 \,\}\\ 
 & \cup \{\,w[0,2,\dots,k,1,\dots,j]\,|\,1\le j<k\le m+n-1\,\}\\
 & \cup \{\,w[0,2,\dots,k,1,\dots,k]\,|\,m \le k \le m+n-1\,\}.
\end{split}
\end{equation}}
Then we have a bijection 
\begin{equation}\label{eq:good Lyndon to positive root}
\begin{split}
 \xymatrixcolsep{2pc}\xymatrixrowsep{0pc}\xymatrix{
  \mc{GL}(\g) \ \ar@{->}[r] &\ \Phi^+ \\
  w=w[i_1,\dots ,i_r] \ \ar@{|->}[r] &\ |w|:=\alpha_{i_1}+\dots+\alpha_{i_r}\\
  l(\beta) \ \ar@{<-|}[r] &\ \beta }.
\end{split}
\end{equation}
where $l(\,\cdot\,)$ denotes its inverse.
We denote by $\prec$ the linear order on $\Phi^+$ induced from that on $\mc{GL}(\g)$ under \eqref{eq:good Lyndon to positive root}.

Let $\beta\in \Phi^+$ be given. Put ${\rm ht}(\beta)=r$ if $\beta=\alpha_{i_1}+\cdots+\alpha_{i_r}$ for some $i_1,\dots,i_r\in I$.
Let us define $\bff_\beta$ by using induction on ${\rm ht}(\beta)$.
If $r=1$, that is, $\beta=\alpha_i$ for some $i\in I$, then put $\bff_\beta=f_i$. 
Suppose that $r\ge 2$. 
Then it is straightforward to check that 
\begin{equation*}
l(\beta)=l(\beta_1)l(\beta_2), 
\end{equation*}
where $l(\beta_1)l(\beta_2)= \max\{\,l(\gamma_1)l(\gamma_2)\,|\,\gamma_1,\gamma_2\in \Phi^+,\, \gamma_1+\gamma_2=\beta,\,l(\gamma_1)<l(\gamma_2) \,\}$.
For such $\beta_1$ and $\beta_2$, let us assume that $l(\beta_1)$ is maximal.
Define
\begin{equation} \label{eq:PBW vector}
 \bff_\beta 
 = \frac{1}{[r_{\beta}+1]_{\beta_s}}\left[\bff_{\beta_2},\bff_{\beta_1}\right]_{\bq}
 = \frac{1}{[r_{\beta}+1]_{\beta_s}}\left(\,\bff_{\beta_2}\bff_{\beta_1}-\bq(\beta_1,\beta_2)^{-1}\bff_{\beta_1}\bff_{\beta_2}\,\right),
\end{equation}
where $r_{\beta}=\max\{\,p \in \Z_{\ge 0} \,|\,\beta_1-p\beta_2\in \Phi \,\}$, $\beta_s$ is such that $(\beta_s|\beta_s)=\min\{\,(\beta_1|\beta_1), (\beta_2|\beta_2)\,\}$, and $[k]_{\beta_s}=[k]_{i}$ when $(\beta_s|\beta_s)=(\alpha_i|\alpha_i)$. There is one exception when $\g=\mf{d}_{m|n}$ and $\beta=(m+1,m+1)$ (see Remark \ref{rem:root vectors}(2)).

\begin{rem}\label{rem:root vectors}
{\rm
(1)
In \cite{CHW}, $\bff_\beta$ is defined by letting $[k]_{\beta_s} = \{k\}_i$ in \eqref{eq:PBW vector} when $\beta_s$ is a non-isotropic odd root. The results in Section \ref{sec:main section} still hold with respect to this definition of $\bff_\beta$.

(2) The formula \eqref{eq:PBW vector} is based on \cite{CHW} (see also \cite[Section 4]{BKM}). According to \eqref{eq:PBW vector}, the root vector $\bff_{-2\de_{u+1}}$ for $\g = \mf{d}_{m|n}$ is given by
\begin{equation}\label{eq:two cases for root vector}
	\begin{cases}
		\displaystyle
		\frac{1}{[2]} \left[ f_{u}, \bff_{-\de_u-\de_{u+1}} \right]_{\bq} & \text{if $m+1\le u \le m+n-1$,} \\ 
		\left[ f_m, \bff_{-\de_{m}-\de_{m+1}} \right]_{\bq} & \text{if $u = m$,}
	\end{cases}
\end{equation}
which is divided into two cases since $-2\de_m$ is not a root. 
But in what follows, we take
\begin{equation*}
 \bff_{-2\de_{u+1}} = \frac{1}{[2]} \left[ f_u, \bff_{-\de_{u}-\de_{u+1}} \right]_{\bq}\quad (m\le u \le m+n-1),
\end{equation*} 
for $u=m$, as its definition, since it provides a simple and more uniform description of the results in Sections \ref{sec:PBW type basis} and \ref{sec:main section} compared to \eqref{eq:two cases for root vector}.}

\end{rem}

\begin{prop}\label{prop:PBW basis}
Let
\begin{equation*}
{B}=\Bigg\{\,\prod^{\rightarrow}_{\alpha\in\Phi^+}\bff_\alpha^{m_\alpha} \,\Bigg|\, \text{$m_{\alpha}\in\Z_{\ge 0}$ $(\alpha\in \Phi^+_{\ov 0}\cup \Phi^+_{\rm{n\text{-}iso}})$,\quad $m_{\alpha}=0,1$  $(\alpha\in \Phi^+_{{\rm iso}})$}\,\Bigg\},
\end{equation*} 
where the product is taken in increasing order with respect to $\prec$. Then
${B}$ is a $\Bbbk$-basis of $\U(\g)^-$.
\end{prop}
\pf Associated to $\mc{GL}(\mf{g})$ \eqref{eq:good Lyndon words}, there is a PBW type basis $B'$ of $U_q(\mf{g})^-:=\langle\, F_i \,|\,i\in I\,\rangle$ \cite[Section 5]{CHW}, which is constructed following \cite{Le}, where each root vector is defined in the same manner as in \eqref{eq:PBW vector} except the commutator $[\ ,\ ]_{\bf q}$ (cf.~\cite{Ya94}). 
Then one can check 
\begin{equation*}
	\tau(B') = 
	\Bigg\{\, {\bm \sigma}  \prod^{\rightarrow}_{\alpha\in\Phi^+}\bff_\alpha^{m_\alpha} \,\Bigg|\, \text{$m_{\alpha}\in\Z_{\ge 0}$ $(\alpha\in \Phi^+_{\ov 0}\cup \Phi^+_{\rm{n\text{-}iso}})$,\quad $m_{\alpha}=0,1$  $(\alpha\in \Phi^+_{{\rm iso}})$}\,\Bigg\},
\end{equation*}
where ${\bm \sigma}$ is a monomial in $\pm \sigma_a$'s depending on $(m_{\alpha})_{\alpha \in \Phi^+}$.
Hence our assertion follows (cf.~\cite[Proposition 2.5]{JKU}).
\qed\smallskip

For $\alpha\in \Phi^+$ and $k\ge 1$, we put
\begin{equation*}
 \bff_{\alpha}^{(k)} = 
\begin{cases}
  \bff_{\alpha}^{k}/[k]_i! & \text{if $\alpha\in \Phi^+_{\ov 0}$ and $(\alpha|\alpha)=(\alpha_i|\alpha_i)$},\\
  \bff_{\alpha}^{k}/\{k\}_i! & \text{if $\alpha\in \Phi^+_{\rm{n\text{-}iso}}$ and  $(\alpha|\alpha)=(\alpha_i|\alpha_i)$}.
\end{cases}
\end{equation*}

\begin{rem}{\rm 
We have $\bff_\alpha^2=0$ for $\alpha\in \Phi^+_{\rm{iso}}$ \cite[Corollary 5.2]{CHW}.
}
\end{rem}

Let $\U(\mf{l})$ be the subalgebra of $\U(\g)$ generated by $k_\mu, e_i, f_i$ for $\mu\in P$ and $i\in I\setminus\{0\}$, and $\U(\mf{l})^-=\U(\mf{l})\cap \U(\mf{g})^-$.
Let $\U(\mf{u}^-)$ be the $\Bbbk$-space spanned by the vectors in $B$ (see Proposition \ref{prop:PBW basis}) such that the product is over $\beta\in \Phi^+(\mf{u})$.
Since $\beta\prec\alpha$ for $\beta\in \Phi^+(\mf{u}), \alpha\in  \Phi^+(\mf{l})$ by \eqref{eq:good Lyndon words}, we have as a $\Bbbk$-space
\begin{equation}\label{eq:parabolic decomp}
 \U(\g)^-\cong \U(\mf{u}^-)\ot \U(\mf{l})^-,
\end{equation}
by Proposition \ref{prop:PBW basis}.

\subsection{Subalgebra $\mc{N}$ associated with $\Phi^+(\mf{u})$} \label{subsec:subalg N}
The roots in $\Phi^+(\mf{u})$ are of particular importance in this paper. Indeed, the space $\U(\mf{u}^-)$ can be viewed as a super analogue of a quantum nilpotent subalgebra associated to $\Phi^+(\mf{u})$. In this subsection, we show that it is a well-defined subalgebra of $\U(\g)^-$ with an explicit presentation.

For $\beta\in \Phi^+(\mf{u})$, let us simply write  
\begin{equation*}
 \beta= 
\begin{cases}
 (i,i) & \text{if $\beta = -2\de_i/r_{\g}$ for $1\le i\le m+n$},\\
 (i,j) & \text{if $\beta = -\de_i-\de_j$ for $1\le i<j\le m+n$}.
\end{cases}
\end{equation*}
By \eqref{eq:roots in radical}, we have
\begin{equation}\label{eq:nilradical roots}
 \Phi^+(\mf{u})=
\begin{cases}
 \{\,(i,i)\,|\,1\le i\le m+n\,\} \cup \{\,(i,j)\,|\,1\le i<j\le m+n\,\} & \text{if $\g=\mf{b}_{m|n}$},\\ 
  \{\,(i,i)\,|\,1\le i\le m\,\} \cup \{\,(i,j)\,|\,1\le i<j\le m+n\,\} & \text{if $\g=\mf{c}_{m|n}$},\\
   \{\,(i,i)\,|\,m+1\le i\le m+n\,\} \cup \{\,(i,j)\,|\,1\le i<j\le m+n\,\} & \text{if $\g=\mf{d}_{m|n}$}.
\end{cases}
\end{equation}
 
\begin{ex} \label{ex: roots of u}
{\em Let us illustrate $\Phi^+(\mf{u})$ for $\mf{b}_{2|3}$, $\mf{c}_{2|3}$, and $\mf{d}_{3|3}$, where each root is located on an triangle-shaped arrangement of $\scalebox{1.3}{$\odot$}$'s as shown below with respect to the linear order $\prec$. 
The notation $\underset{\beta}{\overset{k}{\scalebox{1.3}{$\odot$}}}$ means the $k$-th root $\beta$ with respect to $\prec$ and $\scalebox{1.3}{$\odot$}$ denotes the type of $\beta$, where $\scalebox{1.3}{$\odot$} = \bigcirc$ (even), $\ot$ (isotropic) or $\raisebox{-0.1cm}{{\scalebox{2.3}{$\bullet$}}}$ (non-isotropic odd).\smallskip

For $\g = \mf{b}_{2|3}$, we have
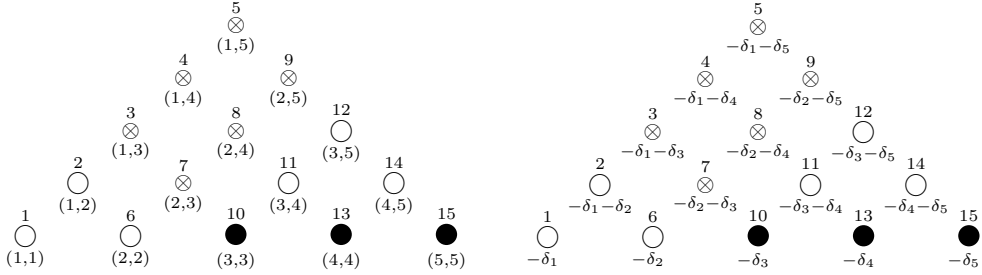
\begin{equation*}
	\begin{tikzpicture}[baseline=(current  bounding  box.center), every node/.style={scale=0.9}]
		\node (glw_00f) at (0,0) {$\underset{(1,1)}{\overset{\tiny 1}{\bigcirc}}$};
		\node (glw_01m) at (0.7,0.7) {$\underset{(1,2)}{\overset{\tiny 2}{\bigcirc}}$};
		\node (glw_02m) at (1.4,1.4) {$\underset{(1,3)}{\overset{\tiny 3}{\ot}}$};
		\node (glw_03m) at (2.1,2.1) {$\underset{(1,4)}{\overset{\tiny 4}{\ot}}$};
		\node (glw_04m) at (2.8,2.8) {$\underset{(1,5)}{\overset{\tiny 5}{\ot}}$};
		\node (glw_01f) at (1.4,0) {$\underset{(2,2)}{\overset{\tiny 6}{\bigcirc}}$};
		\node (glw_12m) at (2.1,0.7) {$\underset{(2,3)}{\overset{\tiny 7}{\ot}}$};
		\node (glw_13m) at (2.8,1.4) {$\underset{(2,4)}{\overset{\tiny 8}{\ot}}$};
		\node (glw_14m) at (3.5,2.1) {$\underset{(2,5)}{\overset{\tiny 9}{\ot}}$};
		\node (glw_02f) at (2.8,0) {$\underset{(3,3)}{\overset{\tiny 10}{{\scalebox{2.3}{$\bullet$}}}}$};
		\node (glw_23m) at (3.5,0.7) {$\underset{(3,4)}{\overset{\tiny 11}{\bigcirc}}$};		
		\node (glw_24m) at (4.2,1.4) {$\underset{(3,5)}{\overset{\tiny 12}{\bigcirc}}$};
		\node (glw_03f) at (4.2,0) {$\underset{(4,4)}{\overset{\tiny 13}{{\scalebox{2.3}{$\bullet$}}}}$};
		\node (glw_34m) at (4.9,0.7) {$\underset{(4,5)}{\overset{\tiny 14}{\bigcirc}}$};
		\node (glw_04m) at (5.6,0) {$\underset{(5,5)}{\overset{\tiny 15}{{\scalebox{2.3}{$\bullet$}}}}$};
		\end{tikzpicture}
		\,\,\quad\,\,
		\begin{tikzpicture}[baseline=(current  bounding  box.center), every node/.style={scale=0.9}]
		\node (glw_00f) at (0,0) {$\underset{-\delta_1\,\,\,}{\overset{\tiny 1}{\bigcirc}}$};
		\node (glw_01m) at (0.7,0.7) {$\underset{-\delta_1-\delta_2}{\overset{\tiny 2}{\bigcirc}}$};
		\node (glw_02m) at (1.4,1.4) {$\underset{-\delta_1-\delta_3}{\overset{\tiny 3}{\ot}}$};
		\node (glw_03m) at (2.1,2.1) {$\underset{-\delta_1-\delta_4}{\overset{\tiny 4}{\ot}}$};
		\node (glw_04m) at (2.8,2.8) {$\underset{-\delta_1-\delta_5}{\overset{\tiny 5}{\ot}}$};
		\node (glw_01f) at (1.4,0) {$\underset{-\delta_2\,\,\,}{\overset{\tiny 6}{\bigcirc}}$};
		\node (glw_12m) at (2.1,0.7) {$\underset{-\delta_2-\delta_3}{\overset{\tiny 7}{\ot}}$};
		\node (glw_13m) at (2.8,1.4) {$\underset{-\delta_2-\delta_4}{\overset{\tiny 8}{\ot}}$};
		\node (glw_14m) at (3.5,2.1) {$\underset{-\delta_2-\delta_5}{\overset{\tiny 9}{\ot}}$};
		\node (glw_02f) at (2.8,0) {$\underset{-\delta_3\,\,\,}{\overset{\tiny 10}{{\scalebox{2.3}{$\bullet$}}}}$};
		\node (glw_23m) at (3.5,0.7) {$\underset{-\delta_3-\delta_4}{\overset{\tiny 11}{\bigcirc}}$};		
		\node (glw_24m) at (4.2,1.4) {$\underset{-\delta_3-\delta_5}{\overset{\tiny 12}{\bigcirc}}$};
		\node (glw_03f) at (4.2,0) {$\underset{-\delta_4\,\,\,}{\overset{\tiny 13}{\scalebox{2.3}{$\bullet$}}}$};
		\node (glw_34m) at (4.9,0.7) {$\underset{-\delta_4-\delta_5}{\overset{\tiny 14}{\bigcirc}}$};
		\node (glw_04m) at (5.6,0) {$\underset{-\delta_5\,\,\,}{\overset{\tiny 15}{{\scalebox{2.3}{$\bullet$}}}}$};
		\end{tikzpicture}
\end{equation*}
\smallskip

For $\g = \mf{c}_{2|3}$, we have
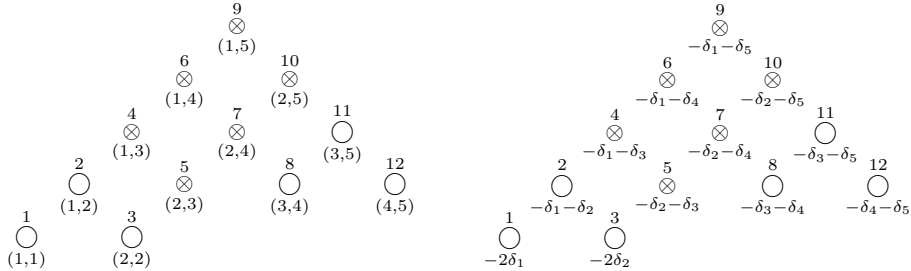
\begin{equation*}
	\begin{tikzpicture}[baseline=(current  bounding  box.center), every node/.style={scale=0.9}]
		\node (glw_00f) at (0,0) {$\underset{(1,1)}{\overset{\tiny 1}{\bigcirc}}$};
		\node (glw_01m) at (0.7,0.7) {$\underset{(1,2)}{\overset{\tiny 2}{\bigcirc}}$};
		\node (glw_02m) at (1.4,1.4) {$\underset{(1,3)}{\overset{\tiny 4}{\ot}}$};
		\node (glw_03m) at (2.1,2.1) {$\underset{(1,4)}{\overset{\tiny 6}{\ot}}$};
		\node (glw_04m) at (2.8,2.8) {$\underset{(1,5)}{\overset{\tiny 9}{\ot}}$};
		\node (glw_01m) at (1.4,0) {$\underset{(2,2)}{\overset{\tiny 3}{\bigcirc}}$};
		\node (glw_12m) at (2.1,0.7) {$\underset{(2,3)}{\overset{\tiny 5}{\ot}}$};
		\node (glw_13m) at (2.8,1.4) {$\underset{(2,4)}{\overset{\tiny 7}{\ot}}$};
		\node (glw_14m) at (3.5,2.1) {$\underset{(2,5)}{\overset{\tiny 10}{\ot}}$};
		\node (glw_23m) at (3.5,0.7) {$\underset{(3,4)}{\overset{\tiny 8}{\bigcirc}}$};		
		\node (glw_24m) at (4.2,1.4) {$\underset{(3,5)}{\overset{\tiny 11}{\bigcirc}}$};
		\node (glw_34m) at (4.9,0.7) {$\underset{(4,5)}{\overset{\tiny 12}{\bigcirc}}$};
		\end{tikzpicture}
		\,\,\,\quad\,\,\,
		\begin{tikzpicture}[baseline=(current  bounding  box.center), every node/.style={scale=0.9}]
		\node (glw_00f) at (0,0) {$\underset{-2\delta_1\,\,\,}{\overset{\tiny 1}{\bigcirc}}$};
		\node (glw_01m) at (0.7,0.7) {$\underset{-\delta_1-\delta_2}{\overset{\tiny 2}{\bigcirc}}$};
		\node (glw_02m) at (1.4,1.4) {$\underset{-\delta_1-\delta_3}{\overset{\tiny 4}{\ot}}$};
		\node (glw_03m) at (2.1,2.1) {$\underset{-\delta_1-\delta_4}{\overset{\tiny 6}{\ot}}$};
		\node (glw_04m) at (2.8,2.8) {$\underset{-\delta_1-\delta_5}{\overset{\tiny 9}{\ot}}$};
		\node (glw_01f) at (1.4,0) {$\underset{-2\delta_2\,\,\,}{\overset{\tiny 3}{\bigcirc}}$};
		\node (glw_12m) at (2.1,0.7) {$\underset{-\delta_2-\delta_3}{\overset{\tiny 5}{\ot}}$};
		\node (glw_13m) at (2.8,1.4) {$\underset{-\delta_2-\delta_4}{\overset{\tiny 7}{\ot}}$};
		\node (glw_14m) at (3.5,2.1) {$\underset{-\delta_2-\delta_5}{\overset{\tiny 10}{\ot}}$};
		\node (glw_23m) at (3.5,0.7) {$\underset{-\delta_3-\delta_4}{\overset{\tiny 8}{\bigcirc}}$};		
		\node (glw_24m) at (4.2,1.4) {$\underset{-\delta_3-\delta_5}{\overset{\tiny 11}{\bigcirc}}$};
		\node (glw_34m) at (4.9,0.7) {$\underset{-\delta_4-\delta_5}{\overset{\tiny 12}{\bigcirc}}$};
		\end{tikzpicture}
\end{equation*}
\smallskip

For $\g = \mf{d}_{3|3}$, we have
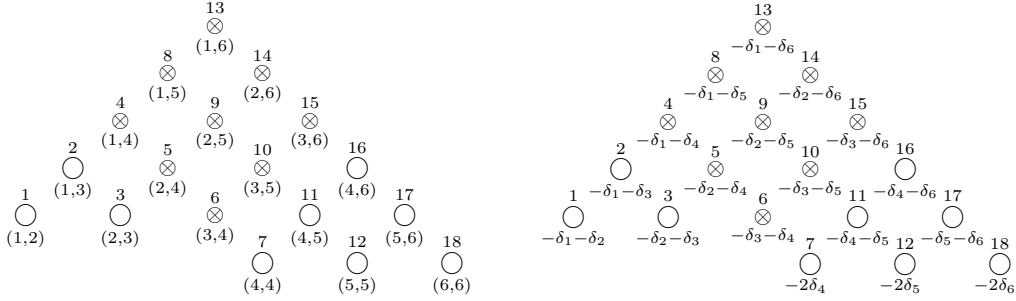
\begin{equation*}
	\begin{tikzpicture}[baseline=(current  bounding  box.center), every node/.style={scale=0.9}, scale=0.9]
		\node (glw_00f) at (0,0) {$\underset{(1,2)}{\overset{\tiny 1}{\bigcirc}}$};
		\node (glw_01m) at (0.7,0.7) {$\underset{(1,3)}{\overset{\tiny 2}{\bigcirc}}$};
		\node (glw_02m) at (1.4,1.4) {$\underset{(1,4)}{\overset{\tiny 4}{\ot}}$};
		\node (glw_03m) at (2.1,2.1) {$\underset{(1,5)}{\overset{\tiny 8}{\ot}}$};
		\node (glw_04m) at (2.8,2.8) {$\underset{(1,6)}{\overset{\tiny 13}{\ot}}$};
		\node (glw_01f) at (1.4,0) {$\underset{(2,3)}{\overset{\tiny 3}{\bigcirc}}$};
		\node (glw_12m) at (2.1,0.7) {$\underset{(2,4)}{\overset{\tiny 5}{\ot}}$};
		\node (glw_13m) at (2.8,1.4) {$\underset{(2,5)}{\overset{\tiny 9}{\ot}}$};
		\node (glw_14m) at (3.5,2.1) {$\underset{(2,6)}{\overset{\tiny 14}{\ot}}$};
		\node (glw_02f) at (2.8,0) {$\underset{(3,4)}{\overset{\tiny 6}{\ot}}$};
		\node (glw_23m) at (3.5,0.7) {$\underset{(3,5)}{\overset{\tiny 10}{\ot}}$};		
		\node (glw_24m) at (4.2,1.4) {$\underset{(3,6)}{\overset{\tiny 15}{\ot}}$};
		\node (glw_03f) at (3.5,-0.7) {$\underset{(4,4)}{\overset{\tiny 7}{\bigcirc}}$};
		\node (glw_03f) at (4.2,0) {$\underset{(4,5)}{\overset{\tiny 11}{\bigcirc}}$};
		\node (glw_34m) at (4.9,0.7) {$\underset{(4,6)}{\overset{\tiny 16}{\bigcirc}}$};
		\node (glw_04m) at (4.9,-0.7) {$\underset{(5,5)}{\overset{\tiny 12}{\bigcirc}}$};
		\node (glw_04m) at (5.6,0) {$\underset{(5,6)}{\overset{\tiny 17}{\bigcirc}}$};
		\node (glw_04m) at (6.3,-0.7) {$\underset{(6,6)}{\overset{\tiny 18}{\bigcirc}}$};
		\end{tikzpicture}
		\,\,\,\quad\,\,\,
	\begin{tikzpicture}[baseline=(current  bounding  box.center), every node/.style={scale=0.9}, scale=0.9]
		\node (glw_00f) at (0,0) {$\underset{-\delta_1-\delta_2}{\overset{\tiny 1}{\bigcirc}}$};
		\node (glw_01m) at (0.7,0.7) {$\underset{-\delta_1-\delta_3}{\overset{\tiny 2}{\bigcirc}}$};
		\node (glw_02m) at (1.4,1.4) {$\underset{-\delta_1-\delta_4}{\overset{\tiny 4}{\ot}}$};
		\node (glw_03m) at (2.1,2.1) {$\underset{-\delta_1-\delta_5}{\overset{\tiny 8}{\ot}}$};
		\node (glw_04m) at (2.8,2.8) {$\underset{-\delta_1-\delta_6}{\overset{\tiny 13}{\ot}}$};
		\node (glw_01f) at (1.4,0) {$\underset{-\delta_2-\delta_3}{\overset{\tiny 3}{\bigcirc}}$};
		\node (glw_12m) at (2.1,0.7) {$\underset{-\delta_2-\delta_4}{\overset{\tiny 5}{\ot}}$};
		\node (glw_13m) at (2.8,1.4) {$\underset{-\delta_2-\delta_5}{\overset{\tiny 9}{\ot}}$};
		\node (glw_14m) at (3.5,2.1) {$\underset{-\delta_2-\delta_6}{\overset{\tiny 14}{\ot}}$};
		\node (glw_02f) at (2.8,0) {$\underset{-\delta_3-\delta_4}{\overset{\tiny 6}{\ot}}$};
		\node (glw_23m) at (3.5,0.7) {$\underset{-\delta_3-\delta_5}{\overset{\tiny 10}{\ot}}$};		
		\node (glw_24m) at (4.2,1.4) {$\underset{-\delta_3-\delta_6}{\overset{\tiny 15}{\ot}}$};
		\node (glw_03f) at (3.5,-0.7) {$\underset{-2\delta_4\,\,\,}{\overset{\tiny 7}{\bigcirc}}$};
		\node (glw_03f) at (4.2,0) {$\underset{-\delta_4-\delta_5}{\overset{\tiny 11}{\bigcirc}}$};
		\node (glw_34m) at (4.9,0.7) {$\underset{-\delta_4-\delta_6}{\overset{\tiny 16}{\bigcirc}}$};
		\node (glw_04m) at (4.9,-0.7) {$\underset{-2\delta_5\,\,\,}{\overset{\tiny 12}{\bigcirc}}$};
		\node (glw_04m) at (5.6,0) {$\underset{-\delta_5-\delta_6}{\overset{\tiny 17}{\bigcirc}}$};
		\node (glw_04m) at (6.3,-0.7) {$\underset{-2\delta_6\,\,\,}{\overset{\tiny 18}{\bigcirc}}$};
	\end{tikzpicture}\end{equation*}
	\smallskip
\noindent
}
\end{ex}
 
We give explicit formulas of the $q$-commutator for $\bff_\alpha$ and $\bff_\beta$ for $\alpha,\, \beta \in \Phi^+(\mf{u})$, which play an important role in this paper.

\begin{prop}\label{prop:commutation rel}
For $\alpha,\, \beta \in \Phi^+(\mf{u})$ with $\alpha=(i_1,j_1) \prec \beta=(i_2,j_2)$, we have the following formulas for $[\bff_\beta, \bff_\alpha]_{\bq}$:
\begin{itemize}
 \item[(1)] $\mf{g}=\mf{b}_{m|n}$
\begin{equation*}
\begin{cases}
 (q+q^{-1})\bff_{(i_1,i_2)}  & \text{if $i_1=j_1<i_2=j_2$}, \\ 
 (q^{-2}-q^2) \bff_{(i_1,j_2)}\bff_{(i_2,i_2)} & \text{if $i_1=j_1<i_2<j_2$}, \\ 
 (q^{-2}-q^2) \bff_{(i_1,i_2)}\bff_{(j_1,j_1)} & \text{if $i_1<j_1<i_2=j_2$}, \\ 
  (q^{-2}-q^2) \bff_{(i_1,j_2)}\bff_{(i_2,j_1)} & \text{if $i_1<i_2<j_1<j_2$}, \\ 
 (q^{-1}-q)({\tt q}_{j_1}^{-1}-1)\bff_{(i_1,j_2)}\bff_{(j_1,j_1)}^2 & \text{if $i_1<j_1=i_2<j_2$}, \\ 
    a\, \bff_{(i_1,j_2)}\bff_{(j_1,i_2)}
+ b\, \bff_{(i_1,i_2)}\bff_{(j_1,j_2)}
+ c\, \bff_{(i_1,j_2)}\bff_{(j_1,j_1)}\bff_{(i_2,i_2)} & \text{if $i_1<j_1<i_2<j_2$}, \\ 
  0  & \text{otherwise},
\end{cases}
\end{equation*}
where $a=q^{-4} - q^{-2} - 2 + q^{2} + q^{4}$, $b=q^{-2} - q^2$, and $c=q^{-3} - q^{-1} - q^{1} + q^{3}$,

\item[(2)] $\mf{g}=\mf{c}_{m|n}$ 
\begin{equation*}
\begin{cases}
 (q^{-2}-1)(q+q^{-1})^{-1}\bff_{(i_1,i_2)}^{2}  & \text{if $i_1=j_1<i_2=j_2\le m$}, \\ 
 (q^{-2}-1) \bff_{(i_1,i_2)}\bff_{(j_1,j_2)} & \text{if $i_1=j_1<i_2<j_2$, $i_1\le m$}, \\ 
 (q^{-2}-1) \bff_{(i_1,i_2)}\bff_{(j_1,j_2)} & \text{if $i_1<j_1<i_2=j_2\le m$}, \\ 
 (q^{-1}-q) \bff_{(i_2,j_1)}\bff_{(i_1,j_2)} & \text{if $i_1<i_2<j_1<j_2$}, \\ 
 (q^{-2}-q^2)  \bff_{(i_2,j_1)}\bff_{(i_1,j_2)} & \text{if $i_1<j_1=i_2<j_2$, $j_1\le m$}, \\ 
 (q^{-1}-q) \bff_{(i_1,i_2)}\bff_{(j_1,j_2)}+
 (q^{-2}-1) \bff_{(j_1,i_2)}\bff_{(i_1,j_2)} & \text{if $i_1<j_1<i_2<j_2$}, \\ 
  0  & \text{otherwise}, 
\end{cases}
\end{equation*}

\item[(3)] $\mf{g}=\mf{d}_{m|n}$
\begin{equation*}
\begin{cases}
 (1-q^{2})(q+q^{-1})^{-1}\bff_{(i_1,i_2)}^{2}  & \text{if $m<i_1=j_1<i_2=j_2$}, \\ 
 (1-q^{2}) \bff_{(i_1,i_2)}\bff_{(j_1,j_2)} & \text{if $m<i_1=j_1<i_2<j_2$}, \\ 
 (1-q^{2}) \bff_{(i_1,i_2)}\bff_{(j_1,j_2)} & \text{if $i_1<j_1<i_2=j_2$, $m<i_2$}, \\ 
 (q^{-1}-q) \bff_{(i_2,j_1)}\bff_{(i_1,j_2)} & \text{if $i_1<i_2<j_1<j_2$}, \\ 
 (q^{-2}-q^2)  \bff_{(i_2,j_1)}\bff_{(i_1,j_2)} & \text{if $i_1<j_1=i_2<j_2$, $m<j_1$}, \\ 
 (q^{-1}-q) \bff_{(i_1,i_2)}\bff_{(j_1,j_2)}+
 (q^{2}-1) \bff_{(j_1,i_2)}\bff_{(i_1,j_2)} & \text{if $i_1<j_1<i_2<j_2$}, \\ 
  0  & \text{otherwise}. 
\end{cases}
\end{equation*}
\end{itemize}
\end{prop}

\begin{cor}
The subspace $\U(\mf{u}^-)$ is a $\Bbbk$-subalgebra of $\U(\g)^-$ generated by $\bff_\beta$ for $\beta\in \Phi^+(\mf{u})$, which we denote by $\mc{N}(\mf{g})=\mc{N}$.
\end{cor}

\begin{rem}\label{rem:commutation relation}
{\rm (1) Suppose that $\mf{g}=\mf{c}_{m|n}, \mf{d}_{m|n}$. By Proposition \ref{prop:commutation rel}(2) and (3), a basis element of $B$ in Proposition \ref{prop:PBW basis} is equal to the one whose product is taken with respect to $\prec'$, where $(i,j)\prec'(k,l)$ if and only if $i<k$ or ($i=k$ and $j<l$) as in $\mf{b}_{m|n}$.

(2) The formulas in Proposition \ref{prop:commutation rel} are also available when $n=0$, which provide explicit coefficients in the Levendorskii-Soibelman formula \cite{LS}, but which we could not find in the literature.
}
\end{rem}

\begin{cor}
For $(\mf{x},\mf{y}) =(\mf{c},\mf{d}), (\mf{d},\mf{c})$, there exists an anti-isomorphism of $\Q$-algebras $\omega : \mc{N}(\mf{x}_{m|n})\longrightarrow \mc{N}(\mf{y}_{n|m})$ such that
\begin{equation*}
\begin{split}
\omega(q)&=-q^{-1},\\
\omega\left(\bff_{(i,j)}\right)&=
\begin{cases}
\bff_{(m+n+1-j,m+n+1-i)} & \text{if $i\neq j$},\\
-\bff_{(m+n+1-i,m+n+1-i)} & \text{if $i= j$}.
\end{cases}
\end{split} 
\end{equation*}
\end{cor}\smallskip
 
In order to prove Proposition \ref{prop:commutation rel}, we consider an embedding of $\U(\g)^-$ into a quantum shuffle algebra following \cite{CHW,Le}, and then use the image of $\bff_\alpha, \bff_\beta$ ($\alpha, \beta \in \Phi^+(\mf{u})$) to compute $[\bff_\beta, \bff_\alpha]_{\bq}$.

Let $\mc{F}$ be the free associative algebra over $\Bbbk$ generated by $\mc{A}=\{\,w_0<\dots<w_{m+n-1}\,\}$ in Section \ref{subsec:root vectors}, where $\mc{W}$ can be viewed as a $\Bbbk$-basis of $\mc{F}$. 
Define a bilinear map $\ast$ on $\mc{F}$ inductively as follows:

\begin{equation*} 
 (ww_a)\ast (w'w_{a'}) =(w\ast w'w_{a'}) w_a + \bq(|ww_a|,|w_{a'}|)^{-1} (ww_a\ast w')w_{a'}
\end{equation*}
for words $w, w'\in \mc{W}$ and letters $a,a'\in I$, where we assume that $w\ast w_a=w_a\ast w=w_a$ for an empty word $w$, and $|w|$ is defined by extending $|\cdot |$ in \eqref{eq:good Lyndon to positive root} to $\mc{W}$. Then $\mc{F}$ is an associative $\Bbbk$-algebra with respect to $\ast$, which we denote by $(\mc{F},\ast)$. We have the following, which can be seen in the same way as in \cite[Corollary 3.4]{CHW}.

\begin{prop}
There exists an injective homomorphism of $\Bbbk$-algebras $\Psi: \U(\mf{g})^-\longrightarrow (\mc{F},\ast)$ such that $\Psi(f_i)=w_i=w[i]$ for $i\in I$.
\end{prop}
\smallskip

We can check the following by using \eqref{eq:PBW vector} and induction on ${\rm ht}(\beta)$ for $\beta\in \Phi^+(\mf{u})$, which is an analogue of \cite[Sections 8.2--8.4]{Le}.

\begin{lem}\label{lem:image of root vector}
For $\beta=(i,j)\in \Phi^+(\mf{u})$, $\Psi(\bff_{\beta})$ is given as follows:
\begin{itemize}
\item[(1)] $\mf{g}=\mf{b}_{m|n}$
\begin{equation*}
\begin{cases}
 w[0] & \text{if $(i,j)=(1,1)$},\\
 a_{i-1} w[0,1, \dots , i-1] & \text{if $2\le i=j$},\\
 (q^{-1}-q)a_{j-1}\, w[0,0,1, \dots ,j-1]  & \text{if $1=i<j$},\\
 (q^{-1}-q) a_{i-1}a_{j-1}\,  w[0](w[1,\dots ,i-1] \ast  w[0,\dots ,j-1]) & \text{if $2\le i<j$},
\end{cases}
\end{equation*}

\item[(2)] $\mf{g}=\mf{c}_{m|n}$
\begin{equation*}
\begin{cases}
 w[0] & \text{if $(i,j)=(1,1)$},\\
 (1+q^2)a_{j-1}\, w[0,1, \cdots ,j-1] & \text{if $1=i<j$},\\
 (1+q^2) a_{i-1}a_{j-1}\,  w[0](w[1,\dots, i-1] \ast  w[1,\dots ,j-1]) & \text{if $2\le i<j$},\\
  q a_{i-1}^2\,  w[0](w[1,\dots ,i-1] \ast  w[1,\dots ,i-1]) & \text{if $2\le i=j\le m$},
\end{cases}
\end{equation*}

\item[(3)] $\mf{g}=\mf{d}_{m|n}$
\begin{equation*}
\begin{cases}
 w[0] & \text{if $(i,j)=(1,2)$},\\
 b_{j-1}\, w[0, 2,3, \dots, j-1] & \text{if $1=i<2<j$},\\
 \begin{split}
 & (1-q^2) b_{i-1}b_{j-1}\,w[0] \cdot \\ 
 & \quad \left(w[1,\dots, i-1] \ast w[2,\dots ,j-1] - q\, w[2,\dots ,i-1] \ast  w[1,\dots, j-1]\right)
 \end{split} & \text{if $2\le i<j$},\\
  q (1-q^2) b_{i-1}^2\,  w[0](w[1,\dots ,i-1] \ast  w[2,\dots, i-1]) & \text{if $m+1\le i=j$},
\end{cases}
\end{equation*}
\end{itemize}
where $a_i=(1-{\tt q}^2_{1})\dots (1-{\tt q}^2_{i})$ for $1\le i< m+n$ and $b_i=(1-{\tt q}^2_{2})\dots (1-{\tt q}^2_{i})$ for $2\le i< m+n$. 
\end{lem}

\begin{rem}{\rm
For $\beta=(i,j)\in \Phi^+(\mf{u})$, the maximal word in $\Psi(\bff_{\beta})$ is $l(\beta)$ \cite[Theorem 5.1]{CHW}. The words in $\mc{GL}(\mf{g})$ are called the {\em good Lyndon words} (see \cite{Le} and the references therein).
}
\end{rem}

\noindent{\em Proof of Proposition \ref{prop:commutation rel}.}
By \cite[Theorem 5.5]{CHW} and \eqref{eq:iso tau standard}, we have
\begin{equation}\label{eq:LS type formula}
 [\bff_\beta, \bff_\alpha]_{\bq} = 
 \sum_{\boldsymbol{\gamma}=(\gamma_1\preceq \dots\preceq\gamma_r)}
 c_{\boldsymbol{\gamma}} \bff_{{\bf \gamma}_1}\cdots \bff_{{\bf \gamma}_r},
\end{equation}
for some $c_{\boldsymbol{\gamma}}\in \Bbbk$,
where the sum is over $\boldsymbol{\gamma}=(\gamma_1\preceq \dots\preceq\gamma_r)$ such that $l(\alpha) l(\beta)\le l(\gamma_r) \cdots l(\gamma_1)< l(\beta)l(\alpha)$. 
Applying $\Psi$ to \eqref{eq:LS type formula}, we get
\begin{equation*}\label{eq:LS type formula-shufflw}
\Psi(\bff_\beta) \ast \Psi(\bff_\alpha) - \bq(\alpha,\beta)^{-1} \Psi(\bff_\alpha)\ast \Psi(\bff_\beta)  
 = 
 \sum_{\boldsymbol{\gamma}=(\gamma_1\preceq \dots\preceq\gamma_r)}c_{\boldsymbol{\gamma}} 
\Psi(\bff_{{\bf \gamma}_1})\ast \cdots\ast\Psi(\bff_{{\bf \gamma}_r}).
\end{equation*}
First, we find the possible $\boldsymbol{\gamma}$'s appearing in \eqref{eq:LS type formula}, say $\boldsymbol{\gamma}^{(1)}, \dots , \boldsymbol{\gamma}^{(s)}$ with $\boldsymbol{\gamma}^{(i)}=(\gamma^{(i)}_1\preceq \dots\preceq\gamma^{(i)}_{r_i})$. 
For $1\le i\le s$, the maximal word in 
$\Psi(\bff_{{\bf \gamma}^{(i)}_1})\ast \cdots\ast\Psi(\bff_{{\bf \gamma}^{(i)}_{r_i}})$
is $w_{\boldsymbol{\gamma}^{(i)}}:=l(\gamma^{(i)}_{r_i})\dots l(\gamma^{(i)}_{1})$ by \cite[Theorem 5.1]{CHW}. 
We assume that $w_{\boldsymbol{\gamma}^{(1)}}>\dots>w_{\boldsymbol{\gamma}^{(s)}}$.
Then we obtain $c_{\boldsymbol{\gamma}^{(i)}}$ inductively by using Lemma \ref{lem:image of root vector} and computing the coefficient of $w_{\boldsymbol{\gamma}^{(i)}}$ in
\begin{equation*}
\begin{split}
  \Psi(\bff_\beta) \ast \Psi(\bff_\alpha) -& \bq(\alpha,\beta)^{-1} \Psi(\bff_\alpha)\ast \Psi(\bff_\beta) \\
 &- c_{\boldsymbol{\gamma}^{(1)}}\Psi(\bff_{{\bf \gamma}^{(1)}_1})\ast \cdots\ast\Psi(\bff_{{\bf \gamma}^{(1)}_{r_1}})
 \dots 
 - c_{\boldsymbol{\gamma}^{(i-1)}}\Psi(\bff_{{\bf \gamma}^{(i-1)}_1})\ast \cdots\ast\Psi(\bff_{{\bf \gamma}^{(i-1)}_{r_{i-1}}}).
\end{split}
\end{equation*}
We leave the details to the reader since it is straightforward computation.
 
\subsection{Quantum adjoint action on $\mc{N}$} 
Let $\mc{N}=\mc{N}(\g)$ be given in Section \ref{subsec:subalg N}.
Let us give a natural $\U(\mf{l})$-module structure on $\mc{N}$ via quantum adjoint.

Let ${\rm ad}_q: \U(\mf{g})  \longrightarrow {\rm End}_{\Bbbk}(\U(\mf{g}))$ be the adjoint representation of $\U(\g)$ with respect to \eqref{eq:comult-1},
where for $\mu \in P$, $i \in I$, and $u \in \U(\g)$,
\begin{equation} \label{eq: quantum adjoint}
\begin{split}
	& {\rm ad}_q(k_\mu)(u) = k_\mu u k_\mu^{-1}, \ \
	{\rm ad}_q(e_i)(u) = \left( e_i u - u e_i \right) k_i, \ \
	{\rm ad}_q(f_i)(u) = f_i u - k_i u k_i^{-1} f_i.
\end{split}
\end{equation}
We write ${\rm ad}_q(x)(u)=x\cdot u$ for simplicity.

%
\begin{lem}\label{lem:e-derivation and adjoint}
For $i\in I\setminus\{0\}$ and $u\in \mc{N}$, we have
\begin{equation*}
 e_i\cdot u = -\frac{\bq(\alpha_i,\alpha_i+|u|)^{-1}}{q_i-q_i^{-1}}\, e'_i(u).
\end{equation*}
\end{lem}
\pf For $u\in \U(\mf{g})^-$ and $i\in I$, we have
\begin{equation}\label{eq:e_i dot action}
 e_i\cdot u 
 = \left( \frac{k_ie_i''(u)-k_i^{-1}e_i'(u)}{q_i-q_i^{-1}} \right)k_i
 =\frac{k_ie_i''(u)k_i-k_i^{-1}e_i'(u)k_i}{q_i-q_i^{-1}} 
\end{equation}
(see \eqref{eq: derivation}). 
On the other hand, we have $e_i''(u)=0$ for $u\in \mc{N}$ and $i\in I\setminus\{0\}$, which can be proved inductively by using \eqref{eq:PBW vector} and
\begin{equation*}
e''_i[x,y]_{\bq}=
[e''_i(x),y]_{\bq} 
+ \bq(\alpha_i,|x|)^{-1}x e''_i(y) 
- \bq(|x|,|y|)^{-1}e''_i(y)x.
\end{equation*}
Hence we have the required identity from \eqref{eq:e_i dot action}. 
\qed
\smallskip

\begin{prop}\label{prop:e dot}
For $i\in I\setminus\{0\}$ and $\beta=(k,l)\in \Phi^+(\mf{u})$, we have
\begin{itemize}
 \item[(1)] $\mf{g}=\mf{b}_{m|n}$
\begin{equation*}
 e_{i}\cdot \bff_\beta=
\begin{cases}
 \bff_{(i,i)} & \text{if $(k,l) = (i+1, i+1)$},\\
 \bff_{(i,l)} & \text{if $k<l$ and $k=i+1$},\\
 \bff_{(k,i)} & \text{if $l=i+1$ and $k<i$},\\
 (\text{\em $\ttq_{i}^{-1}-1)(q+q^{-1})^{-1}\bff_{(i,i)}^{2}$} & \text{if $(k,l)=(i,i+1)$},\\
 0 & \text{otherwise},
\end{cases}
\end{equation*}

\item[(2)] $\mf{g}=\mf{c}_{m|n}$, $\mf{d}_{m|n}$
\begin{equation*}
 e_{i}\cdot \bff_\beta=
\begin{cases}
 c\,\bff_{(i,l)} & \text{if $k=i+1$},\\
 \bff_{(k,i)} & \text{if $l=i+1$ and $k < i$},\\  
 (q+q^{-1})\bff_{(i,i)} & \text{if $(k,l)=(i,i+1)$ with $i \le m$ and $\mf{g}=\mf{c}_{m|n}$},\\
 (q+q^{-1})\bff_{(i,i)} & \text{if $(k,l)=(i,i+1)$ with $i > m$ and $\mf{g}=\mf{d}_{m|n}$},\\
 0 & \text{otherwise},
\end{cases}
\end{equation*} 
where $c=-1$ for $\mf{d}_{m|n}$ and $(k,l)=(i+1,i+1)$, and $c=1$ otherwise.
\end{itemize}
\end{prop}
\pf
By Lemma \ref{lem:e-derivation and adjoint}, it is enough to consider $e_i' (\bff_\beta)$. Then we apply induction on ${\rm ht}(\beta)$ by using
\begin{equation}\label{eq:e'-formula}
e'_i[x,y]_{\bq}=
\bq(|x|,\alpha_i)[x,e'_i(y)]_{\bq}  + e'_i(x)y  - \bq(|x|-\alpha_i,|y|)^{-1}ye'_i(x).
\end{equation}
\qed

\begin{prop}\label{prop:f dot}
For $i\in I\setminus\{0\}$ and $\beta=(k,l)\in \Phi^+(\mf{u})$, we have
\begin{itemize}
 \item[(1)] $\mf{g}=\mf{b}_{m|n}$
\begin{equation*}
 f_{i}\cdot \bff_\beta=
\begin{cases}
  \bff_{(i+1,i+1)} & \text{if $(k,l) = (i,i)$},\\
  \bff_{(i+1,l)} & \text{if $k=i$ and $l > i+1$},\\
  \bff_{(k,i+1)} & \text{if $k<i$ and $l=i$},\\  
  (\text{\em $\ttq_{i+1}^{-1}-1)(q+q^{-1})^{-1}\bff_{(i+1,i+1)}^2$} & \text{if $(k,l) = (i, i+1)$},\\
  0 & \text{otherwise},
\end{cases}
\end{equation*}

\item[(2)] $\mf{g}=\mf{c}_{m|n}$, $\mf{d}_{m|n}$
\begin{equation*}
 f_{i}\cdot \bff_\beta=
\begin{cases}
  \bff_{(i+1,l)} & \text{if $k=i$ and $l > i+1$},\\
   c\,\bff_{(k,i+1)} & \text{if $l=i$},\\  
  (q+q^{-1})\bff_{(i+1,i+1)} & \text{if $(k,l)=(i,i+1)$ with $i\le m-1$ and $\mf{g}=\mf{c}_{m|n}$},\\
  (q+q^{-1})\bff_{(i+1,i+1)} & \text{if $(k,l)=(i,i+1)$ with $i> m-1$ and $\mf{g}=\mf{d}_{m|n}$},\\
 0 & \text{otherwise},
\end{cases}
\end{equation*}
where $c=-1$ for $\mf{d}_{m|n}$ and $(k,l)=(i,i)$, and $c=1$ otherwise.
\end{itemize}
\end{prop}
\pf Since $\beta \prec \alpha_i$, we may apply the same argument in the proof of Proposition \ref{prop:commutation rel}. We leave the details to the reader.
\qed

\begin{prop}\label{prop:poly alg}
The subalgebra $\mc{N}$ is a $\U(\mf{l})$-submodule of $\U(\mf{g})$ with respect to $ad_q$ such that
\begin{equation*}
	x\cdot m_{\mc N}(u_1\ot u_2) = m_{\mc N}(\Delta(x)\cdot(u_1\ot u_2))
\end{equation*}
for $x\in \U(\mf{l})$ and $u_1\ot u_2\in \mc{N}\ot \mc{N}$, where $m_{\mc N}$ denotes the multiplication on $\mc{N}$ and $\Delta$ is the comultiplication of $\U(\mf{g})$.
\end{prop}
\pf It follows from Propositions \ref{prop:e dot} and \ref{prop:f dot} that $\mc{N}$ is stable under the action of $\U(\mf{l})$. For $u=u_1u_2 \in \mc{N}$, we have
\begin{equation*}
\begin{split}
k_\mu\cdot u &= (k_\mu\cdot u_1)(k_\mu\cdot u_2),\\
e_i \cdot u  &= (e_i\cdot u_1)(k_i^{-1}\cdot u_2) + u_1(e_i\cdot u_2),\\
f_i \cdot u  &= (f_i\cdot u_1)u_2 + (k_i\cdot u_1)(f_i\cdot u_2), 
\end{split}
\qquad (\mu\in P,\, i\in I\setminus\{0\}).
\end{equation*}
This implies $x\cdot m_{\mc N}(u_1\ot u_2) = m_{\mc N}(\Delta(x)\cdot(u_1\ot u_2))$ for all $x\in \U(\mf{l})$.
\qed

\section{Crystal base of a parabolic Verma module}\label{sec:main section}

\subsection{Polynomial representations of $\U(\mf{l})$}\label{subsec:poly repn}
 Let us briefly review the results on the polynomial representations and their crystal bases for the quantum superalgebra $U_q(\gl(m|n))$ \cite{BKK}, which are slightly adjusted in terms of $\U(\mf{l})$-modules (cf.~\cite{JKU,KO}).

Let $P_{\geq0}=\sum_{a\in \I}\Z_{\ge 0}\delta _a$. 
Let $\cO_{\geq0}$ be the category of $\U(\mf{l})$-modules $V$ such that
$V=\bigoplus_{\mu\in P_{\ge 0}}V_\mu$ with $\dim V_\mu < \infty$.
We call $V\in \mc{O}_{\ge 0}$ a polynomial representation of $\U(\mf{l})$.

Let $\cP$ be the set of all partitions $\la=(\la_i)_{i\ge 1}$ with $\la_1\ge\la_2\ge\dots$, and let $\cP_{m|n}$ be the set of $\la\in\cP$ such that $\la_{m+1}\leq n$.
We identify $\la\in \cP_{m|n}$ with $\La_\la\in P_{\geq 0}$, where
\begin{equation}\label{eq:highest weight correspondence}
\La_\la = \la_1\de_1 +\cdots +\la_m\de_m + \mu_1\de_{m+1}+\cdots+\mu_n\de_{m+n},
\end{equation}
where $(\mu_1,\ldots,\mu_n)$ is the conjugate of the partition $(\la_{m+1},\la_{m+2},\ldots)$. Let $V_{\mf{l}}(\la)$ be the irreducible highest weight $\U(\mf{l})$-module with highest weight $\La_\la$. Then any irreducible $\U(\mf{l})$-module in $\cO_{\geq 0}$ is isomorphic to $V_{\mf{l}}(\la)$ for some $\la\in \cP_{m|n}$ by \cite[Proposition 4.5]{BKK}.

\begin{prop}\label{prop:semisimplicity}
For $V\in \cO_{\geq 0}$, $V$ is semisimple.
\end{prop}
\pf For $V\in \mc{O}_{\ge 0}$, let $V_d=\bigoplus_{\mu\in P_{\ge 0}, |\mu|=d}V_\mu$, where $|\mu|=\sum_{a\in \I}\mu_a$ for $\mu=\sum_{a\in \I}\mu_a\de_a\in P_{\ge 0}$. Then $V_d$ is a $\U(\mf{l})$-submodule of $V$, which is finite-dimensional since there are only finitely many $\mu\in P_{\ge 0}$ such that $|\mu|=d$. We have $V=\bigoplus_{d\in \Z_{\ge 0}}V_d$. Let $\mc{O}_{\ge 0,d}$ be the subcategory of $\mc{O}_{\ge 0}$ consisting of $V$ such that $V=V_d$. We have $\mc{O}_{\ge 0}=\bigoplus_{d\in \Z_{\ge 0}}\mc{O}_{\ge 0,d}$.

Let $V\in \mc{O}_{\ge 0,d}$ be given. 
We claim that $V$ is semisimple.
By Remark \ref{rem:comparison of reps of two algs}, we may assume that $V$ is a representation of $U_q(\mf{gl}(m|n))$, the subalgebra of $U_q(\g)$ generated by $K_\mu, E_i, F_i$ for $\mu\in P$ and $i\in I\setminus\{0\}$.
Suppose that we have an exact sequence
\begin{equation}\label{eq:extension of two poly}
 0 \longrightarrow U \longrightarrow V \longrightarrow W\longrightarrow 0,
\end{equation}
where $U=V_{\mf{l}}(\la)$ and $W=V_{\mf{l}}(\mu)$ 
for some $\la,\mu\in \cP_{m|n}$. Let $v_\la$ be a highest weight vector of $U\subset V$ and let $v_\mu\in V$ be given such that its image at $V/U$ is a highest weight vector of $W$.

Let $A=\mathbb{Q}[q,q^{-1}]$. Let $U_{A}=\sum_{{\bf i}}A F_{\bf i }v_\la$,
where the sum is over ${\bf i}=(i_1,\dots,i_r)$ for $r\ge 0$ and $i_k\in I\setminus\{0\}$, and $F_{\bf i}=F_{i_1}\dots F_{i_r}$ with $F_{\bf i}=1$ for $r=0$. We have $U_A=\bigoplus_{\nu\in P_{\ge 0}}U_{\nu,A}$, where $U_{\nu,A}=U_A\cap U_\nu$.
Recall that $U_{\nu,A}$ is a free $A$-module and $U_{\nu,A}\otimes_A \Bbbk = U_{\nu}$ for $\nu$ (cf.~\cite[Chapter 5]{Ja}). Let $W_A$ be defined in the same way.
Also we let
\begin{equation*}
V_{A}=\sum_{{\bf i}}A F_{\bf i}v_\la + \sum_{{\bf j}}A F_{\bf j}v_\mu,\end{equation*}
where the sum is over all ${\bf i}=(i_1,\dots,i_r)$ and ${\bf j}=(j_1,\dots,j_s)$ with $r,s\ge 0$ and $i_k,j_l\in I\setminus\{0\}$.
Here we may assume that $E_i v_\mu \in U_A$ for all $i\in I\setminus\{0\}$ by replacing $v_\mu$ with $f(q)v_\mu$ for some $f(q)\in \mathbb{Q}[q]$.

Let $X$ denote one of $U$, $V$ and $W$.
Then the ${A}$-module $X_{A}$ is invariant under $E_i$, $F_i$, $K_{\pm \de_a}$ and $H_a=\frac{K_{\de_a}-K_{-\de_a}}{v-v^{-1}}$ ($v=q^{r_\g}$) for $i\in I\setminus\{0\}$ and $a\in \I$. 
Let $\ov{X}=X_{A}\otimes_{A}\mathbb{C}$ and $\ov{X}_{\nu}=X_{\nu,A}\otimes_{A}\mathbb{C}$, where $\mathbb{C}$ is an ${A}$-module by $f(q)\cdot 1=f(1)$ for $f(q)\in A$. We have $\ov{X}=\bigoplus_{\mu\in P_{\ge 0}}\ov{X}_\mu$ with $\dim_{\mathbb{C}}\ov{X}_\mu={\rm rank}_{A}X_{\mu,{A}}$.

Let $U(\gl(m|n))$ be the universal enveloping algebra of $\gl(m|n)$.
Then $\ov{X}$ becomes a $U(\gl(m|n))$-module with respect to the actions of $\ov{E_i}$, $\ov{F_i}$ and $\ov{H_a}$, the $\mathbb{C}$-linear endomorphisms on $\ov{X}$ induced from $E_i$, $F_i$ and $H_a$ for $i\in I\setminus\{0\}$ and $a\in \I$ (cf.~\cite[Section 2.6]{K14}). By \eqref{eq:extension of two poly}, we have an exact sequence of $U(\mf{gl}(m|n))$-modules
\begin{equation}\label{eq:extension of two poly at q=1}
 0 \longrightarrow \ov{U} \longrightarrow \ov{V} \longrightarrow \ov{W} \longrightarrow 0,
\end{equation}
where $\ov{U}$ and $\ov{W}$ are irreducible highest weight modules with highest weights $\La_\la$ and $\La_\mu$, respectively.
By \cite[Cor 3.1]{CK}, the short exact sequence \eqref{eq:extension of two poly at q=1} splits, and hence so does \eqref{eq:extension of two poly}.

Now, we use induction on the length of composition series to show that $V\in \mc{O}_{\ge 0,d}$ is a direct sum of $V_{\mf l}(\la)$'s for $\la\in \cP_{m|n}$. This completes the proof.
\qed\smallskip

Let $V$ be a $\U(\mf{l})$-module which has a weight space decomposition in $P$, and let $u\in V$ be a weight vector.
For $i\in I\setminus\{0\}$, we define $\tilde{e}_i u$ and $\tilde{f}_i u$ as follows:\smallskip

\begin{itemize}
 \item[(1)] Suppose that $i\neq m$. We have $u=\sum_{k \geq 0} f_i^{(k)}u_k$, where $f_i^{(k)}=f_i^k/[k]_i!$ and $e_iu_k=0$ for $k\geq 0$. We define
\begin{gather} 
\te_iu=\sum_{k\geq1}f_i^{(k-1)}u_k,\quad \tf_iu=\sum_{k\geq0}f_i^{(k+1)}u_k \quad (i<m),\label{eq:crystal operator <m}\\
\tilde{e}_iu=\sum_{k\geq1}q_i^{-l_k+2k-1}f_i^{(k-1)}u_k,\quad 
\tilde{f}_iu=\sum_{k\geq 0}q_i^{l_k-2k-1}f_i^{(k+1)}u_k\quad (i>m),\label{eq:crystal operator >m}
\end{gather}
where $l_k= - ({\rm wt}(u_k)|\alpha_i)/r_{\mf{g}}$. Here we follow the convention in \cite[Section 4.2]{JKU}.

\item[(2)] Suppose that $i=m$. We define
\begin{equation}\label{eq:crystal operator m}
\tilde{e}_m u =\eta(f_m) u =q_m^{-1}k_me_m u,\quad \tilde{f}_m u=f_m u,
\end{equation}
where $\eta$ is the anti-involution on $\U(\mf{g})$ defined by
$\eta(k_\mu)=k_\mu$, $\eta(e_i)=q_if_ik^{-1}_i$, and $\eta(f_i)=q^{-1}_ik_ie_i$ for $\mu \in  P$ and $i\in I$.
\end{itemize}

Let $A_0$ be the subring of $f(q)\in \Bbbk$ regular at $q=0$. 
We call a pair $(L,B)$ a {\em crystal base of $V$} if it satisfies the following conditions:
\begin{itemize}
\item[(1)] $L$ is an $A_0$-lattice of $V$ and $L=\bigoplus_{\mu \in P}L_\mu$, where $L_\mu=L \cap V_\mu$,

\item[(2)] $B$ is a signed basis of $L/qL$, that is $B=\mathbf{B}\cup -\mathbf{B}$ where $\mathbf{B}$ is a $\Q$-basis of $L/qL$,

\item[(3)] $B= \bigsqcup_{\mu \in P} B_\mu$ where $B_\mu \subset (L/qL)_\mu$,

\item[(4)] $\tilde{e}_i L \subset L, \tilde{f}_i L \subset L$ and $\tilde{e}_iB \subset B \cup \{0\}, \tilde{f}_iB \subset B \cup \{0\}$ for $i \in I \setminus \{0\}$,

\item[(5)] $\tilde{f}_ib=b'$ if and only if $\tilde{e}_ib'=\pm b$ for $i\in I \setminus \{0\}$ and $b, b' \in B$.
\end{itemize}
The set $B/\{\pm 1\}$ has a colored oriented graph structure, where $b \stackrel{i}{\rightarrow} b'$ if and only if $\tf_ib=b'$ for $b,b'\in B/\{\pm 1\}$ and $i\in I\setminus\{0\}$. We call $B$ an {\em $\mf{l}$-crystal} or simply {\em crystal of $V$}. We identify $B$ with $B/\{\pm 1\}$ when we consider $B$ as a crystal, if there is no confusion.

For $\la\in\cP_{m|n}$, let
\begin{equation}\label{eq:crystal base of poly}
\begin{split}
\ms{L}(V_{\mf{l}}(\la))&=\sum_{r\geq 0,\, i_1,\ldots,i_r\in I \setminus \{0\}}A_0 \tilde{x}_{i_1}\cdots\tilde{x}_{i_r}v_\la, \\
\ms{B}(V_{\mf{l}}(\la))&=\{\,\pm \,{\tilde{x}_{i_1}\cdots\tilde{x}_{i_r}v_\la}\!\!\! \pmod{q \ms{L}(\la)}\,|\,r\geq 0, i_1,\ldots,i_r\in I\setminus\{0\}\,\}\setminus\{0\},
\end{split}
\end{equation}
where $v_\la$ is a highest weight vector in $V_{\mf l}(\la)$ and $x=e, f$ for each $i_k$.

\begin{thm}[{\cite[Theorem 5.1]{BKK}}]\label{thm:crystal base of irr poly}
The pair $(\ms{L}(V_{\mf{l}}(\la)),\ms{B}(V_{\mf{l}}(\la)))$ is a crystal base of $V_{\mf{l}}(\la)$.
\end{thm}

\begin{prop}[{\cite[Proposition 2.8]{BKK}}]\label{prop:tensor product thm}
 Let $V_k\in \cO_{\geq 0}$ be given with a crystal base $(L_k,B_k)$ for $k=1,2$.
Then $(L_1\otimes L_2, B_1\otimes B_2)$ is a crystal base of $V_1\otimes V_2$, where $\te_i$ and $\tf_i$ $(i \in I\setminus\{0\})$ act on $B_1\otimes B_2$ as follows:

{\allowdisplaybreaks
\begin{equation*}\label{eq:tensor product rule for even -}
\begin{split}
&\te_i(b_1\otimes b_2)= \begin{cases}
\te_i b_1 \otimes b_2 & \text{if $\varphi_i(b_1)\geq\varepsilon_i(b_2)$}, \\ 
 b_1 \otimes \te_ib_2 & \text{if $\varphi_i(b_1)<\varepsilon_i(b_2)$},\\
\end{cases}
\\
&\tf_i(b_1\otimes b_2)=
\begin{cases}
\tf_ib_1 \otimes  b_2 & \text{if $\varphi_i(b_1)>\varepsilon_i(b_2)$}, \\
b_1 \otimes \tf_i  b_2 & \text{if $\varphi_i(b_1)\leq\varepsilon_i(b_2)$}, 
\end{cases}
\end{split}\quad\quad (i<m)\quad\quad 
\end{equation*}}

{\allowdisplaybreaks
\begin{equation*}\label{eq:tensor product rule for even +}
\begin{split}
&\te_i(b_1\otimes b_2)= \begin{cases}
 b_1 \otimes \te_ib_2 & \text{if $\varphi_i(b_2)\geq\varepsilon_i(b_1)$}, \\ 
\te_ib_1 \otimes  b_2 & \text{if $\varphi_i(b_2)<\varepsilon_i(b_1)$},\\
\end{cases}
\\
&\tf_i(b_1\otimes b_2)=
\begin{cases}
 b_1 \otimes \tf_ib_2 & \text{if $\varphi_i(b_2)>\varepsilon_i(b_1)$}, \\
\tf_ib_1 \otimes b_2 & \text{if $\varphi_i(b_2)\leq\varepsilon_i(b_1)$}, 
\end{cases}
\end{split}\quad\quad (i>m)\quad
\end{equation*}}

\begin{equation*}\label{eq:tensor product rule for odd +}
\begin{split}
\te_m(b_1\otimes b_2)=&
\begin{cases}
\te_m  b_1\otimes b_2 & \text{if $({\rm wt}(b_1)|\alpha_m)\neq 0$}, \\ 
 b_1\otimes  \te_m b_2 & \text{if $({\rm wt}(b_1)|\alpha_m)= 0$},
\end{cases}
\\
\tf_m(b_1\otimes b_2)=&
\begin{cases}
 \tf_m b_1\otimes b_2 & \text{if $({\rm wt}(b_1)|\alpha_m)\neq 0$}, \\ 
b_1\otimes  \tf_m b_2 & \text{if $({\rm wt}(b_1)|\alpha_m)= 0$},
\end{cases}
\end{split}\quad\quad\quad\quad\quad\quad 
\end{equation*}
\noindent up to multiplication by $\{\pm 1\}$, where $\varepsilon_i(b)=\max\{k\geq 0 \,|\ \te_i^k b  \neq 0 \}$ and $\varphi_i(b)=\max\{k\geq 0 \,|\ \tf_i^k b \neq 0 \}$ for $b\in B_1, B_2$.
\end{prop}

\begin{rem}{\rm
(1) The map $\tau$ in \eqref{eq:iso tau standard} is not an isomorphisms of Hopf algebras, but the above formula also holds for the crystal bases of $\U(\mf{l})$-modules with respect to \eqref{eq:comult-1}.

(2) The proof of Theorem \ref{thm:crystal base of irr poly} in case of $\U(\g)$ can be obtained in almost the same way as in \cite{BKK} (cf.~\cite[Theorem 4.12]{KO}).
} 
\end{rem}

Let $SST_{m|n}(\la)$ be the set of all $(m|n)$-hook semistandard tableaux of shape $\la$ with letters in $\I$ (cf.~\cite[Definition 4.1]{BKK}). Let $H_\la$ be the unique element in $SST_{m|n}(\la)$ with ${\rm wt}(H_\la)=\La_\la$, which is given by filling the $i$-th row with $i$ for $1\le i\le m$, and then the $j$-th column (of the remaining diagram) with $m+j$ for $1\le j\le n$.

Let $B$ be the crystal $\ms{B}(V_{\mf l}(\de_1))$ which is given by
\begin{equation*}
1\ \stackrel{^{1}}{\longrightarrow}\ 2\ \stackrel{^{2}}{\longrightarrow}
\cdots\stackrel{^{m+n-1}}{\longrightarrow}\  m+n 
\end{equation*}
(see \cite[Section 3.2]{BKK}). For $T\in SST_{m|n}(\la)$, let $w_1 \cdots  w_\ell$ be its column word, that is, the word given by reading the entries in $T$ column by column from right to left and from top to bottom in each column.
Then there is a well-defined $I\setminus\{0\}$-colored oriented graph structure on $SST_{m|n}(\la)$ induced from the map sending $T$ to $w_1\ot \dots \ot w_\ell\in B^{\ot \ell}$. 

\begin{thm}[{\cite[Theorems 4.8 and 5.1]{BKK}}]\label{thm:Blambda cnn}
As an $I\setminus\{0\}$-colored oriented graph, $SST_{m|n}(\la)$ is connected and isomorphic to $\ms{B}(V_{\mf l}(\la))$. 
\end{thm}

The element $H_\la$ is called a genuine highest weight vector in \cite{BKK}.
Conversely, for $w=w_1\ot \dots \ot w_\ell\in B^{\ot \ell}$, there exists a unique $P(w) \in SST_{m|n}(\la)$ for some $\la\in \cP_{m|n}$ such that the connected component of $w$ in $B^{\ot \ell}$ is isomorphic to $SST_{m|n}(\la)$ under $w \mapsto T$. Indeed, $P(w)$ is determined by a column insertion (see \cite[Section 4.5]{BKK} for more details).

\subsection{Crystal base of $\mc{N}$}
In this subsection, we construct a crystal base of $\mc{N}$ as a $\U(\mf{l})$-module with respect to \eqref{eq: quantum adjoint}.
We first observe the following.
\begin{lem}\label{lem:N poly}
We have $\mc{N}\in\mc{O}_{\ge 0}$.
\end{lem}
\pf For $\beta\in \Phi^+(\mf{u})$, we have $k_\mu\cdot \bff_\beta =\bq(-\beta,\mu) \bff_\beta$, that is, ${\rm wt}(\bff_\beta)=-\beta\in P_{\ge 0}$. Hence the weight of a product of $\bff_\beta$'s belongs to $P_{\ge 0}$. Also it is easy to see that $\dim \mc{N}_\mu <\infty$ for $\mu\in P_{\ge 0}$.
\qed
\smallskip

We will construct a crystal base of $\mc{N}$ using a PBW type basis in Proposition \ref{prop:PBW basis}. 
In order to have a crystal base with respect to $\tf_i$ for $(\alpha_i|\alpha_i)<0$ (an upper crystal base in \cite{Kas91}), we need to normalize each root vector as follows. 
For $\beta\in \Phi^+(\mf{u})$ and $k\in \Z_{\ge 0}$, let
\begin{equation*}\label{eq:divided power of root vector}
\bfF_\beta^{(k)} = 
\begin{cases}
 \bff_\beta^{(k)} & \text{if $(\beta|\beta)\ge 0$},\\
 q^{-\frac{k(k-1)}{2}r_{{\mf x}}}\bff_\beta^{(k)} & \text{if $(\beta|\beta)< 0$ and $\beta\ne(i,i)$},\\
 q^{-\frac{k(k-1)}{2}}\bff_\beta^{(k)} & \text{if $\mf{g}=\mf{b}_{m|n}$, $(\beta|\beta)< 0$, and $\beta=(i,i)$},\\
 q^{-k^2}\bff_\beta^{(k)} & \text{if $\mf{g}=\mf{d}_{m|n}$, $(\beta|\beta)< 0$, and $\beta=(i,i)$}.
\end{cases}
\end{equation*}

Let $i\in I\setminus\{0\}$ be given. Let $\U(\mf{l})_i$ be the subalgebra of $\U(\mf{l})$ generated by $e_i$, $f_i$, and $k_{\de_{a}}^{\pm 1}$ ($a=i,i+1$). One may consider a crystal base of a $\U(\mf{l})_i$-module with respect to \eqref{eq:crystal operator <m}--\eqref{eq:crystal operator m}.

For $1\le k<i$ and $c\in\Z_{\ge 0}$, put
\begin{equation*}\label{eq:N i-comp}
\begin{split}
 \mc{N}_k(i;c) &= \bigoplus_{a+b=c}\Bbbk\, \bfF_{(k,i)}^{(a)}\bfF_{(k,i+1)}^{(b)},\\
  \ms{L}(\mc{N}_k(i;c)) &= \bigoplus_{a+b=c}A_0\, \bfF_{(k,i)}^{(a)}\bfF_{(k,i+1)}^{(b)},\\
 \ms{B}(\mc{N}_k(i;c)) &= \left\{\, \pm\bfF_{(k,i)}^{(a)}\bfF_{(k,i+1)}^{(b)}\!\!\!\pmod{q\ms{L}(\mc{N}_k(i;c))}\,|\,a+b=c\,\right\}.
\end{split}
\end{equation*}

\begin{lem}\label{lem:i-comp submodule-1}
$\mc{N}_k(i;c)$ is an irreducible $\U(\mf{l})_i$-submodule of $\mc{N}$, which has a crystal base $\left(\ms{L}(\mc{N}_k(i;c)),\ms{B}(\mc{N}_k(i;c))\right)$.
\end{lem}
\pf It follows from Propositions \ref{prop:commutation rel}, \ref{prop:e dot} and \ref{prop:f dot}.
\qed
\smallskip

Similarly, for $i+1<l\le m+n$ and $c\in\Z_{\ge 0}$, put 
\begin{equation*}
\begin{split}
 \mc{N}^l(i;c) &= \bigoplus_{a+b=c}\Bbbk\, \bfF_{(i,l)}^{(a)}\bfF_{(i+1,l)}^{(b)},\\
 \ms{L}\left(\mc{N}^l(i;c)\right) &= \bigoplus_{a+b=c}A_0\,  \bfF_{(i,l)}^{(a)}\bfF_{(i+1,l)}^{(b)}, \\
 \ms{B}\left(\mc{N}^l(i;c)\right) &= \left\{\,\pm \bfF_{(i,l)}^{(a)}\bfF_{(i+1,l)}^{(b)}\!\!\!\pmod{q\ms{L}(\mc{N}^l(i;c))} \,|\,a+b=c\,\right\}.
\end{split}
\end{equation*}
As in Lemma \ref{lem:i-comp submodule-1}, we can prove
 
\begin{lem}\label{lem:i-comp submodule-2}
$\mc{N}^l(i;c)$ is an irreducible $\U(\mf{l})_i$-submodule of $\mc{N}$, which has a crystal base $\left(\ms{L}(\mc{N}^l(i;c)),\ms{B}(\mc{N}^l(i;c))\right)$.
\end{lem}

\begin{rem}\label{rem:crystal operator on i-comp-1}
{\rm
The crystal operator $\tf_i$ acts on $\ms{B}(\mc{N}_k(i;c))$ and  $\ms{B}(\mc{N}^l(i;c))$ by
\begin{equation*}
\begin{split}
 &\bfF_{(k,i)}^{(a)}\bfF_{(k,i+1)}^{(b)}\
 \stackrel{i}{\longrightarrow}\
 \bfF_{(k,i)}^{(a-1)}\bfF_{(k,i+1)}^{(b+1)},\qquad
 \bfF_{(i,l)}^{(a)}\bfF_{(i+1,l)}^{(b)}\
 \stackrel{i}{\longrightarrow}\
 \bfF_{(i,l)}^{(a-1)}\bfF_{(i+1,l)}^{(b+1)},
\end{split}\quad (a\ge 1),
\end{equation*}
respectively, where we assume $b = 0$ (resp. $a=1$) if $(k,i+1)$ (resp.~$(i,l)$) is an isotropic root.
} 
\end{rem}
\smallskip

Next, for $d\in\Z_{\ge 0}$, put
\begin{equation*}\label{eq:N i-comp}
\begin{split}
 \mc{N}_{\Delta}(i;d)&=\bigoplus_{a+r_{\mf g}b+c=d}\Bbbk \, \bfF_{(i,i)}^{(a)}\bfF_{(i,i+1)}^{(b)}\bfF_{(i+1,i+1)}^{(c)},\\
  \ms{L}\left(\mc{N}_{\Delta}(i;d)\right)&=\bigoplus_{a+r_{\mf g}b+c=d}A_0\, \bfF_{(i,i)}^{(a)}\bfF_{(i,i+1)}^{(b)}\bfF_{(i+1,i+1)}^{(c)},\\
 \ms{B}\left(\mc{N}_{\Delta}(i;d)\right) &= \left\{\, \pm \bfF_{(i,i)}^{(a)}\bfF_{(i,i+1)}^{(b)}\bfF_{(i+1,i+1)}^{(c)}\!\!\!\pmod{q\ms{L}\left(\mc{N}_{\Delta}(i;d)\right)}\,|\,a+r_{\mf g}b+c=d \,\right\},
\end{split}
\end{equation*}
where we assume that $\bfF_{(k,l)}=0$ for $(k,l)\not \in \Phi^+(\mf{u})$.
As in Lemmas \ref{lem:i-comp submodule-1} and \ref{lem:i-comp submodule-2}, we have
\begin{lem}
$\mc{N}_{\Delta}(i;d)$ is a $\U(\mf{l})_i$-submodule of $\mc{N}$.
\end{lem}

Note that $\mc{N}_{\Delta}(i;d)$ is not irreducible in general.

\begin{prop}\label{prop:i-comp submodule-3}
$\left(\ms{L}(\mc{N}_{\Delta}(i;d)),\ms{B}(\mc{N}_{\Delta}(i;d))\right)$ is a crystal base of $\mc{N}_{\Delta}(i;d)$.
\end{prop}
\pf Put 
\begin{equation}\label{eq:S}
 S=\{\,(i,i), (i,i+1), (i+1,i+1)\,\}\cap \Phi^+(\mf{u}),
\end{equation}
and let 
$N=\oplus_{d \ge 0}\, \mc{N}_{\Delta}(i;d)$ and $(L,B)=\oplus_{d \ge 0}\,\left(\ms{L}(\mc{N}_{\Delta}(i;d)),\ms{B}(\mc{N}_{\Delta}(i;d))\right)$. We prove that $(L,B)$ is a crystal base of $N$, which implies the case of $\mc{N}_{\Delta}(i;d)$.

\noindent{\em Case 1}. $S=\{\,(i,i+1)\,\}$. In this case, $N$ is a trivial representation. So it is clear.\smallskip

\noindent{\em Case 2}. $S=\{\,(i,i), (i,i+1)\,\}$ or $\{\,(i,i+1), (i+1,i+1)\,\}$. In this case, we have $i=m$ (cf.~ Example \ref{ex: roots of u}). It is also easy to see that $(L,B)$ is a crystal base of $N$.\smallskip

\noindent{\em Case 3}. $S=\{\,(i,i), (i,i+1), (i+1,i+1)\,\}$ and $\mf{g}=\mf{b}_{m|n}$, $\mf{c}_{m|n}$ with $i<m$.
Let $\texttt{U}=\U(\mf{x}_{2|0})$ ($\mf{x}=\mf{b}, \mf{c}$), and let $\texttt{N}$ be the subalgebra of $\texttt{U}^-$ generated by $\bff_1=\bff_{(1,1)}$, $\bff_2=\bff_{(1,2)}$, $\bff_3=\bff_{(2,2)}$. 

By Proposition \ref{prop:commutation rel}, there exists an isomorphism of $\Bbbk$-algebras $\rho: N \longrightarrow \texttt{N}$ such that $\rho(\bff_{(i,i)})=\bff_1$, $\rho(\bff_{(i,i+1)})=\bff_2$, and $\rho(\bff_{(i+1,i+1)})=\bff_3$. If we identify $\U(\mf{l})_i$ with $\texttt{U}_1=\langle\, e_1, f_1, k_{\de_a}^{\pm 1}\,|\, a=1,2\,\rangle$, both of which are isomorphic to $U_q(\mf{gl}_2)$, then $\rho$ is also an isomorphism of $U_q(\mf{gl}_2)$-modules in the sense of Proposition \ref{prop:poly alg}. 

Note that $\texttt{U}^-\cong \texttt{N}\ot \texttt{U}_1^-$ as a $\Bbbk$-space.
 Let $\pi_1: \texttt{U}^- \longrightarrow \texttt{P}$ be the natural projection, where $\texttt{P}=\texttt{U}^-/\texttt{U}^-\texttt{U}^-_1$, and let $\pi_2 : \texttt{P} \longrightarrow \texttt{N}$ be the $\Bbbk$-linear isomorphism induced from \eqref{eq:parabolic decomp}. 
 Put $\pi=\pi_2\circ\pi_1 : \texttt{U}^- \longrightarrow \texttt{N}$. Note that $\pi$ commutes with the action of $f_1$, where we regard $\texttt{P}$ as $P(0)$ in Section \ref{subsec:parabolic Ver} and apply Proposition \ref{prop:par Verma and N are iso}.

 Let $(\texttt{L}(\infty),\texttt{B}(\infty))$ be the crystal base of $\texttt{U}^-$ \cite{Kas91}. Let $(\texttt{L},\texttt{B})$ be the image of $(L,B)$ under $\rho$. Indeed, $\texttt{L}(\infty)$ is equal to $A_0$-span of $\bff_1^{(c_1)}\bff_2^{(c_2)}\bff_3^{(c_3)}f_1^{(c_4)}$ (see \cite{Lu90,Lu90-2,S94}), while $\texttt{L}$ is the $A_0$-span of $\bff_1^{(c_1)}\bff_2^{(c_2)}\bff_3^{(c_3)}$, so that $\pi$ maps $(\texttt{L}(\infty),\texttt{B}(\infty))$ onto $(\texttt{L},\texttt{B})$. 
Recall that the crystal operator $\tf_1$ on $\texttt{U}^-$ is given by 
$\tilde{f}_1u=\sum_{k\geq 0}f_1^{(k+1)}u_k$ for $u=\sum_{k\geq 0}f_1^{(k)}u_k$ with $e'_1 (u_k)=0$. 
On the other hand, the crystal operator $\tf_1$ on $\texttt{N}$ is given by \eqref{eq:crystal operator <m} with respect to \eqref{eq: quantum adjoint}.

For ${\bf c}=(c_1,c_2,c_3)\in \Z_{\ge 0}^3$, let $\bff^{({\bf c})}=\bff_1^{(c_1)}\bff_2^{(c_2)}\bff_3^{(c_3)}$. Regarding $\bff^{({\bf c})}\in \texttt{B}(\infty)$,
let $\texttt{B}_0=\{\,\bff^{({\bf c})}\,|\, \te_1 \bff^{({\bf c})}=0\,\}=\{\,\bff^{({\bf c})}\,|\,{\bf c}=(a,0,c)\ (a\ge c)\,\}$ (see Remark \ref{rem:crystal operator on i-comp-2}).
For each ${\bf c}_0 = (a, 0, c)$ with $a\ge c$, there exists $E_{{\bf c}_0}\in \texttt{L}$ such that $e_1\cdot E_{{\bf c}_0}=0$ and $E_{{\bf c}_0}\equiv \bff^{({\bf c}_0)}\pmod{q\texttt{L}}$. Indeed, by using Proposition \ref{prop:e dot}, one can find that
\begin{equation*}
 E_{{\bf c}_0}=E_{a,c}=
 \begin{cases}
 	\bff_1^{(a)}\bff_3^{(c)} + \sum_{r=1}^{c} X_r \bff_1^{(a-r)}\bff_2^{(r)}\bff_3^{(c-r)} & \text{if $\mf{g} = \mf{b}_{m|n}$,} \\
	\bff_1^{(a)}\bff_3^{(c)} + \sum_{r=1}^{c} Y_r \bff_1^{(a-r)}\bff_2^{(2r)}\bff_3^{(c-r)} & \text{if $\mf{g} = \mf{c}_{m|n}$,}
 \end{cases}
\end{equation*}
where
\begin{equation*}
\begin{split}
X_r &= (-1)^r q^{-cr+\frac{r^2+3r}{2}} \frac{(1+q^2) (1+q^4) \cdots (1+q^{2r})}{(1-q^2)^r [a] \cdots [a-r+1]} \in q^{r(a-c+1)} A_0 \\
Y_r &= (-1)^r q^{-2cr+2r^2+r} \frac{[2r-1][2r-3] \cdots [1]}{[2a][2a-2] \cdots [2a-2r+2]} \in q^{2r(a-c+1)} A_0
\end{split} \quad (1 \le r \le c),
\end{equation*}
(cf.~Appendix \ref{subsec:b i>m}).

Regarding $E_{{\bf c}_0}$ as an element of $\texttt{L}(\infty)$, 
we have $e'_1(E_{{\bf c}_0})=0$ by Lemma \ref{lem:e-derivation and adjoint}. 
Since $(\texttt{L}(\infty),\texttt{B}(\infty))$ is a crystal base of $\texttt{U}^-$, we have $\tf_1^s E_{{\bf c}_0}= f_1^{(s)}E_{{\bf c}_0}\in \texttt{L}(\infty)$ and $\tf_1^s E_{{\bf c}_0}\equiv \tf_1^s\bff^{({\bf c}_0)}\pmod{q\texttt{L}(\infty)}$ for $s\ge 0$, where $\tf_1^s\bff^{({\bf c}_0)}\equiv \bff^{({\bf c})}\pmod{q\texttt{L}(\infty)}$ for some ${\bf c}$. 
Hence we have 
\begin{equation*}
 \pi(f_1^{(s)}E_{{\bf c}_0}) = f_1^{(s)}\cdot E_{{\bf c}_0} = \tf_1^s E_{{\bf c}_0}\in \texttt{L},\quad 
 \tf_1^s E_{{\bf c}_0}\equiv \bff^{({\bf c})}\pmod{q\texttt{L}}.
\end{equation*}
This implies that $\texttt{L}=\bigoplus_{\bff^{({\bf c}_0)}\in \texttt{B}_0,s\in\Z_{\ge 0}}A_0 \tf_1^s\bff^{({\bf c}_0)}$, and $(\texttt{L},\texttt{B})$ is a crystal base of $\texttt{N}$ as a $\texttt{U}_1$-module.
\smallskip

\noindent{\em Case 4}. $S=\{\,(i,i), (i,i+1), (i+1,i+1)\,\}$,  and $\mf{g}=\mf{b}_{m|n}$ with $i\ge m$ or $\mf{g}=\mf{d}_{m|n}$ with $i> m$. In this case, we may not apply the arguments in {\em Case 3} directly. The proof is given in Appendix \ref{appendix-1}.
\qed
\smallskip

\begin{rem}\label{rem:crystal operator on i-comp-2}
{\rm
The crystal operator $\tf_i$ on $\ms{B}(\mc{N}_\Delta(i;d))$ can be described explicitly as follows. Let $S$ be as in \eqref{eq:S}.

Suppose that $|S|=2$. Then we have $\mf{g} = \mf{c}_{m|n}$ or $\mf{d}_{m|n}$ with $i = m$, where
\begin{equation*}
\begin{cases}
 \bfF_{(m,m)}^{(a)}\
 \stackrel{m}{\longrightarrow}\
 \bfF_{(m,m)}^{(a-1)}\bfF_{(m,m+1)} & \text{for $\mf{g} = \mf{c}_{m|n}$ with $a \ge 1$,} \\
 \bfF_{(m,m+1)}\bfF_{(m+1,m+1)}^{(b)}\
 \stackrel{m}{\longrightarrow}\
 \bfF_{(m+1,m+1)}^{(b+1)} & \text{for $\mf{g} = \mf{d}_{m|n}$,}
\end{cases}
\end{equation*} 
for $a,b\in \Z_{\ge 0}$.

Suppose that $|S|=3$. We identify $\bfF_{(i,i)}^{(a)}\bfF_{(i,i+1)}^{(b)}\bfF_{(i+1,i+1)}^{(c)}$ with $(a,b,c)$ for simplicity.

 If $i<m$, then we have $\g=\mf{b}_{m|n}, \mf{c}_{m|n}$. The $\mf{sl}_2$-crystal structure on $\ms{B}(\mc{N}_\Delta(i;d))$ can be deduced from the crystal of the negative half of the quantized enveloping algebra of type $B_2$ and $C_2$ by applying the method of folding \cite{Kas96} to that of type $A_3$ (see for example \cite[Proposition 4.5]{K18}). 

If $i\ge m$, then we have $\g=\mf{b}_{m|n}$ with $i\ge m$ or $\g=\mf{c}_{m|n}$ with $i<m$. The formulas of $\tf_i$ are obtained in the proof of Proposition \ref{prop:i-comp submodule-3} in Appendix \ref{appendix-1}. We should remark that the formula of $\tf_i$ for $\mf{b}_{m|n}$ (resp.~$\mf{d}_{m|n}$) with $i>m$ is not equal to that of $\tf_i$  for $\mf{b}_{m|n}$ (resp.~$\mf{c}_{m|n}$) with $i<m$.

It is summarized as follows.\smallskip
 
\noindent{\em Case 1}. $\mf{g}=\mf{b}_{m|n}$.
If $i<m$, then we have
\begin{equation*}
\tf_i (a,b,c)=
\begin{cases}
(a-2,b+1,c) & \text{if $a-c\ge 2$},\\
(a-1,b,c+1) & \text{if $a-c=1$},\\
(a,b-1,c+2) & \text{if $a\le c$ and $b\ge 1$},\\
0 & \text{otherwise}.
\end{cases}
\end{equation*}
If $i>m$, then we have
\begin{equation*}
\tf_i (a,b,c)=
\begin{cases}
(a-1,b,c+1) & \text{if $a\ge 1$},\\
0 & \text{otherwise}.\\ 
\end{cases}
\end{equation*}
If $i=m$, then we have
\begin{equation*}
\tf_i (a,b,c)=
\begin{cases}
(a-2,1,c) & \text{if $a\ge 2$ and $b=0$},\\
(0,0,c+1) & \text{if $a= 1$ and $b=0$},\\
0 & \text{otherwise}.\\ 
\end{cases}
\end{equation*}

\smallskip

\noindent{\em Case 2}. $\mf{g}=\mf{c}_{m|n}$. We have $i<m$ and 
\begin{equation*}
\tf_i (a,b,c)=
\begin{cases}
(a-1,b+1,c) & \text{if $a-c\ge 1$},\\
(a,b-1,c+1) & \text{if $a\le c$ and $b\ge 1$},\\
0 & \text{otherwise}.
\end{cases}
\end{equation*}
\smallskip

\noindent{\em Case 3}. $\mf{g}=\mf{d}_{m|n}$. We have $i>m$ and
\begin{equation*}
\tf_i (a,b,c)=
\begin{cases}
(a-1,b+1,c) & \text{if $a\ge 1$ and $b$ is even},\\
(a,b-1,c+1) & \text{if $b\ge 1$ and odd},\\
0 & \text{otherwise}.
\end{cases}
\end{equation*}
}
\end{rem}
\smallskip

Let
\begin{equation*}
	M(\mf{u}) = 
	\left\{\, (c_\beta)_{\beta\in \Phi^+(\mf{u})}\in \Z_{\ge 0}^{N} \, | \,
	\text{$c_{\beta}\in\Z_{\ge 0}$ $(\beta\in \Phi^+_{\ov 0}\cup \Phi^+_{\rm{n\text{-}iso}})$,\quad $c_{\beta}=0,1$  $(\beta\in \Phi^+_{{\rm iso}})$} 
	\,\right\},
\end{equation*}
where $N=|\Phi^+(\mf{u})|$.
For ${\bf c}=(c_\beta)_{\beta\in \Phi^+(\mf{u})}\in M(\mf{u})$, let 
\begin{equation*}
\bfF^{(\bf c)} = \prod^{\rightarrow}_{\beta\in\Phi^+(\mf{u})}\bfF_\beta^{(c_\beta)}
\end{equation*}
(cf.~Remark \ref{rem:commutation relation}).
Then 
$B(\mc{N})=\{\,\bfF^{(\bf c)}\,|\, {\bf c}\in M(\mf{u}) \,\}$ 
is a $\Bbbk$-basis of $\mc{N}$.
Let 
\begin{equation}\label{eq:crystal base of N}
\begin{split}
 \ms{L}\left(\mc{N}\right) &= \bigoplus_{{\bf c}\in M(\mf{u})} A_0 \bfF^{(\bf c)},\\
 \ms{B}\left(\mc{N}\right) &= \left\{\,\pm  \bfF^{(\bf c)}\!\!\!\pmod{q\ms{L}\left(\mc{N}\right)}\,|\,{\bf c}\in M(\mf{u})\,\right\}.
\end{split} 
\end{equation}

\begin{thm}\label{thm:crystal base of N}
The pair $(\ms{L}\left(\mc{N}\right),\ms{B}\left(\mc{N}\right))$ is a crystal base of $\mc{N}$ as a $\U(\mf{l})$-module.
\end{thm}
\pf Suppose that $i\in I\setminus\{0\}$ is given.
Let
\begin{equation*}\label{eq:N i-comp-2}
\mc{N}(i)=
\mc{N}_{{}^{\searrow}}(i)\mc{N}_{\Delta}(i)\mc{N}_{{}^{\nearrow}}(i),
\end{equation*}
where
\begin{equation}\label{eq:notations for N(i)'s}
\begin{split}
 \mc{N}_{{}^{\searrow}}(i) &=\bigoplus_{d_1,\dots,d_{i-1}\in\Z_{\ge 0}} \mc{N}_1(i;d_1)\mc{N}_2(i;d_2)\dots \mc{N}_{i-1}(i;d_{i-1}),\\
 \mc{N}_{{}^{\nearrow}}(i) &=\bigoplus_{d_{i+2},\dots,d_{m+n}\in\Z_{\ge 0}} \mc{N}^{i+2}(i;d_{i+2})\mc{N}^{i+3}(i;d_{i+3})\dots \mc{N}^{m+n}(i;d_{m+n}),\\
 \mc{N}_{\Delta}(i)&=\bigoplus_{d\in\Z_{\ge 0}} \mc{N}_{\Delta}(i;d).
\end{split}
\end{equation}
Here we assume that $\mc{N}_{\ast}(i)=\Bbbk$ ($\ast={\searrow},{\nearrow}$) whenever it is not defined.
It is a $\U(\mf{l})_i$-submodule of $\mc{N}$ by Proposition \ref{prop:poly alg}.
Let $S_i$ be a subset of $\beta \in \Phi^+(\mf{u})$ such that $\bfF_\beta$ appears as a factor of a monomial in $B(\mc{N})\cap \mc{N}(i)$.

Let ${\bf c}=(c_\beta)\in M(\mf{u})$ be given, and let ${\bf c}^\circ=(c_\beta^\circ)\in M(\mf{u})$ be such that $c^\circ_\beta=c_\beta$ for $\beta\in S_i$ and $0$ otherwise.
Consider a $\Bbbk$-linear map
\begin{equation}\label{eq:i_c}
 \xymatrixcolsep{2pc}\xymatrixrowsep{0pc}\xymatrix{
 \iota_{{\bf c}} : \mc{N}(i) \ \ar@{->}[r] &\ \mc{N} \\
 \quad  \bfF^{({\bf c}')} \ \ar@{|->}[r] &\  \bfF^{({\bf c}'')}}, 
\end{equation}
where for ${\bf c}'=(c'_\beta)_{\beta\in S_i}$,  ${\bf c}''=(c''_\beta) \in M(\mf{u})$ is given by $c''_\beta=c'_\beta$ for $\beta\in S_i$ and $c''_\beta=c_\beta$ otherwise. We have $\iota_{{\bf c}}(\bfF^{({\bf c}^\circ)})=\bfF^{({\bf c})}$ by definition.
Since $e_i\cdot \bfF_\beta=f_i\cdot \bfF_\beta=0$ and $k_\mu\cdot \bfF_\beta=\bfF_\beta$ ($\mu=\pm\de_i,\pm\de_{i+1}$) for $\beta\not\in S_i$ by Propositions \ref{prop:e dot} and  \ref{prop:f dot}, it follows from  Proposition \ref{prop:poly alg} that $\iota_{{\bf c}}$ is a $\U(\mf{l})_i$-linear map.

On the other hand, since $\mc{N}_k(i;d_k)$, $\mc{N}^l(i;d_l)$, and $\mc{N}_\Delta(i;d)$ have crystal bases, it follows from Propositions \ref{prop:poly alg} and \ref{prop:tensor product thm} that $\mc{N}(i)$ has a crystal base $(\ms{L}(\mc{N}(i)),\ms{B}(\mc{N}(i)))$, where 
\begin{equation*}
 \ms{L}(\mc{N}(i))=\ms{L}(\mc{N}) \cap \mc{N}(i),\quad \ms{B}(\mc{N}(i))=\ms{B}(\mc{N})\cap (\ms{L}(\mc{N}(i))/q\ms{L}(\mc{N}(i))).
\end{equation*}

Finally, since $\iota_{{\bf c}}$ is $\U(\mf{l})_i$-linear and $\iota_{{\bf c}}$ maps $(\ms{L}(\mc{N}(i)),\ms{B}(\mc{N}(i)))$ into $(\ms{L}(\mc{N}),\ms{B}(\mc{N}))$, we conclude that $\tilde{x}_i \bfF^{({\bf c})} \in \ms{L}(\mc{N})$ and $\tilde{x}_i \bfF^{({\bf c})} \in \ms{B}\left(\mc{N}\right) \cup \{ 0 \} \text{ (mod $q\ms{L}(\mc{N})$)}$, where $x = e, f$. This implies that $(\ms{L}(\mc{N}),\ms{B}(\mc{N}))$ is a crystal base of $\mc{N}$.
\qed
\medskip

\begin{rem}
{\em 
The $\mf{l}$-crystal structure on $\ms{B}(\mc{N})$ can be described explicitly by Proposition \ref{prop:tensor product thm} and Remarks \ref{rem:crystal operator on i-comp-1} and \ref{rem:crystal operator on i-comp-2}.
}
\end{rem}

\begin{ex} \label{ex:crystal structure on N}
{\rm
Let us give examples of $\tf_i {\bf c}$ for ${\bf c}\in M(\mf{u})$, where we identify ${\bf c}$ with $\bfF^{({\bf c})}\in \ms{B}(\mc{N})$.
\smallskip

\noindent
(1) Suppose that $\mf{g}=\mf{b}_{2|3}$.
Let ${\bf c}=\left( c_{(i,j)} \right)_{(i,j)\in \Phi^+(\mf{u})}\in M(\mf{u})$ be given by 

\begin{equation*}
	{\bf c} = 
	\begin{tikzpicture}[baseline=(current  bounding  box.center), every node/.style={scale=0.85}, scale=0.95]
		\node (glw_00f) at (0,0) {$c_{(1,1)}$};
		\node (glw_01m) at (0.7,0.7) {$c_{(1,2)}$};
		\node (glw_02m) at (1.4,1.4) {$c_{(1,3)}$};
		\node (glw_03m) at (2.1,2.1) {$c_{(1,4)}$};
		\node (glw_04m) at (2.8,2.8) {$c_{(1,5)}$};
		\node (glw_01f) at (1.4,0) {$c_{(2,2)}$};
		\node (glw_12m) at (2.1,0.7) {$c_{(2,3)}$};
		\node (glw_13m) at (2.8,1.4) {$c_{(2,4)}$};
		\node (glw_14m) at (3.5,2.1) {$c_{(2,5)}$};
		\node (glw_02f) at (2.8,0) {$c_{(3,3)}$};
		\node (glw_23m) at (3.5,0.7) {$c_{(3,4)}$};		
		\node (glw_24m) at (4.2,1.4) {$c_{(3,5)}$};
		\node (glw_03f) at (4.2,0) {$c_{(4,4)}$};
		\node (glw_34m) at (4.9,0.7) {$c_{(4,5)}$};
		\node (glw_04m) at (5.6,0) {$c_{(5,5)}$};
	\end{tikzpicture}
	 = 
	\begin{tikzpicture}[baseline=(current  bounding  box.center), every node/.style={scale=0.85}, scale=0.95]
		\node (glw_00f) at (0,0) {$2$};
		\node (glw_01m) at (0.7,0.7) {$2$};
		\node (glw_02m) at (1.4,1.4) {$1$};
		\node (glw_03m) at (2.1,2.1) {$0$};
		\node (glw_04m) at (2.8,2.8) {$1$};
		\node (glw_01f) at (1.4,0) {$1$};
		\node (glw_12m) at (2.1,0.7) {$0$};
		\node (glw_13m) at (2.8,1.4) {$1$};
		\node (glw_14m) at (3.5,2.1) {$0$};
		\node (glw_02f) at (2.8,0) {$0$};
		\node (glw_23m) at (3.5,0.7) {$2$};		
		\node (glw_24m) at (4.2,1.4) {$1$};
		\node (glw_03f) at (4.2,0) {$1$};
		\node (glw_34m) at (4.9,0.7) {$0$};
		\node (glw_04m) at (5.6,0) {$1$};
	\end{tikzpicture}
\end{equation*} (cf.~Example \ref{ex: roots of u}).\smallskip

	First, let us compute $\tf_1 {\bf c}$.
	According to the proof of Theorem \ref{thm:crystal base of N} (cf.~\eqref{eq:i_c}), it is enough to consider the bold-faced entries, denoted by ${\bf c}^\circ$, while the other entries are given in gray:
	
	\begin{equation*}
	\qquad
	{\bf c} = 
		\begin{tikzpicture}[baseline=(current  bounding  box.center), every node/.style={scale=0.85}, scale=0.95]
		\node (glw_00f) at (0,0) {\bf 2};
		\node (glw_01m) at (0.7,0.7) {\bf 2};
		\node (glw_02m) at (1.4,1.4) {\bf 1};
		\node (glw_03m) at (2.1,2.1) {\bf 0};
		\node (glw_04m) at (2.8,2.8) {\bf 1};
		\node (glw_01f) at (1.4,0) {\bf 1};
		\node (glw_12m) at (2.1,0.7) {\bf 0};
		\node (glw_13m) at (2.8,1.4) {\bf 1};
		\node (glw_14m) at (3.5,2.1) {\bf 0};
		\node (glw_02f) at (2.8,0) {\color{gray}{$0$}};
		\node (glw_23m) at (3.5,0.7) {\color{gray}{$2$}};		
		\node (glw_24m) at (4.2,1.4) {\color{gray}{$1$}};
		\node (glw_03f) at (4.2,0) {\color{gray}{$1$}};
		\node (glw_34m) at (4.9,0.7) {\color{gray}{$0$}};
		\node (glw_04m) at (5.6,0) {\color{gray}{$1$}};
	\end{tikzpicture}
	\quad
	{\bf c}^\circ = 
	\begin{tikzpicture}[baseline=(current  bounding  box.center), every node/.style={scale=0.85}, scale=0.95]
		\node (glw_00f) at (0,0) {\bf 2};
		\node (glw_01m) at (0.7,0.7) {\bf 2};
		\node (glw_02m) at (1.4,1.4) {\bf 1};
		\node (glw_03m) at (2.1,2.1) {\bf 0};
		\node (glw_04m) at (2.8,2.8) {\bf 1};
		\node (glw_01f) at (1.4,0) {\bf 1};
		\node (glw_12m) at (2.1,0.7) {\bf 0};
		\node (glw_13m) at (2.8,1.4) {\bf 1};
		\node (glw_14m) at (3.5,2.1) {\bf 0};
		\node (glw_02f) at (2.8,0) {};
		\node (glw_23m) at (3.5,0.7) {};		
		\node (glw_24m) at (4.2,1.4) {};
		\node (glw_03f) at (4.2,0) {};
		\node (glw_34m) at (4.9,0.7) {};
		\node (glw_04m) at (5.6,0) {};
	\end{tikzpicture}
	\end{equation*}\smallskip
	
\noindent In fact, ${\bf c}^\circ$ can be understood as an element in a tensor product of crystals (cf.~\eqref{eq:notations for N(i)'s}) 
	$$
		{\bf c}^\circ 
		= \underbrace{\left( c_{(1,1)}, c_{(1,2)}, c_{(2,2)} \right)}_{{\bf c}^\circ_{1,1}} \ot 
		  \underbrace{\left( c_{(1,3)}, c_{(2,3)} \right)}_{{\bf c}^\circ_{1,2}} \ot
		  \underbrace{\left( c_{(1,4)}, c_{(2,4)} \right)}_{{\bf c}^\circ_{1,3}} \ot
		  \underbrace{\left( c_{(1,5)}, c_{(2,5)} \right)}_{{\bf c}^\circ_{1,4}},
	$$
	where   ${\bf c}^\circ_{1,1}\in \ms{B}\left(\mc{N}_{\Delta}(1;5)\right)$ and ${\bf c}^\circ_{1,l}\in \ms{B}\left(\mc{N}^l(1;1)\right)$ ($l=2,3,4$).
	
	Now we apply the {\em signature rule} in the tensor product of crystals (cf.~\cite{BS}).
	Let 
\begin{equation*}
 \sigma({\bf c}^\circ) = 
		( 
		\underbrace{-\dots -}_{\varepsilon_1({\bf c}^\circ_{1,1})}, \underbrace{+\dots +}_{\varphi_1({\bf c}^\circ_{1,1})}, 
		\underbrace{-\dots -}_{\varepsilon_1({\bf c}^\circ_{1,2})}, \underbrace{+\dots +}_{\varphi_1({\bf c}^\circ_{1,2})}, 
		\underbrace{-\dots -}_{\varepsilon_1({\bf c}^\circ_{1,3})}, \underbrace{+\dots +}_{\varphi_1({\bf c}^\circ_{1,3})}, 
		\underbrace{-\dots -}_{\varepsilon_1({\bf c}^\circ_{1,4})}, \underbrace{+\dots +}_{\varphi_1({\bf c}^\circ_{1,4})}
		)
		=
		(\sigma_1, \sigma_2, \dots)
\end{equation*}
	be a finite sequence of $\pm$'s, where each $\varepsilon_1({\bf c}^\circ_{1,k})$ and $\varphi_1({\bf c}^\circ_{1,k})$ can be obtained by Remarks \ref{rem:crystal operator on i-comp-1} and \ref{rem:crystal operator on i-comp-2}.

	We replace a pair $(\sigma_i, \sigma_j) = (+,-)$ by $(\,\cdot\,,\,\cdot\,)$, where $i < j$ and $\sigma_k = \cdot$ for $i < k < j$, and then repeat this process until we get a sequence $\sigma^{\rm red}({\bf c}^\circ)$ with no $-$ placed to the right of $+$.
	We have
	\begin{equation*}
		\sigma({\bf c}^\circ) = (-,-,+,+,+,+,-,+), \quad
		\sigma^{\rm red}({\bf c}^\circ) = (-,-,+,+,+,\,\cdot\,,\,\cdot\,,+).
	\end{equation*}
	Then $\tf_1$ (resp.~$\te_1$) acts on ${\bf c}^\circ_{1,k}$ corresponding to the left-most $+$ (resp.~right-most $-$) in $\sigma^{\rm red}({\bf c}^\circ)$, so that $k = 1$ and 
	
	\begin{equation*}
		\tf_1 {\bf c} = 
		\begin{tikzpicture}[baseline=(current  bounding  box.center), every node/.style={scale=0.85}, scale=0.75]
		\node (glw_00f) at (0,0) {\red{\bf 1}};
		\node (glw_01m) at (0.7,0.7) {\red{\bf 2}};
		\node (glw_02m) at (1.4,1.4) {\bf 1};
		\node (glw_03m) at (2.1,2.1) {\bf 0};
		\node (glw_04m) at (2.8,2.8) {\bf 1};
		\node (glw_01f) at (1.4,0) {\red{\bf 2}};
		\node (glw_12m) at (2.1,0.7) {\bf 0};
		\node (glw_13m) at (2.8,1.4) {\bf 1};
		\node (glw_14m) at (3.5,2.1) {\bf 0};
		\node (glw_02f) at (2.8,0) {\color{gray}{$0$}};
		\node (glw_23m) at (3.5,0.7) {\color{gray}{$2$}};		
		\node (glw_24m) at (4.2,1.4) {\color{gray}{$1$}};
		\node (glw_03f) at (4.2,0) {\color{gray}{$1$}};
		\node (glw_34m) at (4.9,0.7) {\color{gray}{$0$}};
		\node (glw_04m) at (5.6,0) {\color{gray}{$1$}};
	\end{tikzpicture}
	\end{equation*}
where the entries in red denote the ones in ${\bf c}^\circ_{1,k}$ on which $\tf_1$ acts.
	Note that $\varphi_1({\bf c})$ (resp.~$\varepsilon_1({\bf c})$) is the number of $+$'s (resp.~$-$'s) in $\sigma^{\rm red}({\bf c}^\circ)$.
	In this case, $\varphi_1({\bf c}) = 4$ and $\varepsilon_1({\bf c}) = 2$.
	
	Next, let us compute $\tf_2 {\bf c}$. In this case, we consider the following bold-faced entries:
	\begin{equation*}
	\qquad
	{\bf c} = 
		\begin{tikzpicture}[baseline=(current  bounding  box.center), every node/.style={scale=0.85}, scale=0.95]
		\node (glw_00f) at (0,0) {\color{gray}{$2$}};
		\node (glw_01m) at (0.7,0.7) {\bf 2};
		\node (glw_02m) at (1.4,1.4) {\bf 1};
		\node (glw_03m) at (2.1,2.1) {\color{gray}{$0$}};
		\node (glw_04m) at (2.8,2.8) {\color{gray}{$1$}};
		\node (glw_01f) at (1.4,0) {\bf 1};
		\node (glw_12m) at (2.1,0.7) {\bf 0};
		\node (glw_13m) at (2.8,1.4) {\bf 1};
		\node (glw_14m) at (3.5,2.1) {\bf 0};
		\node (glw_02f) at (2.8,0) {\bf 0};
		\node (glw_23m) at (3.5,0.7) {\bf 2};		
		\node (glw_24m) at (4.2,1.4) {\bf 1};
		\node (glw_03f) at (4.2,0) {\color{gray}{$1$}};
		\node (glw_34m) at (4.9,0.7) {\color{gray}{$0$}};
		\node (glw_04m) at (5.6,0) {\color{gray}{$1$}};
		\end{tikzpicture}
		\quad
		{\bf c}^\circ =
		\begin{tikzpicture}[baseline=(current  bounding  box.center), every node/.style={scale=0.85}, scale=0.95]
		\node (glw_00f) at (0,0) {};
		\node (glw_01m) at (0.7,0.7) {\bf 2};
		\node (glw_02m) at (1.4,1.4) {\bf 1};
		\node (glw_03m) at (2.1,2.1) {};
		\node (glw_04m) at (2.8,2.8) {};
		\node (glw_01f) at (1.4,0) {\bf 1};
		\node (glw_12m) at (2.1,0.7) {\bf 0};
		\node (glw_13m) at (2.8,1.4) {\bf 1};
		\node (glw_14m) at (3.5,2.1) {\bf 0};
		\node (glw_02f) at (2.8,0) {\bf 0};
		\node (glw_23m) at (3.5,0.7) {\bf 2};		
		\node (glw_24m) at (4.2,1.4) {\bf 1};
		\node (glw_03f) at (4.2,0) {};
		\node (glw_34m) at (4.9,0.7) {};
		\node (glw_04m) at (5.6,0) {};
		\end{tikzpicture}
	\end{equation*}\smallskip
where we regard
\begin{equation*}
 {\bf c}^\circ 
		= \underbrace{\left( c_{(1,2)}, c_{(1,3)} \right)}_{{\bf c}^\circ_{2,1}} \ot 
		  \underbrace{\left( c_{(2,2)}, c_{(2,3)}, c_{(3,3)} \right)}_{{\bf c}^\circ_{2,2}} \ot
		  \underbrace{\left( c_{(2,4)}, c_{(3,4)} \right)}_{{\bf c}^\circ_{2,3}} \ot
		  \underbrace{\left( c_{(2,5)}, c_{(3,5)} \right)}_{{\bf c}^\circ_{2,4}}.
\end{equation*}
	with each ${\bf c}^\circ_{2,k}$ being an element of a crystal ($1\le k\le 4$).
 We have $\tf_2 {\bf c} = 0$, since $( {\rm wt}( {\bf c}^\circ_{2,1} ) \,|\, \alpha_2) \neq 0$.
	
	Finally, let us compute $\tf_3 {\bf c}$. In this case, we consider
	\begin{equation*}
	\qquad
	{\bf c} = 
	\begin{tikzpicture}[baseline=(current  bounding  box.center), every node/.style={scale=0.85}, scale=0.95]
		\node (glw_00f) at (0,0) {\color{gray}{$2$}};
		\node (glw_01m) at (0.7,0.7) {\color{gray}{$2$}};
		\node (glw_02m) at (1.4,1.4) {\bf 1};
		\node (glw_03m) at (2.1,2.1) {\bf 0};
		\node (glw_04m) at (2.8,2.8) {\color{gray}{$1$}};
		\node (glw_01f) at (1.4,0) {\color{gray}{$1$}};
		\node (glw_12m) at (2.1,0.7) {\bf 0};
		\node (glw_13m) at (2.8,1.4) {\bf 1};
		\node (glw_14m) at (3.5,2.1) {\color{gray}{$0$}};
		\node (glw_02f) at (2.8,0) {\bf 0};
		\node (glw_23m) at (3.5,0.7) {\bf 2};		
		\node (glw_24m) at (4.2,1.4) {\bf 1};
		\node (glw_03f) at (4.2,0) {\bf 1};
		\node (glw_34m) at (4.9,0.7) {\bf 0};
		\node (glw_04m) at (5.6,0) {\color{gray}{$1$}};
	\end{tikzpicture}
	\quad
	{\bf c}^\circ =
	\begin{tikzpicture}[baseline=(current  bounding  box.center), every node/.style={scale=0.85}, scale=0.95]
		\node (glw_00f) at (0,0) {};
		\node (glw_01m) at (0.7,0.7) {};
		\node (glw_02m) at (1.4,1.4) {\bf 1};
		\node (glw_03m) at (2.1,2.1) {\bf 0};
		\node (glw_04m) at (2.8,2.8) {};
		\node (glw_01f) at (1.4,0) {};
		\node (glw_12m) at (2.1,0.7) {\bf 0};
		\node (glw_13m) at (2.8,1.4) {\bf 1};
		\node (glw_14m) at (3.5,2.1) {};
		\node (glw_02f) at (2.8,0) {\bf 0};
		\node (glw_23m) at (3.5,0.7) {\bf 2};		
		\node (glw_24m) at (4.2,1.4) {\bf 1};
		\node (glw_03f) at (4.2,0) {\bf 1};
		\node (glw_34m) at (4.9,0.7) {\bf 0};
		\node (glw_04m) at (5.6,0) {};
	\end{tikzpicture}\ ,
	\end{equation*}
	where we regard
	\begin{equation*}
		{\bf c}^\circ = 
		\underbrace{\left( c_{(1,3)}, c_{(1,4)} \right)}_{{\bf c}^\circ_{3,1}} \ot
		\underbrace{\left( c_{(2,3)}, c_{(2,4)} \right)}_{{\bf c}^\circ_{3,2}} \ot
		\underbrace{\left( c_{(3,3)}, c_{(3,4)}, c_{(4,4)} \right)}_{{\bf c}^\circ_{3,3}} \ot
		\underbrace{\left( c_{(3,5)}, c_{(4,5)} \right)}_{{\bf c}^\circ_{3,4}}.
	\end{equation*}
	We apply the signature rule to ${\bf c}^{\circ}$ as in the case of $i=1$, but now in a reverse way since $(\alpha_3|\alpha_3)<0$. Let
\begin{equation*}
\zeta({\bf c}^\circ) = 
		( 
		\underbrace{+\dots +}_{\varphi_3({\bf c}^\circ_{3,1})}, 
		\underbrace{-\dots -}_{\varepsilon_3({\bf c}^\circ_{3,1})},
		\underbrace{+\dots +}_{\varphi_3({\bf c}^\circ_{3,2})}, 
		\underbrace{-\dots -}_{\varepsilon_3({\bf c}^\circ_{3,2})},
		\underbrace{+\dots +}_{\varphi_3({\bf c}^\circ_{3,3})}, 
		\underbrace{-\dots -}_{\varepsilon_3({\bf c}^\circ_{3,3})}, 
		\underbrace{+\dots +}_{\varphi_3({\bf c}^\circ_{3,4})},
		\underbrace{-\dots -}_{\varepsilon_3({\bf c}^\circ_{3,4})}
		)
		=
		(\zeta_1, \zeta_2, \dots).
\end{equation*}
	We replace a pair $(\zeta_i, \zeta_j) = (-,+)$ by $(\,\cdot\,,\,\cdot\,)$, where $i < j$ and $\zeta_k = \cdot$ for $i < k < j$, and then repeat this process until we get a sequence $\zeta^{\rm red}({\bf c}^\circ)$ with no $-$ placed to the left of $+$.
	We have
\begin{equation*}
 \zeta({\bf c}^\circ) = (+,-,-,+), \quad
		\zeta^{\rm red}({\bf c}^\circ) = (+,-,\,\cdot\,,\,\cdot\,).
\end{equation*}
	Then $\tf_3$ (resp.~$\te_3$) acts on ${\bf c}^\circ_{3,k}$ corresponding to the right-most $+$ (resp.~left-most $-$) in $\zeta^{\rm red}({\bf c}^\circ)$ so that $k = 1$ and 
	\begin{equation*}
	\tf_3 {\bf c} = 
	\begin{tikzpicture}[baseline=(current  bounding  box.center), every node/.style={scale=0.85}, scale=0.75]
		\node (glw_00f) at (0,0) {\color{gray}{$2$}};
		\node (glw_01m) at (0.7,0.7) {\color{gray}{$2$}};
		\node (glw_02m) at (1.4,1.4) {\red{\bf 0}};
		\node (glw_03m) at (2.1,2.1) {\red{\bf 1}};
		\node (glw_04m) at (2.8,2.8) {\color{gray}{$1$}};
		\node (glw_01f) at (1.4,0) {\color{gray}{$1$}};
		\node (glw_12m) at (2.1,0.7) {\bf 0};
		\node (glw_13m) at (2.8,1.4) {\bf 1};
		\node (glw_14m) at (3.5,2.1) {\color{gray}{$0$}};
		\node (glw_02f) at (2.8,0) {\bf 0};
		\node (glw_23m) at (3.5,0.7) {\bf 2};		
		\node (glw_24m) at (4.2,1.4) {\bf 1};
		\node (glw_03f) at (4.2,0) {\bf 1};
		\node (glw_34m) at (4.9,0.7) {\bf 0};
		\node (glw_04m) at (5.6,0) {\color{gray}{$1$}};
	\end{tikzpicture}\quad .
	\end{equation*}
	Note that $\varphi_3({\bf c})$ (resp.~$\varepsilon_3({\bf c})$) is the number of $+$'s (resp.~$-$'s) in $\zeta^{\rm red}({\bf c}^\circ)$.
	In this case, $\varphi_3({\bf c}) = 1$ and $\varepsilon_3({\bf c}) = 1$.
	\smallskip

\noindent
(2) Suppose that $\mf{g}=\mf{c}_{2|3}$. Let us compute $\tf_i^k {\bf c}$  ($1\le k\le \varphi_i({\bf c})$). We leave the details to the reader since it can be done in the same way as in (1).

Let ${\bf c}=\left( c_{(i,j)} \right)_{(i,j)\in \Phi^+(\mf{u})}\in M(\mf{u})$ be given by
\begin{equation*}
	\qquad
	{\bf c} = 
	\begin{tikzpicture}[baseline=(current  bounding  box.center), every node/.style={scale=0.9}]
		\node (glw_00f) at (0,0) {$c_{(1,1)}$};
		\node (glw_01m) at (0.7,0.7) {$c_{(1,2)}$};
		\node (glw_02m) at (1.4,1.4) {$c_{(1,3)}$};
		\node (glw_03m) at (2.1,2.1) {$c_{(1,4)}$};
		\node (glw_04m) at (2.8,2.8) {$c_{(1,5)}$};
		\node (glw_01m) at (1.4,0) {$c_{(2,2)}$};
		\node (glw_12m) at (2.1,0.7) {$c_{(2,3)}$};
		\node (glw_13m) at (2.8,1.4) {$c_{(2,4)}$};
		\node (glw_14m) at (3.5,2.1) {$c_{(2,5)}$};
		\node (glw_23m) at (3.5,0.7) {$c_{(3,4)}$};		
		\node (glw_24m) at (4.2,1.4) {$c_{(3,5)}$};
		\node (glw_34m) at (4.9,0.7) {$c_{(4,5)}$};
	\end{tikzpicture}
	=
	\begin{tikzpicture}[baseline=(current  bounding  box.center), every node/.style={scale=0.9}]
		\node (glw_00f) at (0,0) {$1$};
		\node (glw_01m) at (0.7,0.7) {$0$};
		\node (glw_02m) at (1.4,1.4) {$1$};
		\node (glw_03m) at (2.1,2.1) {$0$};
		\node (glw_04m) at (2.8,2.8) {$1$};
		\node (glw_01m) at (1.4,0) {$0$};
		\node (glw_12m) at (2.1,0.7) {$0$};
		\node (glw_13m) at (2.8,1.4) {$1$};
		\node (glw_14m) at (3.5,2.1) {$0$};
		\node (glw_23m) at (3.5,0.7) {$0$};		
		\node (glw_24m) at (4.2,1.4) {$2$};
		\node (glw_34m) at (4.9,0.7) {$1$};
	\end{tikzpicture}
\end{equation*}

For $i=1$, we have $\varphi_1({\bf c}) = 3$, $\varepsilon_1({\bf c}) = 0$, and
\begin{equation*}
\begin{split}
	\tf_1^{\, 3} {\bf c}&=\, \tf_1^2
	\begin{tikzpicture}[baseline=(current  bounding  box.center), every node/.style={scale=0.9}, scale=0.75]
		\node (glw_00f) at (0,0) {\red{\bf 0}};
		\node (glw_01m) at (0.7,0.7) {\red{\bf 1}};
		\node (glw_02m) at (1.4,1.4) {\bf 1};
		\node (glw_03m) at (2.1,2.1) {\bf 0};
		\node (glw_04m) at (2.8,2.8) {\bf 1};
		\node (glw_01m) at (1.4,0) {\bf 0};
		\node (glw_12m) at (2.1,0.7) {\bf 0};
		\node (glw_13m) at (2.8,1.4) {\bf 1};
		\node (glw_14m) at (3.5,2.1) {\bf 0};
		\node (glw_23m) at (3.5,0.7) {\color{gray}{$0$}};		
		\node (glw_24m) at (4.2,1.4) {\color{gray}{$2$}};
		\node (glw_34m) at (4.9,0.7) {\color{gray}{$1$}};
	\end{tikzpicture}
	\!
	=
	\,\tf_1
	\begin{tikzpicture}[baseline=(current  bounding  box.center), every node/.style={scale=0.9}, scale=0.75]
		\node (glw_00f) at (0,0) {\bf 0};
		\node (glw_01m) at (0.7,0.7) {\red{\bf 0}};
		\node (glw_02m) at (1.4,1.4) {\bf 1};
		\node (glw_03m) at (2.1,2.1) {\bf 0};
		\node (glw_04m) at (2.8,2.8) {\bf 1};
		\node (glw_01m) at (1.4,0) {\red{\bf 1}};
		\node (glw_12m) at (2.1,0.7) {\bf 0};
		\node (glw_13m) at (2.8,1.4) {\bf 1};
		\node (glw_14m) at (3.5,2.1) {\bf 0};
		\node (glw_23m) at (3.5,0.7) {\color{gray}{$0$}};		
		\node (glw_24m) at (4.2,1.4) {\color{gray}{$2$}};
		\node (glw_34m) at (4.9,0.7) {\color{gray}{$1$}};
	\end{tikzpicture}
	\!\\
	& =
	\begin{tikzpicture}[baseline=(current  bounding  box.center), every node/.style={scale=0.9}, scale=0.75]
		\node (glw_00f) at (0,0) {\bf 0};
		\node (glw_01m) at (0.7,0.7) {\bf 0};
		\node (glw_02m) at (1.4,1.4) {\bf 1};
		\node (glw_03m) at (2.1,2.1) {\bf 0};
		\node (glw_04m) at (2.8,2.8) {\red{\bf 0}};
		\node (glw_01m) at (1.4,0) {\bf 1};
		\node (glw_12m) at (2.1,0.7) {\bf 0};
		\node (glw_13m) at (2.8,1.4) {\bf 1};
		\node (glw_14m) at (3.5,2.1) {\red{\bf 1}};
		\node (glw_23m) at (3.5,0.7) {\color{gray}{$0$}};		
		\node (glw_24m) at (4.2,1.4) {\color{gray}{$2$}};
		\node (glw_34m) at (4.9,0.7) {\color{gray}{$1$}};
	\end{tikzpicture}\quad .
	\end{split}
\end{equation*}
For $i=3$, we have $\varphi_3 ( {\bf c} ) = 2$, $\varepsilon_3 ( {\bf c} ) = 1$, and 
\begin{equation*}
	\tf_3^{\,2} {\bf c}=\tf_3 
	\begin{tikzpicture}[baseline=(current  bounding  box.center), every node/.style={scale=0.9}, scale=0.75]
		\node (glw_00f) at (0,0) {\color{gray}{$1$}};
		\node (glw_01m) at (0.7,0.7) {\color{gray}{$0$}};
		\node (glw_02m) at (1.4,1.4) {\bf 1};
		\node (glw_03m) at (2.1,2.1) {\bf 0};
		\node (glw_04m) at (2.8,2.8) {\color{gray}{$1$}};
		\node (glw_01m) at (1.4,0) {\color{gray}{$0$}};
		\node (glw_12m) at (2.1,0.7) {\bf 0};
		\node (glw_13m) at (2.8,1.4) {\bf 1};
		\node (glw_14m) at (3.5,2.1) {\color{gray}{$0$}};
		\node (glw_23m) at (3.5,0.7) {{\bf 0}};		
		\node (glw_24m) at (4.2,1.4) {\red{\bf 1}};
		\node (glw_34m) at (4.9,0.7) {\red{\bf 2}};
	\end{tikzpicture}
	\,=
	\begin{tikzpicture}[baseline=(current  bounding  box.center), every node/.style={scale=0.9}, scale=0.75]
		\node (glw_00f) at (0,0) {\color{gray}{$1$}};
		\node (glw_01m) at (0.7,0.7) {\color{gray}{$0$}};
		\node (glw_02m) at (1.4,1.4) {\red{\bf 0}};
		\node (glw_03m) at (2.1,2.1) {\red{\bf 1}};
		\node (glw_04m) at (2.8,2.8) {\color{gray}{$1$}};
		\node (glw_01m) at (1.4,0) {\color{gray}{$0$}};
		\node (glw_12m) at (2.1,0.7) {\bf 0};
		\node (glw_13m) at (2.8,1.4) {\bf 1};
		\node (glw_14m) at (3.5,2.1) {\color{gray}{$0$}};
		\node (glw_23m) at (3.5,0.7) {{\bf 0}};		
		\node (glw_24m) at (4.2,1.4) {\bf 1};
		\node (glw_34m) at (4.9,0.7) {\bf 2};
	\end{tikzpicture}\quad .
\end{equation*}
\smallskip

\noindent
(3) Suppose that $\mf{g}=\mf{d}_{3|3}$.
Let ${\bf c}=\left( c_{(i,j)} \right)_{(i,j)\in \Phi^+(\mf{u})}\in M(\mf{u})$ be given by
\begin{equation*}
	\qquad
	{\bf c} = 
	\begin{tikzpicture}[baseline=(current  bounding  box.center), every node/.style={scale=0.9}, scale=0.85]
		\node (glw_00f) at (0,0) {$c_{(1,2)}$};
		\node (glw_01m) at (0.7,0.7) {$c_{(1,3)}$};
		\node (glw_02m) at (1.4,1.4) {$c_{(1,4)}$};
		\node (glw_03m) at (2.1,2.1) {$c_{(1,5)}$};
		\node (glw_04m) at (2.8,2.8) {$c_{(1,6)}$};
		\node (glw_01f) at (1.4,0) {$c_{(2,3)}$};
		\node (glw_12m) at (2.1,0.7) {$c_{(2,4)}$};
		\node (glw_13m) at (2.8,1.4) {$c_{(2,5)}$};
		\node (glw_14m) at (3.5,2.1) {$c_{(2,6)}$};
		\node (glw_02f) at (2.8,0) {$c_{(3,4)}$};
		\node (glw_23m) at (3.5,0.7) {$c_{(3,5)}$};		
		\node (glw_24m) at (4.2,1.4) {$c_{(3,6)}$};
		\node (glw_03f) at (3.5,-0.7) {$c_{(4,4)}$};
		\node (glw_03f) at (4.2,0) {$c_{(4,5)}$};
		\node (glw_34m) at (4.9,0.7) {$c_{(4,6)}$};
		\node (glw_04m) at (4.9,-0.7) {$c_{(5,5)}$};
		\node (glw_04m) at (5.6,0) {$c_{(5,6)}$};
		\node (glw_04m) at (6.3,-0.7) {$c_{(6,6)}$};
	\end{tikzpicture}
	=
	\begin{tikzpicture}[baseline=(current  bounding  box.center), every node/.style={scale=0.9}, scale=0.85]
		\node (glw_00f) at (0,0) {$0$};
		\node (glw_01m) at (0.7,0.7) {$2$};
		\node (glw_02m) at (1.4,1.4) {$0$};
		\node (glw_03m) at (2.1,2.1) {$1$};
		\node (glw_04m) at (2.8,2.8) {$0$};
		\node (glw_01f) at (1.4,0) {$1$};
		\node (glw_12m) at (2.1,0.7) {$1$};
		\node (glw_13m) at (2.8,1.4) {$0$};
		\node (glw_14m) at (3.5,2.1) {$1$};
		\node (glw_02f) at (2.8,0) {$1$};
		\node (glw_23m) at (3.5,0.7) {$0$};		
		\node (glw_24m) at (4.2,1.4) {$1$};
		\node (glw_03f) at (3.5,-0.7) {$2$};
		\node (glw_03f) at (4.2,0) {$2$};
		\node (glw_34m) at (4.9,0.7) {$3$};
		\node (glw_04m) at (4.9,-0.7) {$1$};
		\node (glw_04m) at (5.6,0) {$2$};
		\node (glw_04m) at (6.3,-0.7) {$0$};
	\end{tikzpicture}.
\end{equation*}
Let us compute $\tf_4^k {\bf c}$ with $k=\varphi_4({\bf c})$. 
We consider the bold-faced entries
\begin{equation*}
	\qquad
	{\bf c} = 
	\begin{tikzpicture}[baseline=(current  bounding  box.center), every node/.style={scale=0.9}, scale=0.85]
		\node (glw_00f) at (0,0) {\color{gray}{$0$}};
		\node (glw_01m) at (0.7,0.7) {\color{gray}{$2$}};
		\node (glw_02m) at (1.4,1.4) {\bf 0};
		\node (glw_03m) at (2.1,2.1) {\bf 1};
		\node (glw_04m) at (2.8,2.8) {\color{gray}{$0$}};
		\node (glw_01f) at (1.4,0) {\color{gray}{$1$}};
		\node (glw_12m) at (2.1,0.7) {\bf 1};
		\node (glw_13m) at (2.8,1.4) {\bf 0};
		\node (glw_14m) at (3.5,2.1) {\color{gray}{$1$}};
		\node (glw_02f) at (2.8,0) {\bf 1};
		\node (glw_23m) at (3.5,0.7) {\bf 0};		
		\node (glw_24m) at (4.2,1.4) {\color{gray}{$1$}};
		\node (glw_03f) at (3.5,-0.7) {\bf 2};
		\node (glw_03f) at (4.2,0) {\bf 2};
		\node (glw_34m) at (4.9,0.7) {\bf 3};
		\node (glw_04m) at (4.9,-0.7) {\bf 1};
		\node (glw_04m) at (5.6,0) {\bf 2};
		\node (glw_04m) at (6.3,-0.7) {\color{gray}{$0$}};
	\end{tikzpicture}
	\quad
	{\bf c}^\circ = 
		\begin{tikzpicture}[baseline=(current  bounding  box.center), every node/.style={scale=0.9}, scale=0.85]
		\node (glw_00f) at (0,0) {};
		\node (glw_01m) at (0.7,0.7) {};
		\node (glw_02m) at (1.4,1.4) {\bf 0};
		\node (glw_03m) at (2.1,2.1) {\bf 1};
		\node (glw_04m) at (2.8,2.8) {};
		\node (glw_01f) at (1.4,0) {};
		\node (glw_12m) at (2.1,0.7) {\bf 1};
		\node (glw_13m) at (2.8,1.4) {\bf 0};
		\node (glw_14m) at (3.5,2.1) {};
		\node (glw_02f) at (2.8,0) {\bf 1};
		\node (glw_23m) at (3.5,0.7) {\bf 0};		
		\node (glw_24m) at (4.2,1.4) {};
		\node (glw_03f) at (3.5,-0.7) {\bf 2};
		\node (glw_03f) at (4.2,0) {\bf 2};
		\node (glw_34m) at (4.9,0.7) {\bf 3};
		\node (glw_04m) at (4.9,-0.7) {\bf 1};
		\node (glw_04m) at (5.6,0) {\bf 2};
		\node (glw_04m) at (6.3,-0.7) {};
	\end{tikzpicture}
\end{equation*}
where
\begin{equation*}
	{\bf c}^\circ = 
	\underbrace{\left( c_{(1,4)}, c_{(1,5)} \right)}_{{\bf c}^\circ_{4,1}} \ot
	\underbrace{\left( c_{(2,4)}, c_{(2,5)} \right)}_{{\bf c}^\circ_{4,2}} \ot
	\underbrace{\left( c_{(3,4)}, c_{(3,5)} \right)}_{{\bf c}^\circ_{4,3}} \ot
	\underbrace{\left( c_{(4,4)}, c_{(4,5)}, c_{(5,5)} \right)}_{{\bf c}^\circ_{4,4}} \ot
	\underbrace{\left( c_{(5,6)}, c_{(6,6)} \right)}_{{\bf c}^\circ_{4,5}}.
\end{equation*}
Then
\begin{equation*}
\begin{split}
	\zeta_4({\bf c}^\circ) &= (-,+,+,+,+,+,+,-,-,+,+,+,-,-),\\
	\zeta_4^{\rm red}({\bf c}^\circ) &= (\,\,\cdot\,,\,\,\cdot\,,+,+,+,+,+,\,\,\cdot\,,\,\,\cdot\,,\,\,\cdot\,,\,\,\cdot\,,+,-,-),
\end{split}
\end{equation*}
so that $\varphi_4({\bf c}) = 6$ and $\varepsilon_4({\bf c}) = 2$.
We have
\begin{equation*}
	\tf_4^{\,6} {\bf c} = 
	\begin{tikzpicture}[baseline=(current  bounding  box.center), every node/.style={scale=0.9}, scale=0.75]
		\node (glw_00f) at (0,0) {\color{gray}{$0$}};
		\node (glw_01m) at (0.7,0.7) {\color{gray}{$2$}};
		\node (glw_02m) at (1.4,1.4) {\bf 0};
		\node (glw_03m) at (2.1,2.1) {\bf 1};
		\node (glw_04m) at (2.8,2.8) {\color{gray}{$0$}};
		\node (glw_01f) at (1.4,0) {\color{gray}{$1$}};
		\node (glw_12m) at (2.1,0.7) {\bf 1};
		\node (glw_13m) at (2.8,1.4) {\bf 0};
		\node (glw_14m) at (3.5,2.1) {\color{gray}{$1$}};
		\node (glw_02f) at (2.8,0) {\red{\bf 0}};
		\node (glw_23m) at (3.5,0.7) {\red{\bf 1}};		
		\node (glw_24m) at (4.2,1.4) {\color{gray}{$1$}};
		\node (glw_03f) at (3.5,-0.7) {\red{\bf 0}};
		\node (glw_03f) at (4.2,0) {\red{\bf 2}};
		\node (glw_34m) at (4.9,0.7) {\red{\bf 2}};
		\node (glw_04m) at (4.9,-0.7) {\red{\bf 3}};
		\node (glw_04m) at (5.6,0) {\red{\bf 3}};
		\node (glw_04m) at (6.3,-0.7) {\color{gray}{$0$}};
	\end{tikzpicture}\quad .
\end{equation*}
} 
\end{ex}\smallskip

\begin{rem}{\rm
We have the following decomposition of $\mc{N}$ into irreducible polynomial $\U(\mf{l})$-modules:
\begin{equation}\label{eq:decomp of crystal of N}
 \mc{N} \cong \bigoplus_{\mu\in \cP(\g)}V_{\mf l}(\mu),
\end{equation}
where $\cP(\mf{b}_{m|n})$ is the set of $(m|n)$-hook partitions, $\cP(\mf{c}_{m|n})$ is the set of $(m|n)$-hook even partitions, and $\cP(\mf{d}_{m|n})$ is the set of $(m|n)$-hook partitions whose conjugates are even. It follows from the character of $\mc{N}$ given as a sum of characters of $V_{\mf l}(\mu)$, which can be obtained by applying standard argument of superization  to classical Littlewood identities (cf.~\cite[I.\,3.\,Example 23, I.\,5.\,Examples 4,5]{Mac}).

We may obtain \eqref{eq:decomp of crystal of N} more directly using crystals when $\g=\mf{b}_{m|n}, \mf{c}_{m|n}$. In this case, the $\mf l$-crystal
$\ms{B}(\mc{N})$ is equal to the one given in \cite[Section 3.4]{K07}, where the elements equivalent to genuine highest weight vectors are characterized under an analogue of RSK correspondence. This implies the decomposition \eqref{eq:decomp of crystal of N} as an $\mf l$-crystal. 
We also expect a similar result in case of $\g=\mf{d}_{m|n}$ by using a super analogue of crystal theoretic interpretation of Burge correspondence of type $D$ (cf.~\cite{JK19}).
 }
\end{rem}

\subsection{Parabolic Verma module $P(\la)$}\label{subsec:parabolic Ver}
Let $\U(\mf{p})$ be the subalgebra of $\U(\g)$ generated by $\U(\mf{l})$ and $\U(\g)^{+}$.
Let $\la\in \cP_{m|n}$ be given. Let
\begin{equation*}\label{eq:P(la)}
 P(\la)=\U(\g)\ot_{\U(\mf{p})}V_{\mf{l}}(\la),
\end{equation*}
where we regard $V_{\mf{l}}(\la)$ as a $\U(\mf{p})$-module with $e_0 v_\la=0$.
We have an isomorphism of $\Bbbk$-spaces
\begin{equation}\label{eq:P(la)-2}
\mc{N}\ot_{\Bbbk} V_{\mf l}(\la)\cong  P(\la),
\end{equation}
sending $u\ot v$ to itself since 
$\U(\g)\cong \U(\mf{u}^-)\ot \U(\mf{l})^-\ot \U(\g)^{\ge 0}\cong \mc{N} \ot \U(\mf{p})$. It is a highest weight $\U(\mf{g})$-module with highest weight $\Lambda_\la$ \eqref{eq:highest weight correspondence}. We denote by $w_\la=1\ot v_\la$ a highest weight vector of $P(\la)$.

\begin{lem}\label{lem:Plambda poly}
We have $P(\la)\in \mc{O}_{\ge 0}$ as a $\U(\mf{l})$-module.
\end{lem}
\pf Let $w\in P(\la)$ be a weight vector such that $w=u\ot v$ for some $u\in \U(\mf{u}^-)_{\beta}$ and $v\in V_{\mf l}(\la)_\mu$ (see \eqref{eq:P(la)-2}). 
For $\nu\in P$, we have $k_\nu w=k_\nu u k_{\nu}^{-1}\ot k_\nu v = \bq(\beta,\nu)\bq(\mu,\nu)u\ot v=\bq(\beta+\mu,\nu)w$. Since $-\Phi^+(\mf{u}) \subset P_{\ge 0}$, we have $\beta\in P_{\ge 0}$ and hence ${\rm wt}(w)=\beta+\mu\in P_{\ge 0}$. Also we have $\dim P(\la)_\mu<\infty$ for $\mu\in P_{\ge 0}$ (cf.~Lemma \ref{lem:N poly}).
\qed

\begin{prop}\label{prop:par Verma and N are iso}
The $\Bbbk$-linear map in \eqref{eq:P(la)-2}
is an isomorphism of $\U(\mf{l})$-modules.
\end{prop}
\pf It is easy to see from \eqref{eq: quantum adjoint}  that 
\begin{equation*}
\begin{split}
 k_\mu u& = (k_\mu\cdot u_1)\ot k_\mu u_2,\\
 e_i u &= (e_i\cdot u_1)\ot k_i^{-1}u_2 + u_1\ot e_iu_2 ,\\
 f_i u &= (f_i\cdot u_1)\ot u_2 + (k_i\cdot u_1)\ot f_iu_2, 
\end{split}
\end{equation*}
for $\mu\in P$, $i\in I\setminus\{0\}$, and $u = u_1\ot u_2\in P(\la)$ with $u_1 \in \mc{N}$, which coincide with the action of $\U(\mf{l})$ on $\mc{N}\ot V_{\mf{l}}(\la)$. This implies that the map in \eqref{eq:P(la)-2} is $\U(\mf{l})$-linear.
\qed

\subsection{Crystal base of $P(\la)$}


Let $\la\in \cP_{m|n}$ be given. Note that $P(\la)$ is a quotient of $\U(\mf{g})^- \cong \mc{N} \ot \U(\mf{l})^-$ by $\mc{N} \ot I_\la$ as a $\Bbbk$-space, where $I_\la$ is the kernel of the natural projection from $\U(\mf{l})^-$ to $V_{\mf{l}}(\la)$. Since $e'_0$ on $\U(\mf{g})^-$ acts on $\U(\mf{l})^-$ and hence on $I_\la$ trivially, it induces a well-defined map on $P(\la)$.

We define the crystal operators $\te_0$ and $\tf_0$ on $P(\la)$ as in the case of symmetrizable Kac-Moody algebras \cite{Kas91}, that is,  for a weight vector $u\in P(\la)$ with $u=\sum_{k\ge 0}f_0^{(k)}u_k$ for $e'_0 ( u_k ) = 0$, define 
\begin{equation*}
 \te_0 u = \sum_{k\ge 1}f_0^{(k-1)}u_k,\quad \tf_0 u = \sum_{k\ge 0}f_0^{(k+1)}u_k.
\end{equation*}
Then a crystal base of $P(\la)$ is defined in the same way as in Section \ref{subsec:poly repn} with respect to $\te_i$ and $\tf_i$ for $i\in I$.

\begin{lem}\label{lem:e' trivially except f0}
For $\beta\in \Phi^+(\mf{u})\setminus\{\alpha_0\}$, we have $e'_0(\bff_\beta)=0$.
\end{lem}
\pf We can check directly by using \eqref{eq:e'-formula} and induction on ${\rm ht}(\beta)$.
\qed
\smallskip

  
Let
\begin{equation*}
\ms{L}(P(\la))=\ms{L}(\mc{N})\ot \ms{L}(V_{\mf{l}}(\la)),\quad
\ms{B}(P(\la))=\ms{B}(\mc{N})\ot \ms{B}(V_{\mf{l}}(\la)),
\end{equation*}
(see \eqref{eq:crystal base of poly} and \eqref{eq:crystal base of N}). Then we have the following.

\begin{thm}
The pair $(\ms{L}(P(\la)),\ms{B}(P(\la)))$ is a crystal base of $P(\la)$ as a $\U(\mf{g})$-module.
\end{thm}
\pf By Proposition \ref{prop:tensor product thm} and Theorem \ref{thm:crystal base of N}, $(\ms{L}(P(\la)),\ms{B}(P(\la)))$ is a crystal base of $P(\la)$ as a $\U(\mf{l})$-module. 
On the other hand, since $\alpha_0$ is minimal with respect to $\prec$ on $\Phi^+(\mf{u})$, it follows from Lemma \ref{lem:e' trivially except f0} that
\begin{equation}\label{eq:action of crystal f0}
 \tf_0 \bfF^{(\bf c)} 
 = \tf_0 \prod^{\rightarrow}_{\beta\in\Phi^+(\mf{u})}\bfF_\beta^{(c_\beta)}
 = f_0^{(c_{\alpha_0}+1)} \!\!\!\!\!\!\!\!\!\!\! \prod^{\rightarrow}_{\beta\in\Phi^+(\mf{u})\setminus\{\alpha_0\}} \!\!\!\!\! \bfF_\beta^{(c_\beta)},
\end{equation}
for ${\bf c}=(c_\beta)\in M(\mf{u})$. This implies that $(\ms{L}(P(\la)),\ms{B}(P(\la)))$ satisfies the conditions for a crystal base with respect to $\te_0$ and $\tf_0$, and hence it is a crystal base of $P(\la)$ as a $\U(\mf{g})$-module.
\qed\smallskip

\begin{ex}
{\em 
Let us give some examples of $\tf_i b$ for $b\in \ms{B}(P(\la))$.
We continue Example \ref{ex:crystal structure on N}~(1), where $\mf{g}=\mf{b}_{2|3}$.
In this case, $\I_{\ov 0} = \{ 1, 2 \}$ and $\I_{\ov 1} = \{ 3, 4, 5 \}$.
Let $\la = \left( 5, 3, 3, 2, 2 \right) \in \ms{P}_{2|3}$, and $\ms{B}(V_{\mf{l}}(\la)) = SST_{2|3}(\la)$.
%
Let $b = {\bf c} \ot T \in \ms{B}(P(\la))$ be given by
\begin{equation*}
	{\bf c} \ot T = 
	\left(\,\,
	\begin{tikzpicture}[baseline=(current  bounding  box.center), every node/.style={scale=0.85}, scale=0.75]
		\node (glw_00f) at (0,0) {$2$};
		\node (glw_01m) at (0.7,0.7) {$2$};
		\node (glw_02m) at (1.4,1.4) {$1$};
		\node (glw_03m) at (2.1,2.1) {$0$};
		\node (glw_04m) at (2.8,2.8) {$1$};
		\node (glw_01f) at (1.4,0) {$1$};
		\node (glw_12m) at (2.1,0.7) {$0$};
		\node (glw_13m) at (2.8,1.4) {$1$};
		\node (glw_14m) at (3.5,2.1) {$0$};
		\node (glw_02f) at (2.8,0) {$0$};
		\node (glw_23m) at (3.5,0.7) {$2$};		
		\node (glw_24m) at (4.2,1.4) {$1$};
		\node (glw_03f) at (4.2,0) {$1$};
		\node (glw_34m) at (4.9,0.7) {$0$};
		\node (glw_04m) at (5.6,0) {$1$};
	\end{tikzpicture}
	\, \ot \,\,\,\,\,\,
	\raisebox{0.7cm}{\scalebox{0.85}{
	\begin{ytableau}
		1 & 1 & 1 & 2 & 2 \\
		2 & 2 & 3 & \none & \none \\
		3 & 4 & 5 & \none & \none \\
		4 & 5 & \none & \none & \none \\
		4 & 5 & \none & \none & \none \\
	\end{ytableau}}}
	\,\,\right).
\end{equation*}
Let us compute $\tf_i^k ({\bf c} \ot T)$ $(k \ge 1)$ by using Proposition \ref{prop:tensor product thm} (cf.~Example \ref{ex:crystal structure on N}~(1)).\smallskip

\noindent (1) For $i = 0$, we have by definition of $\tf_0$
\begin{equation*}
	\tf_0 ({\bf c} \ot T) = (\tf_0 {\bf c}) \ot T =  
	\left(\,\,
	\begin{tikzpicture}[baseline=(current  bounding  box.center), every node/.style={scale=0.85}, scale=0.75]
		\node (glw_00f) at (0,0) {\red{\bf 3}};
		\node (glw_01m) at (0.7,0.7) {\color{gray}{$2$}};
		\node (glw_02m) at (1.4,1.4) {\color{gray}{$1$}};
		\node (glw_03m) at (2.1,2.1) {\color{gray}{$0$}};
		\node (glw_04m) at (2.8,2.8) {\color{gray}{$1$}};
		\node (glw_01f) at (1.4,0) {\color{gray}{$1$}};
		\node (glw_12m) at (2.1,0.7) {\color{gray}{$0$}};
		\node (glw_13m) at (2.8,1.4) {\color{gray}{$1$}};
		\node (glw_14m) at (3.5,2.1) {\color{gray}{$0$}};
		\node (glw_02f) at (2.8,0) {\color{gray}{$0$}};
		\node (glw_23m) at (3.5,0.7) {\color{gray}{$2$}};		
		\node (glw_24m) at (4.2,1.4) {\color{gray}{$1$}};
		\node (glw_03f) at (4.2,0) {\color{gray}{$1$}};
		\node (glw_34m) at (4.9,0.7) {\color{gray}{$0$}};
		\node (glw_04m) at (5.6,0) {\color{gray}{$1$}};
	\end{tikzpicture}
	\,\ot\,\,\,\,\,\,
	\raisebox{0.7cm}{\scalebox{0.85}{
	\begin{ytableau}
		1 & 1 & 1 & 2 & 2 \\
		2 & 2 & 3 & \none & \none \\
		3 & 4 & 5 & \none & \none \\
		4 & 5 & \none & \none & \none \\
		4 & 5 & \none & \none & \none \\
	\end{ytableau}}}
	\,\,\right).
\end{equation*}\smallskip

\noindent (2) For $i=1$, we have
\begin{equation*}
\begin{split}
	\tf_1^{\,3} \left( {\bf c} \ot T \right) &= \tf_1 \left( (\tf_1^{\,2} {\bf c}) \ot T \right) = \left( \tf_1^{\,2} {\bf c} \right) \ot  \left( \tf_1 T \right) \\
	&= 
	\left(\,\,
	\begin{tikzpicture}[baseline=(current  bounding  box.center), every node/.style={scale=0.85}, scale=0.75]
		\node (glw_00f) at (0,0) {\red{\bf 1}};
		\node (glw_01m) at (0.7,0.7) {\red{\bf 1}};
		\node (glw_02m) at (1.4,1.4) {\bf 1};
		\node (glw_03m) at (2.1,2.1) {\bf 0};
		\node (glw_04m) at (2.8,2.8) {\bf 1};
		\node (glw_01f) at (1.4,0) {\red{\bf 4}};
		\node (glw_12m) at (2.1,0.7) {\bf 0};
		\node (glw_13m) at (2.8,1.4) {\bf 1};
		\node (glw_14m) at (3.5,2.1) {\bf 0};
		\node (glw_02f) at (2.8,0) {\color{gray}{$0$}};
		\node (glw_23m) at (3.5,0.7) {\color{gray}{$2$}};		
		\node (glw_24m) at (4.2,1.4) {\color{gray}{$1$}};
		\node (glw_03f) at (4.2,0) {\color{gray}{$1$}};
		\node (glw_34m) at (4.9,0.7) {\color{gray}{$0$}};
		\node (glw_04m) at (5.6,0) {\color{gray}{$1$}};
	\end{tikzpicture}
	\, \ot \,\,\,\,\,\,
	\raisebox{0.7cm}{\scalebox{0.85}{
	\begin{ytableau}
		1 & 1 & \red{\bf 2} & 2 & 2 \\
		2 & 2 & 3 & \none & \none \\
		3 & 4 & 5 & \none & \none \\
		4 & 5 & \none & \none & \none \\
		4 & 5 & \none & \none & \none \\
	\end{ytableau}}}
	\,\,\right),
\end{split}
\end{equation*}
since $(\varepsilon_1 ({\bf c}), \varphi_1({\bf c})) = (2, 4)$ and $(\varepsilon_1 ( T ), \varphi_1 ( T )) = (2, 1)$.\smallskip

\noindent (3) For $i=2$, we have $\tf_2 ({\bf c} \ot T) = (\tf_2 {\bf c}) \ot T = 0$, since $\left( {\rm wt}({\bf c}) | \alpha_2 \right) \neq 0$.\smallskip

\noindent (4) For $i = 3$,
we have
\begin{equation*}
\begin{split}
	\tf_3 \left( {\bf c} \ot T \right) &= \left( \tf_3 {\bf c} \right) \ot T =
	\left(\,\,
	\begin{tikzpicture}[baseline=(current  bounding  box.center), every node/.style={scale=0.85}, scale=0.75]
		\node (glw_00f) at (0,0) {\color{gray}{$2$}};
		\node (glw_01m) at (0.7,0.7) {\color{gray}{$2$}};
		\node (glw_02m) at (1.4,1.4) {\red{\bf 0}};
		\node (glw_03m) at (2.1,2.1) {\red{\bf 1}};
		\node (glw_04m) at (2.8,2.8) {\color{gray}{$1$}};
		\node (glw_01f) at (1.4,0) {\color{gray}{$1$}};
		\node (glw_12m) at (2.1,0.7) {\bf 0};
		\node (glw_13m) at (2.8,1.4) {\bf 1};
		\node (glw_14m) at (3.5,2.1) {\color{gray}{$0$}};
		\node (glw_02f) at (2.8,0) {\bf 0};
		\node (glw_23m) at (3.5,0.7) {\bf 2};		
		\node (glw_24m) at (4.2,1.4) {\bf 1};
		\node (glw_03f) at (4.2,0) {\bf 1};
		\node (glw_34m) at (4.9,0.7) {\bf 0};
		\node (glw_04m) at (5.6,0) {\color{gray}{$1$}};
	\end{tikzpicture}
	\,\ot\,\,\,\,\,\,
	\raisebox{0.7cm}{\scalebox{0.85}{
	\begin{ytableau}
		1 & 1 & 1 & 2 & 2 \\
		2 & 2 & 3 & \none & \none \\
		3 & 4 & 5 & \none & \none \\
		4 & 5 & \none & \none & \none \\
		4 & 5 & \none & \none & \none \\
	\end{ytableau}}}
	\,\,\right),
\end{split}
\end{equation*}
since $(\varepsilon_3 \left( {\bf c} \right), \varphi_3 \left( {\bf c} \right)) = (1, 1)$ and $(\varepsilon_3 \left( T \right), \varphi_3 \left( T \right)) = (2, 1)$.
}
\end{ex}

\begin{thm}\label{thm:connected of BP(la)}
The crystal $\ms{B}(P(\la))$ is connected.
\end{thm}
\pf We give a proof when $\mf{g}=\mf{b}_{m|n}$. The proofs for the other cases are similar.

Let us identify $\bfF^{(\bf c)}\in \ms{B}(\mc{N})$ with ${\bf c} =(c_{ij})_{1\le i\le j\le m+n}\in M(\mf{u})$ by \eqref{eq:nilradical roots}, where we write $c_{ij} = c_{(i,j)}$ for simplicity.
Let $\mathbb{O}$ denote the zero vector in $M(\mf{u})$.

Let $b={\bf c}\ot b'\in \ms{B}(\mc{N})\ot \ms{B}(V_{\mf l}(\la))$ be given. We claim that $b$ is connected to $\mathbb{O}\ot v_\la$. We use induction on ${\rm ht}(-{\rm wt}({\bf c}))$.

Suppose that ${\rm ht}(-{\rm wt}({\bf c}))=0$ or ${\bf c}=\mathbb{O}$. Then $b$ is connected to $\mathbb{O}\ot v_\la$ by Proposition \ref{prop:tensor product thm} since $\ms{B}(V_{\mf l}(\la))\cong SST_{m|n}(\la)$ is connected by Theorem \ref{thm:Blambda cnn}.

By Proposition \ref{prop:semisimplicity} and Lemma \ref{lem:Plambda poly}, $P(\la)$ is a direct sum of $V_{\mf l}(\mu)$'s. This implies that any connected component of $\ms{B}(P(\la))$ as an $\mf{l}$-crystal is isomorphic to $\ms{B}(V_{\mf l}(\mu))$ for some $\mu\in \cP_{m|n}$.
In other words, $b$ is $\mf l$-crystal equivalent to $T\in SST_{m|n}(\mu)$ for some $\mu\in\cP_{m|n}$ and $T$, that is, there is an isomorphism of $\mf{l}$-crystals from $C(b)$ to $SST_{m|n}(\mu)$ sending $b$ to $T$, where $C(b)$ is the connected component of $b\in \ms{B}(P(\la))$ with respect to $\te_i$ and $\tf_i$ for $i\in I\setminus\{0\}$. So we may assume that $b$ is $\mf l$-crystal equivalent to $H_\mu$, the genuine highest weight vector of shape $\mu$.

Since $\te_i b=0$ for $1\le i\le m-1$, it is not difficult to see from Proposition \ref{prop:tensor product thm} that $c_{11}\ge c_{22}\ge \dots\ge c_{mm} \ge 0$ and $c_{ij}=0$ for $1\le i<j\le m$. If $c_{11}>0$, then $\te_0 b = (\te_0 {\bf c})\ot b' \neq 0$ with ${\rm ht}(-{\rm wt}(\te_0 {\bf c}))<{\rm ht}(-{\rm wt}({\bf c}))$ by \eqref{eq:action of crystal f0}. By induction hypothesis, $\te_0 b$ is connected to $\mathbb{O}\ot v_\la$.
So we may assume
\begin{equation}\label{eq:positive part of c}
 c_{ij}=0 \quad  (1\le i\le j\le m).
\end{equation}

In this case, we claim that ${\bf c}=\mathbb{O}$.
Let $\mc{N}_1$ be the subalgebra generated by $\bfF_{\beta}$ for $\beta=(1,j)$ for $1\le j\le m+n$, and $\mc{N}_{\ge 2}$ the subalgebra generated by $\bfF_{\beta}$ for $\beta=(i,j)$ for $2 \le i \le j\le m+n$, which are well-defined by Proposition \ref{prop:commutation rel}.

Let $\U(\mf{l}')$ be the subalgebra of $\U(\mf{l})$ generated by $k_\mu, e_i, f_i$ for $\mu\in P'$ and $i\in I\setminus\{0,1\}$, where $P'=\bigoplus_{a\in \I\setminus\{1\}}\Z\de_a$. (We may regard $\mf l'$ as a subalgebra of $\mf l$ isomorphic to $\mf{gl}(m-1|n)$.) By Propositions \ref{prop:e dot} -- \ref{prop:poly alg}, we see that $\mc{N}_\ast$ ($\ast= 1, {\ge 2}$) is a $\U(\mf{l}')$-submodules of $\mc{N}$, and
\begin{equation}\label{eq:tensor decomp of N}
 \mc{N} \cong \mc{N}_1 \ot \mc{N}_{\ge 2},
\end{equation}
as a $\U(\mf{l}')$-module.
If we put $\ms{L}(\mc{N}_\ast)=\ms{L}(\mc{N})\cap \mc{N}_\ast$ and $\ms{B}(\mc{N}_\ast)=\ms{B}(\mc{N})\cap \mc{N}_\ast$, then $(\ms{L}(\mc{N}_\ast),\ms{B}(\mc{N}_\ast))$ is a crystal base of $\mc{N}_\ast$ as a $\U(\mf{l}')$-module, and
\begin{equation*}
 \ms{B}(\mc{N}) \cong \ms{B}(\mc{N}_1)\ot \ms{B}(\mc{N}_{\ge 2}),
\end{equation*}
as an $\mf l'$-crystal by \eqref{eq:tensor decomp of N}.
Hence we may write ${\bf c}={\bf c}_1\ot {\bf c}_{\ge 2}$, where ${\bf c}_1=(c_{11},\dots,c_{1\,m+n})$ and ${\bf c}_{\ge 2}=(c_{ij})_{2\le i\le j\le m+n}$.
Since $\mc{N}_1$ is isomorphic to a direct sum of $V_{\mf l'}((r))$ for $r \ge 0$, we may identify ${\bf c}_1\in \ms{B}(\mc{N}_1)$ with $T'_1\in SST_{m-1|n}((r))$, an $(m-1|n)$-hook semistandard tableau of shape $(r)$ with even (resp.~odd) letters in $\{2,\dots,m\}$ (resp.~$\{m+1,\dots,m+n\}$), where $r=c_{12}+\dots+c_{1\,m+n}$ and $c_{1j}$ ($2\le j\le m+n$) is the number of occurrences of $j$ in $T'_1$. 

Let $T'_2$ and $T'_3$ be the $(m-1|n)$-hook semistandard tableaux which are $\mf l'$-equivalent to ${\bf c}_{\ge 2}$ and $b'$.
Since $b$ is ${\mf l}'$-crystal equivalent to $H_{\mu_{\ge 2}}$, where $\mu_{\ge 2}=(\mu_2,\mu_3,\dots)$, the tableau $P$ obtained by inserting $T'_2$ and then $T'_3$ into $T'_1$ is $H_{\mu_{\ge 2}}$ (see \cite[Section 4]{BKK} for more details on insertion of tableaux). Considering the bumping process in the insertion of $T'_2$ and $T'_3$,  we conclude that $c_{1j}=0$ for $m\le j\le m+n$, since $P$ can not be $H_\mu$ otherwise. Hence we have ${\bf c}_1=(0,\dots,0)$ by \eqref{eq:positive part of c}.

Next, we consider ${\bf c}_2:=(c_{2j})_{2\le j\le m+n}$. Note that $c_{2j}=0$ for $2\le j\le m$ by \eqref{eq:positive part of c}. If $c_{2j}=1$ for some $m\le j\le m+n$, then we have $\te_1 b = (\te_1 {\bf c})\ot b'\neq 0$, which is a contradiction. So we have ${\bf c}_2=(0,\dots,0)$. Repeating similar argument, we obtain
\begin{equation}\label{eq:non-negative part of c}
 {\bf c}_i:=(c_{ij})_{i\le j\le m+n}=(0,\dots,0) \quad  (1\le i\le m).
\end{equation}

Finally, let $T_1$ and $T_2$ be the $(m|n)$-hook semistandard tableaux, which are $\mf l$-crystal equivalent to ${\bf c}$ and $b'$, respectively. The tableau obtained by inserting $T_2$ into $T_1$ is $H_{\mu}$. On the other hand, ${\rm wt}({\bf c})\in \Z_{\ge 0}\de_{m+1}\oplus\dots\oplus\Z_{\ge 0}\de_{m+n}$ by \eqref{eq:non-negative part of c}. Again by considering the bumping process in the insertion of $T_2$ into $T_1$ whose resulting tableau should be $H_\mu$, we conclude that ${\bf c}=\mathbb{O}$. This completes the proof.
\qed

\begin{cor}
 The crystal $\ms{B}(\mc{N})$ is connected.
\end{cor}
\pf It follows from $\ms{B}(\mc{N})\cong \ms{B}(P(0))$.
\qed

\begin{cor}
We have
\begin{equation*}
\begin{split}
\ms{L}(P(\la))&=\sum_{r \ge 0, i_1,\dots,i_r \in I}A_0\td{x}_{i_1}\dots\td{x}_{i_r}w_\la, \\ 
\ms{B}(P(\la))&=\left\{\,\pm\td{x}_{i_1}\dots\td{x}_{i_r}w_\la \!\!\!\pmod{q\ms{L}(P(\la))}\,|\,r\ge 0, i_1,\dots,i_r\in I\,\right\}\setminus\{0\}, 
\end{split}
\end{equation*}
where $x = e,f$ for each $i_k$.
\end{cor}
\pf If follows from Theorem \ref{thm:connected of BP(la)} and Nakayama's lemma.
\qed


\appendix

\section{Proof of Proposition \ref{prop:i-comp submodule-3}}\label{appendix-1}

Suppose that  $\mf{g}=\mf{b}_{m|n}$ with $i\ge m$ or $\mf{g}=\mf{d}_{m|n}$ with $i> m$. In this case, we have $S=\{\,(i,i), (i,i+1), (i+1,i+1)\,\}$. 

Let $L=\ms{L}\left(\mc{N}_{\Delta}(i)\right)$ and 
$B=\ms{B}\left(\mc{N}_{\Delta}(i)\right)$, where $\mc{N}_{\Delta}(i)$ is given in \eqref{eq:notations for N(i)'s}.
We first find all weight vectors $E$ in $L$ such that $e_i\cdot E=0$ and $E\equiv b \pmod{qL}$ for some $b\in B$. Let $B_0$ be the set of such $b$'s and denote the corresponding maximal vector by $E_b$. Then we show that $\tf_i^s E_b\in L$ and $\tf_i^s E_b \in B \cup\{0\} \pmod{qL}$ for $b\in B_0$ and $s\ge 0$, and every $b\in B$ is obtained in this way. This shows that $L=\bigoplus_{b\in B_0,s\in\Z_{\ge 0}}A_0 \tf_i^sE_b$ and  $(L,B)$ is a crystal base of $\mc{N}_{\Delta}(i)$, which also implies Proposition \ref{prop:i-comp submodule-3} for $\mc{N}_{\Delta}(i;d)$ for $d\in\Z_{\ge 0}$.

For a statement $P$, $\de_P$ means $1$ and $0$, when $P$ is true and false, respectively.


\subsection{$\mf{g}=\mf{b}_{m|n}$ with $i=m$}
In this case, we consider $\mc{N}_{\Delta}(m)$ corresponding to
\begin{equation*}
	\begin{tikzpicture}[baseline=(current  bounding  box.center), every node/.style={scale=1.1}]
		\node (0) at (0,0) {$\underset{(m,m)}{\overset{\tiny }{\bigcirc}}$};
		\node (1) at (0.9,0.9) {$\underset{(m,m+1)}{\overset{\tiny }{\scalebox{1.2}{$\otimes$}}}$};
		\node (2) at (1.8,0) {$\underset{(m+1,m+1)}{{\overset{\tiny }{{\scalebox{2.3}{$\bullet$}}}}}$};
	\end{tikzpicture}
\end{equation*}
(cf.~ Example \ref{ex: roots of u}).
\noindent
For simplicity, we put
\begin{equation*}
\begin{split}
& A=\bff_{(m,m)},\quad B=\bff_{(m,m+1)},\quad C=\bff_{(m+1,m+1)},\ \
 e=e_m,\quad f=f_m,\quad  k=k_m,\\
& 	L=\ms{L}\left(\mc{N}_{\Delta}(m)\right) =\bigoplus_{a,c\in \Z_{\ge 0}, b \in \{0,1 \} }{A}_0\, q^{{\bf z}(c)} A^{(a)}B^{(b)}C^{(c)}, 
\end{split}
\end{equation*}
where ${\bf z}(c)=-\frac{c(c-1)}{2}$.
Let $\te=\te_m$ and $\tf=\tf_m$ be given in \eqref{eq:crystal operator m}.

Let us first find the weight vectors in ${\rm Ker} \ e \cap L$.
By Propositions \ref{prop:commutation rel}, \ref{prop:e dot}, and \ref{prop:f dot}, we have
\begin{gather*}
	e \cdot A = 0 , \quad e \cdot B = \frac{q^{-2}-1}{[2]} \, A^2, \quad e \cdot C = A, \\
	f \cdot A = C , \quad f \cdot B = \frac{-q^{2}-1}{[2]} \, C^2, \quad f \cdot C = 0, \\
	k \cdot A = q^{2} A , \quad k \cdot B = -q^4 B, \quad k \cdot C = -q^2 \, C, \\
	BA = q^{-2} AB, \quad CB = -q^2 BC, \quad CA-AC = [2] B, \quad B^2=0,
\end{gather*}
and then by induction on $s\ge 1$,
\begin{equation} \label{eq: ef divided m}
	\begin{split}
	e \cdot C^{(s)} & = (-1)^{s-1} q^{-s+1} AC^{(s-1)} + (-1)^s \frac{[2]}{ \{ 2\} } BC^{(s-2)}\de_{s\ge 2}, \\
	f \cdot A^{(s)} & = q^{s-1} A^{(s-1)}C + A^{(s-2)}B \, \de_{s \ge 2},
	\end{split}
\end{equation}
where $A^{(a)}=A^a/[a]!$ and $C^{(c)}=C^c/\{c\}!$ for $a, c \in \Z_{\ge 0}$.

It is clear that $E_{a,0}:=A^{(a)} \in {\rm Ker} \, e \cap L$.
Since
\begin{equation*}
	\begin{split}
		e \cdot A^{(a)} C^{(c)} & = (-1)^{c-1} q^{-c+1} [a+1] A^{(a+1)} C^{(c-1)} \de_{c\ge 1} + (-1)^{c} \frac{[2]}{\{2\}} A^{(a)}BC^{(c-2)}\de_{c \ge 2}, \\
		e \cdot A^{(a)} B C^{(c)} & = (-1)^{c} q^{-2c} \frac{(q^{-2}-1)}{[2]} [a+2][a+1] A^{(a+2)}C^{(c)} \\
		& \hspace{4.6cm} + (-1)^{c-1}q^{-(c+1)}[a+1] A^{(a+1)}BC^{(c-1)} \de_{c \ge 1},
	\end{split}
\end{equation*}
by \eqref{eq: ef divided m}, 
we also have 
\begin{equation*}
E_{a+1,c+1} := q^{{\bf z}(c+1)}
\left(
A^{(a+1)}C^{(c+1)} - \frac{[2]q^{c+1}}{\{2\}[a+1]}A^{(a)}BC^{(c)} \right)\in {\rm Ker}\,e\cap L
\end{equation*}
for $a, c \ge 0$, and then
\begin{equation*}
	E_{a+1,c+1} \equiv q^{{\bf z}(c+1)} A^{(a+1)}C^{(c+1)} \quad ({\rm mod} \ qL),
\end{equation*}
where $\frac{[2]q^{c+1}}{\{2\}[a+1]} \in q^{a+c+1}(1+q{A}_0)$ and ${\bf z}(c+1) + (a+c+1) = {\bf z}(c) +(a+1) > {\bf z}(c) $.
We can check that any weight vector in $\mc{N}_{\Delta}(m)$ belongs to  ${\rm Ker}\, e$ if and only if it is either $E_{a,0}$ or $E_{a+1,c+1}$ ($a,c\ge 0$) up to scalar multiplication.

Next, we apply $\tf=f$ to a maximal vector.
By \eqref{eq: ef divided m}, we have $f \cdot E_{a,0}\in L$ and
\begin{equation*}
	f \cdot E_{a,0} \equiv 
\begin{cases}
 A^{(a-1)}B  & \text{if $a\ge 2$,}\\
 C & \text{if $a=1$,} \\
  0  & \text{if $a=0$,}
\end{cases} \quad ({\rm mod} \ qL).
\end{equation*}
Also, we have
\begin{equation*}
	f \cdot E_{a+1,c+1} = q^{a+{\bf z}(c+1)} \{c+2\} X A^{(a)}C^{(c+2)} + q^{{\bf z}(c+1)} X A^{(a-1)} B C^{(c+1)} \de_{a \ge 1},
\end{equation*}
where $X=1+\frac{q^{a+c+2}[2]\{c+1\}}{\{2\}[a+1]} \in 1+q{A}_0$, which implies that $f \cdot E_{a+1,c+1}\in L$ and
\begin{equation*}
	f \cdot E_{a+1,c+1} \equiv 	
\begin{cases}
 q^{{\bf z}(c+1)} A^{(a-1)} B C^{(c+1)}  & \text{if $a\ge 1$,}\\
q^{{\bf z}(c+2)} C^{(c+2)}  & \text{if $a=0$,}
\end{cases}	\quad ({\rm mod} \ qL),
\end{equation*}
since  $a-{\bf z}(c+1) - c-1 = a- {\bf z}(c+2)$.
Note that the above computation yields to {\it Case 1} for $\mf{g} = \mf{b}_{m|n}$ and $i = m$ in Remark \ref{rem:crystal operator on i-comp-2}.

\subsection{$\mf{g}=\mf{b}_{m|n}$ with $i>m$}\label{subsec:b i>m}
In this case, we consider $\mc{N}_{\Delta}(i)$ corresponding to
\begin{equation*}
	\begin{tikzpicture}[baseline=(current  bounding  box.center), every node/.style={scale=1.1}]
		\node (0) at (0,0) {$\underset{(i,i)}{\overset{\tiny }{\scalebox{2.3}{$\bullet$}}}$};
		\node (1) at (0.9,0.9) {$\underset{(i,i+1)}{\overset{\tiny }{\bigcirc}}$};
		\node (2) at (1.8,0) {$\underset{(i+1,i+1)}{{\overset{\tiny }{{\scalebox{2.3}{$\bullet$}}}}}$};
	\end{tikzpicture}
\end{equation*}
\noindent
Put
\begin{equation}\label{eq:notations for b i>m}
\begin{split}
 & A=\bff_{(i,i)},\quad B=\bff_{(i,i+1)},\quad C=\bff_{(i+1,i+1)},\ \
 e=e_i,\quad f=f_i,\quad  k=k_i,\\
 &L=\ms{L}\left(\mc{N}_{\Delta}(i)\right) =\bigoplus_{a,b,c\in \Z_{\ge 0}}{A}_0\, q^{{\bf z}(a,b,c)} A^{(a)}B^{(b)}C^{(c)},
\end{split} 
\end{equation}
where ${\bf z}(a,b,c)=-\frac{a(a-1)}{2}-b(b-1)-\frac{c(c-1)}{2}$.
Let $\te=\te_i, \tf=\tf_i$ be given in \eqref{eq:crystal operator >m}, that is,
\begin{equation*}
	\tilde{e}_iu=\sum_{k\geq1}q^{2(-l_k+2k-1)}f_i^{(k-1)}u_k,\quad 
	\tilde{f}_iu=\sum_{k\geq 0}q^{2(l_k-2k-1)}f_i^{(k+1)}u_k,
\end{equation*}
where $u=\sum_{k \geq 0} f_i^{(k)}u_k$ with $e_i\cdot u_k=0$ and $l_k= - ({\rm wt}(u_k)|\alpha_i)/2$.
In this case, we have
\begin{gather*}
e \cdot A = 0 , \quad e \cdot B = -q \, A^2, \quad e \cdot C = A, \\
f \cdot A = C , \quad f \cdot B = -q \, C^2, \quad f \cdot C = 0, \\
k \cdot A = -q^{-2} A , \quad k \cdot B = B, \quad k \cdot C = -q^2 \, C, \\
BA = -q^2 AB, \quad CB = -q^2 BC, \quad CA-AC = [2] B,
\end{gather*}
and then by induction on $s\ge 1$,
\begin{equation} \label{eq: ef divided i>m;b}
\begin{split}
e \cdot C^{(s)} & = (-1)^{s-1} q^{-s+1} AC^{(s-1)} + (-1)^s \frac{[2]}{ \{ 2\} } BC^{(s-2)} \, \de_{s\ge 2}, \\
e \cdot B^{(s)} & = -q^{2s-1} A^2 B^{(s-1)}, \\
f \cdot A^{(s)} & = (-1)^{s-1} q^{-s+1} A^{(s-1)}C + (-1)^s \frac{[2]}{ \{ 2\} } A^{(s-2)}B \, \de_{s\ge 2}, \\
f \cdot B^{(s)} & = -q^{2s-1} B^{(s-1)} C^2,
\end{split}
\end{equation}
where $A^{(a)}=A^a/\{a\}!$, $B^{(b)}=B^b/[b]_{q^2}!$, and $C^{(c)}=C^c/\{c\}!$ for $a, b, c \in \Z_{\ge 0}$.

By \eqref{eq: ef divided i>m;b}, it is straightforward to see that a weight vector in $\mc{N}_{\Delta}(i)$ belongs to ${\rm Ker} \ e$ if and only if it is of the following form:
\begin{equation}\label{eq:maximal;b,i>m}
 E_{a, c} := \sum_{r=0}^{c} X_r A^{(a-r)}B^{(r)}C^{(c-r)},
\end{equation}
for $a \ge c$ such that $X_r = q^{c-3r+1} \frac{1+(-q^2)^r}{(1-q^2)\{a-r+1\}} X_{r-1}$ ($1 \le r \le c$), and hence
\begin{equation*}
X_r = q^{rc-\frac{3r^2+r}{2}} \frac{(1+(-q^2)) \cdots (1+(-q^2)^r)}{(1-q^2)^r \{a\} \cdots \{a-r+1\}} X_0
\end{equation*}
with some $X_0 \in \Bbbk^\times$.
 
For $f\in \Bbbk^\times$, let ${\bf d}(f)$ be a maximal integer $d$ such that $f \in q^{d} {A}_0$. 
%
Let $d_r={\bf d}(X_r)$. Then $d_r=d_0+ar+cr-2r^2-r$ for $1\le r\le c$.
If we put $X_0 = q^{-\frac{a(a-1)}{2}-\frac{c(c-1)}{2}+c}$, then 
\begin{equation}\label{eq:d_r}
\begin{split}
d_r&=-\frac{(a-r)(a-r-1)}{2}-r(r-1)-\frac{(c-r)(c-r-1)}{2}+(c-r)\\
   &= {\bf z}(a-r,r,c-r) + c-r, 
\end{split}
\end{equation}
and hence $X_r \in \pm q^{d_r} (1+q {A}_0)$. 
This implies that
$E_{a,c} \in L$ and 
\begin{equation*}
 E_{a,c} \equiv \pm q^{{\bf z}(a-c,c,0)}A^{(a-c)}B^{(c)}\ ({\rm mod} \ qL ).
\end{equation*}
	
Now, let us compute $\tf^s E_{a,c}$ ($s\ge 0$).
Since $2(\alpha_i | {\rm wt}(E_{a,c}))/ (\alpha_i | \alpha_i) = a-c $, we have by \eqref{eq: ef divided i>m;b}
%
%
\begin{equation*}
\tf^s E_{a, c} = \sum_{r=0}^{{\rm min}(a-s,c+s)} X_{s,r} A^{(a-s-r)}B^{(r)}C^{(c+s-r)}\quad (0 \le s \le a-c)
\end{equation*}
for some $X_{s,r} \in \Bbbk$.
Let $d_{s,r}={\bf d}(X_{s,r})$.
Note that $X_{0,r}=X_r$ and $d_{0,r}=d_r $. 

The rest of this subsection is devoted to proving 
\begin{equation}\label{eq:formula for d;b}
\begin{split}
&d_{s,r}={\bf z}(a-s-r,r,c+s-r)+
\begin{cases}
(r-c)^2 & \text{if $r\ge c-1$},\\
c-r & \text{if $r \le c$},
\end{cases}\\
&X_{s,r} \in \pm q^{d_{s,r}} (1+q {A}_0), 
\end{split}
\end{equation}
which implies that $\tf^s E_{a, c}\in L$ for $0 \le s \le a-c$, and 
\begin{equation*}
\tf^s E_{a, c} \equiv \pm q^{{\bf z}(a-s-c,c,s)}A^{(a-s-c)}B^{(c)}C^{(s)}\ ({\rm mod} \ qL ).
\end{equation*}
Note that the above computation yields to {\it Case 1} for $\mf{g} = \mf{b}_{m|n}$ and $i > m$ in Remark \ref{rem:crystal operator on i-comp-2}.

To prove \eqref{eq:formula for d;b}, we consider recursive formulas for $X_{s,r}$.
From $\tf^{s+1}E_{a,c}=\frac{q^{2(a-c-2s-1)}}{{[s+1]_{q^2}}}f \cdot \tf^s E_{a,c}$ and \eqref{eq: ef divided i>m;b}, we obtain
\begin{equation} \label{eq: app f formula;b}
X_{s+1,r} = \frac{q^{2(a-c-2s-1)}}{[s+1]_{q^2}} (x_1X_{s,r}+x_2X_{s,r-1}+x_3X_{s,r+1}),
\end{equation}
where
\begin{equation*}
\begin{split}
& x_1 = 	(-1)^{a-s-1} q^{-a+s+3r+1} \{c+s-r+1\} \ \de_{0\le r \le {\rm min}(a-s-1, c+s)}, \\
& x_2 = (-1)^{a-s-r+1} \frac{[2]}{\{2\}} [r]_{q^2} \ \de_{1\le r \le {\rm min}(a-s-1, c+s+1)}, \\
& x_3 = (-1)^{a-s-r} q^{-2a+2s+4r+3} \{c+s-r\}\{c+s-r+1\} \ \de_{0\le r \le {\rm min}(a-s-1, c+s-1)}.
\end{split}
\end{equation*}
On the other hand, from $ef^{(s)}=f^{(s)}e+f^{(s-1)}\frac{q^{2(s-1)}k-q^{-2(s-1)}k^{-1}}{q^2-q^{-2}}$, we have $\tf^{s-1} E_{a,c} = (-1)^{a+c+1}\frac{q^{2(-a+c+2s-1)}}{[a-c-s+1]_{q^2}} \ e \cdot \tf^s E_{a,c} $, which implies
\begin{equation} \label{eq: app e formula;b}
\begin{split}
X_{s-1,r} = (-1)^{a+c+1} \frac{q^{2(-a+c+2s-1)}}{[a-c-s+1]_{q^2}} (y_1 X_{s,r} + y_2 X_{s,r-1}+ y_3X_{s,r+1}),
\end{split}
\end{equation}
where
\begin{equation*}
\begin{split}
& y_1 = 	(-1)^{c+s-1} q^{-c-s+3r+1} \{a-s-r+1\} \ \de_{0\le r \le {\rm min}(a-s, c+s-1)}, \\
& y_2 = (-1)^{c+s-r+1} \frac{[2]}{\{2\}} [r]_{q^2} \ \de_{1\le r \le {\rm min}(a-s+1, c+s-1)}, \\
& y_3 = (-1)^{c+s-r} q^{-2c-2s+4r+3} \{a-s-r\}\{a-s-r+1\} \ \de_{0\le r \le {\rm min}(a-s-1, c+s-1)}.
\end{split}
\end{equation*}
Note that the formula \eqref{eq: app e formula;b} can be recovered from \eqref{eq: app f formula;b} (and vice versa) up to multiplication $(-1)^{a+c+1}$  by replacing $X_{l,r}$ with $X_{\pi(l),r}$ where $\pi: l \longmapsto a-c-l$.

\begin{lem}\label{lem:symmetry of recursive formulas}
We have $X_{s,r}=\pm X_{a-c-s,r}$ and $d_{s,r}=d_{a-c-s,r}$.
\end{lem}
\pf We have $F_{a,c}:=\sum_{r=0}^{c} X_r A^{(c-r)}B^{(r)}C^{(a-r)} \in {\rm Ker} \, f$, where $X_r$ is as in \eqref{eq:maximal;b,i>m}.
We have $\tf^{a-c} E_{a,c} = Z F_{a,c}$ for some $Z \in \Bbbk^\times$, and hence $Z X_{0,r} = X_{a-c,r}$.
Since the formulas \eqref{eq: app e formula;b} and \eqref{eq: app f formula;b} coincide under $\pi$ up to sign, 
we have $Z X_{s,r} = (-1)^{s(a+c+1)} X_{a-c-s,r}$, especially $ Z X_{a-c,r} = (-1)^{(a-c)(a+c+1)} X_{0,r} = X_{0,r}$.
Therefore we have $X_{0,r}=Z X_{a-c,r}=Z^2 X_{0,r}$, that is, $Z=\pm 1$.
This implies $d_{s,r}=d_{a-c-s,r}$.
\qed
\smallskip	

By Lemma \ref{lem:symmetry of recursive formulas}, it is enough to compute $d_{s,r}$ for $s \le \frac{a-c}{2}$, and $0 \le r \le {\rm min} (a-s,c+s)=c+s$.
We prove \eqref{eq:formula for d;b} in the following steps:
\smallskip

{\em Step 1}. Consider $d_{0,r}$ for $0\le r\le c$. Then it follows from \eqref{eq:d_r}. \smallskip

{\em Step 2}. Consider $d_{s,c+s}$. We use induction on $s$. By \eqref{eq: app f formula;b} we have
	\begin{equation*}
		d_{s+1,c+s+1}=d_{s,c+s}+2a-4c-4s-2
	\end{equation*}
	if $s+1 \le \frac{a-c}{2}$, and $d_{0,c}={\bf z}(a-c,c,0)$ by {\em Step 1}.
Thus we have $d_{s+1,c+s+1}={\bf z}(a-c-2s-2,c+s+1,0)+(s+1)^2$ and $X_{s+1,c+s+1} \in \pm q^{d_{s+1,c+s+1}} (1+q {A}_0)$ by induction hypothesis.
	
Next, consider $d_{s,c+s-1}$ ($c+s\ge 1$). We also use induction on $s$.  It is true for $s=0$ since $d_{0,c-1}={\bf z}(a-c+1,c-1,1)+1$.
By \eqref{eq: app f formula;b} we have
	\begin{equation*}
		d_{s+1,c+s} = 2a-2c-2s-2 + {\rm min} (d_1, d_2)
	\end{equation*}
	when $d_1 \ne d_2$, where $d_1=d_{s,c+s}-a+3c+4s+1$ and  $d_2=d_{s,c+s-1}-2c-2s+2$.
Here we use the fact that if $f, g_1, g_2  \in \Bbbk^\times $ satisfy (1) $f= g_1 +g_2 $,
(2) $d(g_1)<d(g_2)$, (3) $g_1  \in \pm q^{d(g_1)} (1+q{A}_0)$, then $d(f)=d(g_1)$ and $f \in \pm q^{d(g_1)} (1+q{A}_0)$.
If $(s,c)\neq (0,0)$, then we have
$d_1 - d_2 = 2c+4s >0 $ by induction hypothesis, and hence $d_{s+1,c+s}=2a-4c-4s+d_{s,c+s-1}$.

If $c=s=0$, then we have $d_{1,0}=-\frac{(a-1)(a-2)}{2}={\bf z}(a-1,0,1)+(1-1)^2$.
	In both cases, we obtain $d_{s+1,c+s}={\bf z}(a-c-2s-1,c+s,1)+s^2$ and $X_{s+1,c+s} \in \pm q^{d_{s+1,c+s}} (1+q {A}_0)$. The induction completes.\smallskip
	
{\em Step 3}. Consider $d_{s,r}$ for $r\ge c-1$. We use induction on $r$ in decreasing order.
By {\em Step 2}, \eqref{eq:formula for d;b} holds for $r=c+s, c+s-1$.
By \eqref{eq: app e formula;b}, if $r-1 \ge c-1$, then we have
	\begin{equation}\label{eq:recursive relation;b step 3}
		d_{s,r-1} = 2(r-1) + {\rm min} (d_1, d_2, d_3)
	\end{equation}
	if at least one of $d_1, d_2, d_3$ is strictly smaller than the others, where
	$d_1=d_{s,r}-a-c+4r+1$, $d_2=d_{s,r+1}-2a-2c+6r+4$ and $d_3=d_{s-1,r}-2s+2$ (by similar arguments as in {\em Step 2}).
	By induction hypothesis, we have
	$d_1-d_3=-2c+4r > 0$, and $d_2-d_3=2-4c+4r\ge2>0$, that is,  $d_3<d_1,d_2$. Hence we have $d_{s,r-1}=2r-2s+d_{s-1,r}$.
	Again by induction hypothesis, we have $d_{s,r-1}={\bf z}(a-s-r+1,r-1,c+s-r+1)+(r-c-1)^2$ and $X_{s,r} \in \pm q^{d_{s,r}} (1+q {A}_0)$. The induction completes. \smallskip
	
{\em Step 4}. Consider $d_{s,r}$ for $r < c-1$.  
We also use induction on $r$ in decreasing order.
By {\em Step 3}, \eqref{eq:formula for d;b} holds for $r=c, c-1$ since $(r-c)^2=c-r$.

We also use \eqref{eq:recursive relation;b step 3}.
By induction hypothesis, we have $d_1-d_2 =2r \ge 2 >0$ and $d_3-d_2=2c-2r>0$, that is, $d_2<d_1,d_3$.
	We have $d_{s,r-1}=-2a-2c+8r+2+d_{s,r+1}={\bf z}(a-s-r+1,r-1,c+s-r+1)+c-r+1$ and $X_{s,r} \in \pm q^{d_{s,r}} (1+q {A}_0)$. The induction completes.\smallskip
	
By {\em Steps 1--4} and Lemma \ref{lem:symmetry of recursive formulas}, \eqref{eq:formula for d;b} holds for all $s,r$.

\subsection{$\mf{g}=\mf{d}_{m|n}$ with $i>m$}
In this case, we consider $\mc{N}_{\Delta}(i)$ corresponding to
\begin{equation*}
	\begin{tikzpicture}[baseline=(current  bounding  box.center), every node/.style={scale=1.1}]
		\node (0) at (0,0) {$\underset{(i,i)}{\overset{\tiny }{\bigcirc}}$};
		\node (1) at (0.9,0.9) {$\underset{(i,i+1)}{\overset{\tiny }{\bigcirc}}$};
		\node (2) at (1.8,0) {$\underset{(i+1,i+1)}{{\overset{\tiny }{{\bigcirc}}}}$};
	\end{tikzpicture}
\end{equation*}
\noindent
The proof is similar to \ref{subsec:b i>m}.
We keep the same notations as in \eqref{eq:notations for b i>m} except
${\bf z}(a,b,c)=-a^2-\frac{b(b-1)}{2}-c^2$.
In this case, $\te=\te_i$ and $\tf=\tf_i$ are given by
	\begin{equation*}
		\tilde{e}_iu=\sum_{k\geq1}q^{-l_k+2k-1}f_i^{(k-1)}u_k,\quad 
		\tilde{f}_iu=\sum_{k\geq 0}q^{l_k-2k-1}f_i^{(k+1)}u_k,
	\end{equation*}
	where $l_k= - ({\rm wt}(u_k)|\alpha_i)$.
As to the actions of $\U(\mf{l})_i$, we have
	\begin{gather*}
		e \cdot A = 0 , \quad e \cdot B = [2] \, A, \quad e \cdot C = -B, \\
		f \cdot A = -B , \quad f \cdot B = [2] \, C, \quad f \cdot C = 0, \\
		k \cdot A = q^{-2} A , \quad k \cdot B = B, \quad k \cdot C = q^2 \, C, \\
		BA = q^2 AB, \quad CB = q^2 BC, \quad CA-AC = \frac{(1-q^2)}{[2]} B^{2},
	\end{gather*}
and for $s\ge 1$,
	\begin{equation} \label{eq: ef divided type c}
		\begin{split}
		e \cdot C^{(s)} & = -BC^{(s-1)}, \\
		e \cdot B^{(s)} & = [2] q^{s-1} A B^{(s-1)}, \\
		f \cdot A^{(s)} & = -A^{(s-1)}B, \\
		f \cdot B^{(s)} & = [2] q^{s-1} B^{(s-1)} C,
		\end{split}
	\end{equation}
where $A^{(a)}=A^a/[a]_{q^2}!$, $B^{(b)}=B^b/[b]!$ and $C^{(c)}=C^c/[c]_{q^2}!$ for $a, b, c \in \Z_{\ge 0}$.

By \eqref{eq: ef divided type c}, we can show that a weight vector in $\mc{N}_{\Delta}(i)$ belongs to ${\rm Ker} \, e$ if and only if it is of the following form:
\begin{equation*}
 E_{a, c} := \sum_{r=0}^{c} X_r A^{(a-r)}B^{(2r)}C^{(c-r)},
\end{equation*}
for $a\ge c$ such that $X_r = q^{2c-4r+1} \frac{[2r-1]}{[2a-2r+2]} X_{r-1}$ ($1 \le r \le c$), and hence
\begin{equation*}
 X_r = q^{2cr-2r^2-r} \frac{[2r-1][2r-3] \cdots [1]}{[2a][2a-2] \cdots [2a-2r+2]} X_0 \quad (1 \le r \le c),
\end{equation*}
with some $X_0 \in \Bbbk^\times$.

Put $d_r={\bf d}(X_r)$.
Then $d_r=d_0+(2a+2c)r-4r^2$.
Letting $X_0 = q^{-a^2-c^2+c}$, we have 
\begin{equation}\label{eq:d_r;d}
	d_r= {\bf z}(a-r,2r,c-r) + c-r,\quad X_r \in \pm q^{d_r} (1+q \mathbb{A}_0)
\end{equation}
so that $E_{a,c} \in L$ and 
\begin{equation*}
 E_{a,c} \equiv \pm q^{z(a-c,2c,0)}A^{(a-c)}B^{(2c)}\ ({\rm mod} \ qL ).
\end{equation*}
	
Let us compute $\tf^s E_{a,c}$ $(s \ge 0)$.
	Since $2(\alpha_i | {\rm wt}(E_{a,c}))/(\alpha_i | \alpha_i) = 2(a-c) $, we have for $0 \le s \le 2(a-c)$,
	\begin{equation*}
		\tf^{s} E_{a, c} =
		\begin{cases*}
			\sum_{r=0}^{{\rm min}(a-t,c+t-1)} X_{s,r} A^{(a-t-r)}B^{(2r+1)}C^{(c+t-r-1)} & \text{ if $s=2t-1$\  ($1\le t\le a-c$)}, \\
		 	\sum_{r=0}^{{\rm min}(a-t,c+t)} X_{s,r} A^{(a-t-r)}B^{(2r)}C^{(c+t-r)} & \text{ if $s=2t$\  ($0\le t\le a-c$)},
		\end{cases*}
\end{equation*}
for some $X_{s,r} \in \Bbbk$ by \eqref{eq: ef divided type c}.
	Let $d_{s,r}={\bf d}(X_{s,r})$.

Then it remains to prove	
\begin{equation}\label{eq:formula for d;d}
	\begin{split}
	d_{2t-1,r} & ={\bf z}(a-t-r,2r+1,c+t-r-1)+
	\begin{cases}
			2(r-c)^2+(r-c) & \text{if $r \ge c$},\\
			3(c-r) & \text{if $r \le c$},
	\end{cases}\\			
	d_{2t,r} & ={\bf z}(a-t-r,2r,c+t-r)+
	\begin{cases}
			2(r-c)^2-(r-c) & \text{if $r \ge c$},\\
			(c-r) & \text{if $r \le c$},
	\end{cases}\\
	X_{s,r} & \in \pm q^{d_{s,r}} (1+q {A}_0), 
	\end{split}			
\end{equation}
which implies that $\tf^{s} E_{a, c}\in L$ and
\begin{equation*}
	\tf^{s} E_{a, c} \equiv 
	\begin{cases*}
		\pm q^{{\bf z}(a-c-t,2c+1,t-1)}A^{(a-c-t)}B^{(2c+1)}C^{(t-1)} & \text{ if $s=2t-1$,}\\
		\pm q^{{\bf z}(a-c-t,2c,t)}A^{(a-c-t)}B^{(2c)}C^{(t)}  & \text{ if $s=2t$,}
	\end{cases*}\quad	({\rm mod} \ qL ).	
\end{equation*}
Note that the above computation yields to {\it Case 3} for $\mf{g} = \mf{d}_{m|n}$ and $i > m$ in Remark \ref{rem:crystal operator on i-comp-2}.

As in Section \ref{subsec:b i>m}, we use recursive formulas for $X_{s,r}$.
From $\tf^{s+1}E_{a,c}=\frac{q^{2a-2c-2s-1}}{{[s+1]}}f \cdot \tf^s E_{a,c}$ and \eqref{eq: ef divided type c}, we have 
\begin{equation} \label{eq: app f odd formula}
	\begin{split}
	& X_{2t,r} = \frac{q^{2a-2c-4t+1}}{[2t]} (x_1X_{2t-1,r-1}+x_2X_{2t-1,r}), \\
	& x_1 = -  [2r] \ \de_{1\le r \le {\rm min}(a-t, c+t)}, \\
		& x_2 =  q^{-2a+2t+4r} [2c+2t-2r] \ \de_{0\le r \le {\rm min}(a-t, c+t-1)},\\
	\end{split}
\end{equation}
\begin{equation} \label{eq: app f even formula}
	\begin{split}
	& X_{2t+1,r} = \frac{q^{2a-2c-4t-1}}{[2t+1]} (x'_1X_{2t,r}+x'_2 X_{2t,r+1}),\\
	& x'_1 = -  [2r+1] \ \de_{0\le r \le {\rm min}(a-t-1, c+t)}, \\
	& x'_2 = q^{-2a+2t+4r+3} [2c+2t-2r] \ \de_{0\le r \le {\rm min}(a-t-1, c+t-1)}.
	\end{split}
\end{equation}
Similarly, from $ef^{(s)}=f^{(s)}e+f^{(s-1)}\frac{q^{(s-1)}k-q^{-(s-1)}k^{-1}}{q-q^{-1}}$ and $\tf^{s-1} E_{a,c} = - \frac{q^{-2a+2c+2s-1}}{[2a-2c-s+1]} \ e \cdot \tf^s E_{a,c} $,
we have
\begin{equation} \label{eq: app e odd formula}
	\begin{split}
	& X_{2t-2,r} = - \frac{q^{-2a+2c+4t-3}}{[2a-2c-2t+2]} (y_1X_{2t-1,r-1} + y_2X_{2t-1,r} ), \\
	& y_1 = - [2r] \ \de_{1 \le r \le {\rm min}(a-t+1, c+t-1)}, \\
		& y_2 = q^{-2c-2t+4r+2} [2a-2t-2r+2] \ \de_{0\le r \le {\rm min}(a-t, c+t-1)},
	\end{split}
\end{equation}
\begin{equation} \label{eq: app e even formula}
	\begin{split}
	& X_{2t-1,r} = - \frac{q^{-2a+2c+4t-1}}{[2a-2c-2t+1]} (y'_1X_{2t,r} + y'_2X_{2t,r+1}),\\
	& y'_1 = - [2r+1] \ \de_{0 \le r \le {\rm min}(a-t, c+t-1)}, \\
		& y'_2 =  q^{-2c-2t+4r+3} [2a-2t-2r] \ \de_{0\le r \le {\rm min}(a-t-1, c+t-1)}.
	\end{split}
\end{equation}
Note that the formula in \eqref{eq: app e odd formula} (resp. \eqref{eq: app e even formula}) can be obtained from the formula in \eqref{eq: app f odd formula} (resp. \eqref{eq: app f even formula}) under  $t \longmapsto a-c+1-t$ (resp. $t \longmapsto a-c-t$) up to multiplication by $-1$. The following is an analogue of Lemma \ref{lem:symmetry of recursive formulas}.
\begin{lem}\label{lem:symmetry of recursive formulas;d}
  $X_{s,r}=\pm X_{2a-2c-s,r}$ and $d_{s,r}=d_{2a-2c-s,r}$.
\end{lem}
So it is enough to compute $d_{s,r}$ for $0\le s \le a-c $ and 
\begin{equation*}
\begin{cases}
0 \le r \le {\rm min} (a-t,c+t-1)=c+t-1 & \text{if $s=2t-1$},\\ 
0 \le r \le {\rm min} (a-t,c+t)=c+t & \text{if $s=2t$}.
\end{cases}
\end{equation*}
We prove \eqref{eq:formula for d;d} in the following steps:
\smallskip

{\em Step 1}. Consider $d_{0,r}$ for $0\le r\le c$. Then it follows from \eqref{eq:d_r;d}. 
\smallskip 
	
{\em Step 2}. Consider $d_{s,0}$. We use induction on $s$. We have $d_{0,0}={\bf z}(a,0,c)+c$ by {\em Step 1}. 
By \eqref{eq: app f odd formula} and \eqref{eq: app f even formula}, we have
\begin{equation*}
\begin{cases}
	d_{2t+1,0}=d_{2t,0}+2a+2c-2t-1 & \text{if $ 1 \le 2t+1 \le a-c$},\\
	d_{2t,0}=d_{2t-1,0}-4c-2t+1 & \text{if $2 \le 2t \le a-c$}.
\end{cases}
\end{equation*}
By induction hypothesis, we obtain $d_{2t,0}={\bf z}(a-t,0,c+t)+c$ and $d_{2t+1,0} = {\bf z}(a-t-1,1,c+t)+3c$ together with $X_{2t+\varepsilon,0}  \in \pm q^{d_{2t+\varepsilon,0}} (1+q {A}_0)$ ($\varepsilon=0,1$). The induction completes.
\smallskip
	
{\em Step 3}. Consider $d_{s,r}$ for $r\le c$. We use induction on $r$. It is true for $r=0$ by {\em Step 2}.	
By \eqref{eq: app e even formula}, if $r+1 \le c$, we have
\begin{equation*}
	d_{2t,r+1} = 2a+2c-6r-4 + {\rm min} (d_1, d_2),
\end{equation*}
when $d_1 \ne d_2$, where $d_1=d_{2t,r}-2r$ and $d_2=d_{2t-1,r}-2t+1$.
We have $d_1-d_2 = -4(c-r) <0$ by induction hypothesis, and hence $d_{2t,r+1}=2a+2c-8r-4+d_{2t,r}$.
This implies $d_{2t,r+1}={\bf z}(a-t-r-1,2r+2,c+t-r-1)+(c-r-1)$ and $X_{2t,r+1} \in \pm q^{d_{2t,r+1}} (1+q {A}_0)$.

By \eqref{eq: app e odd formula}, if $r+1 \le c$, we have
\begin{equation*}
	d_{2t-1,r+1} = 2a+2c-6r-7 + {\rm min} (d_1, d_2)
\end{equation*}
when $d_1 \ne d_2$, where
$d_1=d_{2t-2,r+1}-2t+2$ and $d_2=d_{2t-1,r}-2r-1$.
Then we have $d_1-d_2 = -2 <0$ by induction hypothesis, and hence $d_{2t-1,r+1} = 2a+2c-6r-2t-5+d_{2t-2,r+1}$.
This implies $d_{2t-1,r+1}={\bf z}(a-t-r-1,2r+3,c+t-r-2)+3(c-r-1)$ and $X_{2t-1,r+1} \in \pm q^{d_{2t-1,r+1}} (1+q {A}_0)$.
The induction completes.
\smallskip
	
{\em Step 4}. We have
\begin{equation*}
\begin{split}
 d_{2t,c+t}&=2a-4c-4t+1+d_{2t-1,c+t-1},\\
 d_{2t+1,c+t}&=2a-4c-4t-1+d_{2t,c+t},
\end{split}
\end{equation*}
by \eqref{eq: app f odd formula} and \eqref{eq: app f even formula}, respectively.
Since $d_{0,c}={\bf z}(a-c,2c,0)$, we have
\begin{equation*}
	\begin{split}
	d_{2t,c+t} & ={\bf z}(a-c-2t,2c+2t,0)+2t^2-t, \\
	d_{2t-1,c+t-1} & = {\bf z}(a-c-2t+1,2c+2t-1,0)+2(t-1)^2+(t-1),
	\end{split}
\end{equation*}
and $X_{2t-\varepsilon,c+t-\varepsilon}  \in \pm q^{d_{2t-\varepsilon,c+t-\varepsilon}} (1+q {A}_0)$ ($\varepsilon=0,1$) by induction hypothesis.
\smallskip
	
{\em Step 5}. Consider $d_{s,r}$ for $r\ge c$.
By {\em Step 4}, \eqref{eq:formula for d;d} holds for $r=c+t$ if $s=2t$ and $r=c+t-1$ if $s=2t-1$. 
We use induction on $s$ and also induction on $r$ in decreasing order for each $s$. 
Suppose that \eqref{eq:formula for d;d} holds for $s=2t-2$.
 
First, we have by \eqref{eq: app e odd formula}
\begin{equation*}
	d_{2t-1,r-1} = 2r-1 + {\rm min} (d_1, d_2) \quad \text{if $r-1 \ge c$},
\end{equation*}
when $d_1 \ne d_2$, where
$d_1=d_{2t-1,r}-2a-2c+6r+1$ and $d_2=d_{2t-2,r}-2t+2$.
We have $d_1-d_2=4(r-c)>0$ by induction hypothesis, and hence $d_{2t-1,r-1}=-2t+2r+1+d_{2t-2,r}$.
This implies $d_{2t-1,r-1}={\bf z}(a-t-r+1,2r-1,c+t-r)+2(r-c-1)^2+(r-c-1)$ and $X_{2t-1,r-1} \in \pm q^{d_{2t-1,r-1}} (1+q {A}_0)$.
	
Next, we have by \eqref{eq: app e even formula}
\begin{equation*}
	d_{2t,r-1} = 2r-2 + {\rm min} (d_1, d_2) \quad \text{if $r-1 \ge c$},
\end{equation*}
when $d_1 \ne d_2$, where
$d_1=d_{2t,r}-2a-2c+6r-2$ and $d_2=d_{2t-1,r-1}-2t+1$.
We have $d_1-d_2=5(r-c-1)+2>0$ by induction hypothesis, and hence $d_{2t,r-1}=d_{2t-1,r-1}+2r-2t-1$.
This implies $d_{2t,r-1}={\bf z}(a-t-r+1,2r-2,c+t-r+1)+2(r-c-1)^2-(r-c-1)$ and $X_{2t,r-1}  \in \pm q^{d_{2t,r-1}} (1+q {A}_0)$.
The induction completes.\smallskip

By {\em Steps 1--5} and Lemma \ref{lem:symmetry of recursive formulas;d}, \eqref{eq:formula for d;d} holds for all $s,r$.
\qed

\end{document}